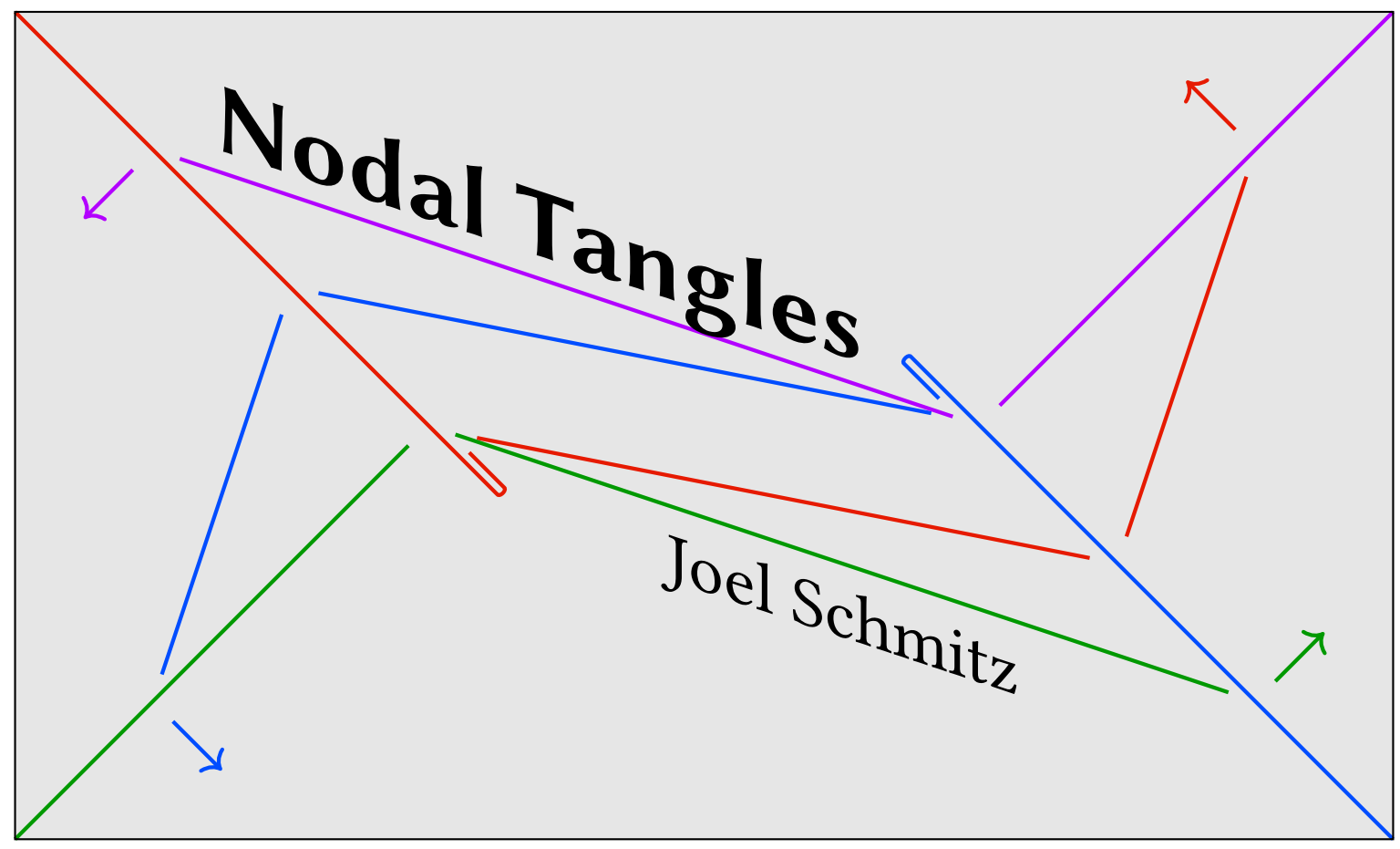

**Abstract**

We study piecewise linear knot diagrams in the base of almost toric fibrations of symplectic four-manifolds. These diagrams translate to deformations of the almost toric fibration. We give several applications to symplectic topology, among them a proof of a conjecture by Symington, simpler counterexamples to Lagrangian Poincaré recurrence in dimension four, the calculation of the displacement energy for many fibres of toric moment maps, and an elementary recipe for building and distinguishing Lagrangian torus knots.

# 1 Introduction

Almost toric fibrations, that were introduced by Symington in [Sym03], have played a central role in symplectic topology over the last years. For instance, they have been used to construct Lagrangian torus knots [Via16; Via17; BHS24], to simplify known and establish new staircases for ellipsoid embeddings [CV22; Mag24; MS24], and to show that some rigidity results in algebraic geometry still persist in symplectic geometry [ES18]. When a symplectic 4-manifold $(X, \omega)$ admits an almost toric fibration, then many properties of $(X, \omega)$ can be read from the base $B$ of the fibration, a two-dimensional integral affine surface with *nodes*. In this paper, we focus on how modifications of the nodal integral affine base $B$ of the almost toric fibration translate to symplectic properties of the total space. Unlike the *mutations* used by Vianna in [Via16; Via17] for the same goal, where the shape of the base diagram changes, we will keep the shape of the base diagram fixed. This will improve or simplify some of the previous applications of almost toric fibrations and lead to new applications. A central concept of the paper are **nodal tangles** (Definition 2.21), which can be thought of as a sequence of *nodal slides* in the base $B$ as defined in [Sym03], or equivalently as a certain smooth deformation of the nodal integral affine structure of $B$. In order to obtain applications of this combinatorial manipulation of $B$ to symplectic geometry, we use three “translation theorems”:

- Theorem 2.27 associates to a nodal integral affine surface a unique almost toric fibration.
- Theorem 2.30 describes how a nodal tangle gives rise to a path of almost toric fibrations.
- Corollary 2.35 describes how symplectic invariants of Lagrangian tori change along nodal tangles.

These three translation theorems are reformulations of results from [Sym03; Eva23; BHS24].

In Section 2 we define nodal integral affine surfaces and nodal tangles, give some tools which allow us to manipulate nodal integral affine surfaces, and formulate the three translation theorems above precisely.

In Section 3 we use the setup of Section 2 and results from elementary tropical geometry in [MS23] to show that a certain class of nodal integral affine surfaces, which includes Delzant polygons, can be deformed by a nodal tangle into a *canonical form*. We use this canonical form to prove three results, summarized in Sections 1.1 to 1.3.

In Section 4 we use nodal tangles to provide an elementary recipe for building and distinguishing Lagrangian torus knots using nodal tangles, summarized in Section 1.4. The same combinatorics lead to an interesting example regarding Lagrangian pinwheels, connecting to the recent work [UZ25] describing which Wall singularities can arise in degenerations of del Pezzo surfaces, described in Section 1.5.

## 1.1 Symington's conjecture

Conjecture 6.8 in [Sym03] states that given two toric moment maps on a four dimensional symplectic manifold, we can always find a family of almost toric fibrations interpolating between them. We prove this conjecture for closed toric symplectic four-manifolds:

**Theorem A** (Corollary 3.26)**.** *Let $\mu_i : X \to \Delta_i$ with $i \in \{0, 1\}$ be two toric moment maps on a closed four dimensional symplectic manifold. Then $\Delta_0$ and $\Delta_1$ are connected by a nodal tangle.*

*In particular, there is a continuous path of almost toric fibrations $\pi_t : X \to B_t$ with $B_0 = \Delta_0$ and $B_1 = \Delta_1$.*

A toric moment map on a symplectic manifold $(X, \omega)$ describes a strong set of symmetries (in the sense of Noether) of the space $(X, \omega)$. Theorem A says that any two such sets of symmetries must be related by a nodal tangle.

For the proof, we first deform both $\Delta_0$ and $\Delta_1$ by a nodal tangle into their canonical form, and then show that the canonical form only depends on the manifold $X$, which means that $\Delta_0$ and $\Delta_1$ have the same canonical form. The path of almost toric fibrations is then obtained using Theorem 2.30, the translation theorem allowing us to "lift" a nodal tangle to a one-parameter family of almost toric fibrations.

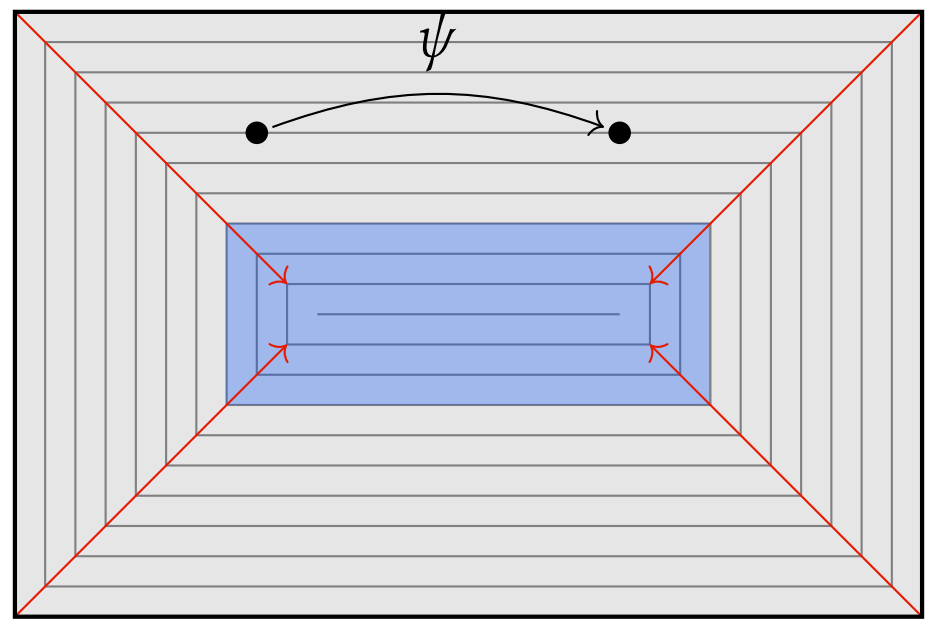

Figure 1: The behaviour of $\psi$ on the base space $\square$ for $S^2 \times S^2$. Nodes are represented by arrowheads, the emerging lines are branch cuts. Level sets of $\mathcal{F}$ are in grey. The shaded region in the middle is $\mathcal{F}^{-1}((M-\varepsilon, M])$.

## 1.2 Lagrangian Poincaré non-recurrence

In [Sch24] we constructed a counterexample to Lagrangian Poincaré recurrence in dimension four using almost toric fibrations. That is, we constructed a Hamiltonian diffeomorphism $\psi$ and a Lagrangian $L$ such that $\psi^n(L) \cap L = \emptyset$ for all $n \in \mathbb{N}$ in a symplectic four manifold. Counterexamples in dimension six and higher had been found shortly before in [BS24]. See [Sch24] for a discussion of related results.

In Section 3.3 we use nodal tangles and the canonical form developed in Section 3.1 to put the construction of [Sch24] in a more general context and also provide a slightly different construction to give more counterexamples. Specifically, we construct counterexamples in all closed non-monotone toric symplectic four manifolds:

**Theorem B.** *For a closed non-monotone toric symplectic four manifold $(X, \omega)$, let $\pi : (X, \omega) \to B$ be the almost toric fibration over the canonical form $B$ of $X$, and let $\mathcal{F} : B \to \mathbb{R}$ be the* height function, *and $M = \max \mathcal{F}$. For any $\varepsilon > 0$ there exists a Hamiltonian diffeomorphism $\psi$ such that for almost all*[1] *points $x$ in the set $\mathcal{F}^{-1}([0, M-\varepsilon])$ the pair $(\pi^{-1}(x), \psi)$ provides a counterexample to Lagrangian Poincaré recurrence.*

The height function $\mathcal{F}$ measures the distance to the boundary $B$, such that $\mathcal{F}^{-1}([0, M-\varepsilon])$ is a set of measure arbitrarily close to the area of $B$. See also the following example.

**Example 1.1.** We quickly describe the effect of the constructed Hamiltonian diffeomorphism $\psi$ in the case of non-monotone $S^2 \times S^2$, illustrated in Figure 1. Take $\omega_{a,b}$ to be the product symplectic form on $S^2 \times S^2$ giving areas $a > b$ to the factors. This admits a unique split toric moment map $\mu : S^2 \times S^2 \to \square$, where $\square \subset \mathbb{R}^2$ is an axis-aligned rectangle of side lengths $a, b$. After a nodal trade at each corner we have a base diagram as in Figure 1, giving a new fibration $\pi : S^2 \times S^2 \to \square$. The height function $\mathscr{F}$ gives the integral affine distance to the boundary of $\square$. Its level sets are also highlighted in Figure 1; its maximum is $M = \frac{b}{2}$. Let $\varepsilon > 0$. We can assume that all nodes are contained in the set $\mathscr{F}^{-1}((M-\varepsilon, M])$, highlighted in blue in Figure 1. Then by Theorem B we may take a Hamiltonian diffeomorphism $\psi$ which is fibred with respect to $\pi$ on the set $\mathscr{F}^{-1}([0, M-\varepsilon])$ and acts as follows on the fibres of $\pi$: Take $h \in [0, M-\varepsilon]$. Note that the integral affine structure on $\square$ induced by $\pi$ (see Section 2) makes the level set $\mathscr{F}^{-1}(h)$ into a straight line (Definition 2.15). The Hamiltonian diffeomorphism $\psi$ sends the fibre over a point $x \in \mathscr{F}^{-1}(h)$ to the fibre over the point $x$ translated by a distance of $2(M-h)$ along $\mathscr{F}^{-1}(h)$ in the clockwise direction. The length of $\mathscr{F}^{-1}(h)$ is $2(a-b) + 8(M-h)$, so if $\frac{M-h}{a-b}$ is irrational, then any fibre over a point $x \in \mathscr{F}^{-1}(h)$ and $\psi$ give a counterexample to Lagrangian Poincaré recurrence.

[1] With respect to the Lebesgue measure on $B$.

Thus we can (partly) classify points $x \in \square$ depending on whether Lagrangian Poincaré recurrence holds for the fibre $\pi^{-1}(x)$:

- If $\frac{M-h}{a-b}$ is irrational, then Lagrangian Poincaré recurrence for fibres over $\mathscr{F}^{-1}(h)$ does not hold.
- If $x$ is in the interior of the line segment $\mathscr{F}^{-1}(M)$, Lagrangian Poincaré recurrence does hold.
- If $x$ is in a level set $\mathscr{F}^{-1}(h)$ with $\frac{M-h}{a-b}$ rational, or $x$ is one of the end points of the line segment $\mathscr{F}^{-1}(M)$, it is unknown whether Lagrangian Poincaré recurrence holds.

The second point is proven in [PS23, Theorem C]. For convenience see also the reformulation [BK25, Theorem 1.15].

For further situations in which a stronger version of Lagrangian Poincaré recurrence holds, see [GG18; JS24].

*Remark* 1.2. For monotone toric symplectic manifolds a similar construction as for Theorem B can be carried out. However in this case it only yields finite orbits of Lagrangian fibres.

## 1.3 Displacement energy of toric fibres

Our third result in Section 3 concerns the displacement energy of fibres of toric moment maps. The displacement energy of a set $A \subset (X, \omega)$ is given by

$$e(A, X) := \inf\left\{ \|H\| \;\middle|\; \begin{array}{l} H \text{ is a time-dependent compactly supported Hamiltonian} \\ \text{on } X \text{ such that its time-1 flow } \varphi_H^1 \text{ satisfies } \varphi_H^1(A) \cap A = \emptyset \end{array} \right\},$$

where

$$\|H\| = \int_0^1 \left( \sup_{x \in X}\{H_t(x)\} - \inf_{x \in X}\{H_t(x)\} \right) dt$$

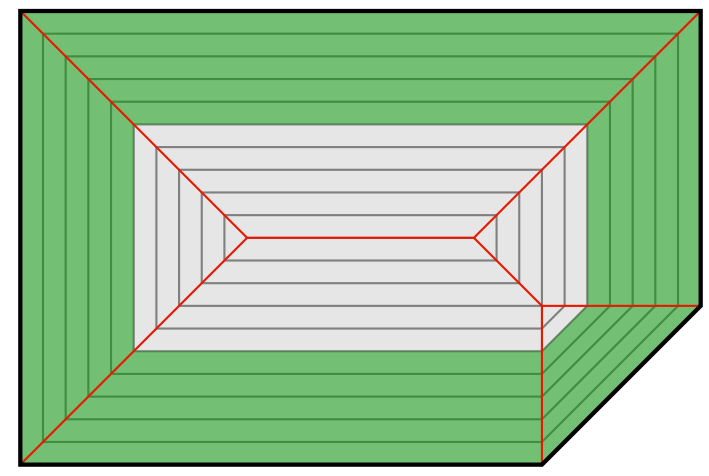

Figure 2: The green area marks the fibres displaceable via Theorem C, the red curve marks the caustic, the excluded one-dimensional subset.

is the Hofer norm. The displacement energy gives a symplectic invariant of the set $A$. In particular, it gives an invariant of Lagrangian tori.

Theorem C calculates the displacement energy of certain fibres of toric symplectic manifolds $\mu : X \to \Delta$, where $\Delta$ is a *Delzant polygonal domain*. This means that $\Delta$ is locally given by the intersection of half-planes

$$\langle \lambda_i, x \rangle + c_i \geq 0,$$

where $\lambda_i \in \mathbb{Z}^2$ are primitive vectors. (See Definitions 3.1 and 3.5 for details.) In particular this includes the case where $X$ is closed and $\Delta$ is a Delzant polygon. Define

$$\mathcal{F}_\Delta : \Delta \to \mathbb{R}_{\geq 0}, \quad \mathcal{F}_\Delta(x) := \min\{\langle \lambda_i, x \rangle + c_i\},$$

where we take the minimum over all the half-planes defining an edge of $\Delta$. Hence $\partial\Delta = \mathcal{F}_\Delta^{-1}(0)$ and $\mathcal{F}_\Delta$ gives a notion of distance to the boundary of $\Delta$.[2]

The canonical form constructed in Section 3.1 gives an almost toric fibration on $X$, and we can use probes in its nodal integral affine base to calculate the displacement energy for toric fibres:

**Theorem C** (Corollary 3.32)**.** *Let $\mu : X \to \Delta$ be a toric symplectic manifold over a Delzant polygonal domain $\Delta$. Except on a one-dimensional subset of $\Delta$, called the* caustic, *if $\mathcal{F}_\Delta(x) < \frac{1}{2} \sup \mathcal{F}_\Delta$ the displacement energy of a fibre is given by*

$$e(\mu^{-1}(x)) = \mathcal{F}_\Delta(x).$$

See Figure 2 for an illustration.

The exclusion of the one-dimensional caustic in Theorem C is often irrelevant for distinguishing Lagrangian tori, as it often suffices to know $e$ on an open dense subset of fibres of $\mu$.

For any $x \in \Delta$, the lower bound $e(\mu^{-1}(x)) \geq \mathcal{F}_\Delta(x)$ is shown to always hold in [Bre23, Proposition 3.2]. The upper bound $e(\mu^{-1}(x)) \leq \mathcal{F}_\Delta(x)$ is not always known to hold. Outside the caustic, it can sometimes be proven in explicit examples using probes, see [McD11], or extended probes, see [ABM14], where we can use [Bre23, Proposition 3.4] to get the displacement energy.

[2] Sometimes $\mathcal{F}_\Delta$ is described as *integral affine distance to* $\partial\Delta$. Since the integral affine length (Definition 2.17) does not induce a metric, this is ambiguous.

*Remark* 1.3. If $\sup \mathscr{F}_\Delta = \infty$ in Theorem C, then every fibre outside the caustic is displaceable by a probe, providing an improvement on the method of extended probes used in [ABM14, Proposition 4.4.4]. For example all fibres in the white region in [ABM14, Figure 4.5.3] can be displaced by a probe using the canonical form as a base diagram.

*Remark* 1.4. For all examples of Delzant polygons I examined, the method of proof of Theorem C can be modified to give $e(\mu^{-1}(x)) = \mathscr{F}_\Delta(x)$ for all fibres outside the caustic, even without the condition $\mathscr{F}_\Delta(x) < \frac{1}{2} \sup \mathscr{F}_\Delta$, see Remark 3.33.

## 1.4 Lagrangian knots

The study of Lagrangian knots in a symplectic manifold $(X, \omega)$ involves classifying embedded Lagrangians in $X$ up to either Lagrangian isotopy, symplectomorphism or Hamiltonian diffeomorphism, see [PS24] for a survey. In this paper we are concerned with the classification up to symplectomorphism. When constructing Lagrangian knots, we often want to rule out that two Lagrangians $L_0, L_1$ are not related by a symplectomorphism for “obvious” reasons, such as being topologically different or having different area classes $[\omega] : H_2(X, L_i) \to \mathbb{R}$. The following definition gives a criterion for when two Lagrangians are sufficiently similar that telling them apart becomes interesting.

**Definition 1.5.** Two Lagrangian embeddings $\varphi_0, \varphi_1 : L \to X$ are **almost Hamiltonian isotopic** if for any neighbourhood $U$ of $0 \in H^1(L; \mathbb{R})$ there is a Lagrangian isotopy $\varphi : L \times [0,1] \to X$ between them such that

$$\forall t \in [0,1], \quad \mathrm{Flux}(\varphi|_{L\times[0,t]}) \in U,$$

where Flux is the *Lagrangian flux map.*

The Lagrangian flux map is defined by measuring the symplectic area swept out by a representative of an element of $H_1(L)$ under the Lagrangian isotopy, see [BHS24, Definition 3.2] for details. If $\mathrm{Flux}(\varphi|_{L\times[0,t]}) = 0$ for all $t$, then the Lagrangian isotopy can be realized by an ambient Hamiltonian isotopy ([BHS24, Lemma 3.3]), and the Lagrangians $\varphi_0(L)$ and $\varphi_1(L)$ are Hamiltonian isotopic. For example Vianna’s tori in [Via16; Via17] are all almost Hamiltonian isotopic.

In Section 4 we describe how nodal integral affine surfaces can be used to construct almost Hamiltonian isotopic Lagrangian torus knots: Almost toric fibrations have been used to construct knotted Lagrangian tori in [Via16; Via17; Bre23; Bre25; BHS24] as follows. This method is due to Vianna in [Via16]. Let $\pi_0 : X \to B_0$ be an almost toric fibration, let $x_0 \in B_0$ be a regular value and let $\mathfrak{a}_1, \ldots, \mathfrak{a}_n \in B_0$ be nodes whose eigenlines run through $x_0$. Sliding the nodes $\mathfrak{a}_1, \ldots, \mathfrak{a}_n$ back and forth over $x_0$ we get a nodal tangle $B_0 \times \mathbb{R}_{\geq 0}$, that is, a one-parameter family of nodal integral affine surfaces $B_t$ with $t \geq 0$, and a corresponding family of almost toric fibrations $\pi_t : X \to B_t$, see Figure 3. By sliding the nodes $\mathfrak{a}_1, \dots, \mathfrak{a}_n$ back and forth over $x_0$ we modify the integral affine structure of $B_0$ near $x_0$, as well as the fibration $\pi_t$. The fibre $\pi_0^{-1}(x_0)$ is a Lagrangian torus. By sliding a node over $x_0$, we modify the fibre $\pi_0^{-1}(x_0)$ by a *la-disc surgery* as in [STW16]. Almost toric fibrations now provide a convenient way to distinguish the modified fibres from the original ones: Let $\mathscr{L}$ denote the space of Lagrangians in $X$ and let $\mathscr{I} : \mathscr{L} \to A$ be an invariant of Lagrangians up to symplectomorphism with values in some set $A$. The germ $[\mathscr{I}]_{x_0}$ given by evaluating $\mathscr{I}$ on fibres near $x_0 \in B_0$ only depends on the fibre $\pi_0^{-1}(x)$ (Lemma 2.37). After having slid $n$ nodes back and forth over $x_0$, the germ $[\mathscr{I}]_{x_0}$ is modified by the *transition map* $\tau_0^n$, a piecewise affine map which describes how the integral affine structure changes near $x_0$ (Corollary 2.35). Then if $[\mathscr{I}]_{x_0}$ and $[\mathscr{I}]_{x_0} \circ \tau_0^n$ are not related by an integral linear transformation, the fibres $\pi_0^{-1}(x_0)$ and $\pi_n^{-1}(\tau_0^n(x_0))$ cannot be related by a symplectomorphism. So when we want to create Lagrangian knots, we can try to create "complicated" transition maps $\tau_0^n$ in order to change the invariant germ $[\mathscr{I}]_{x_0}$ enough to be distinguishable.

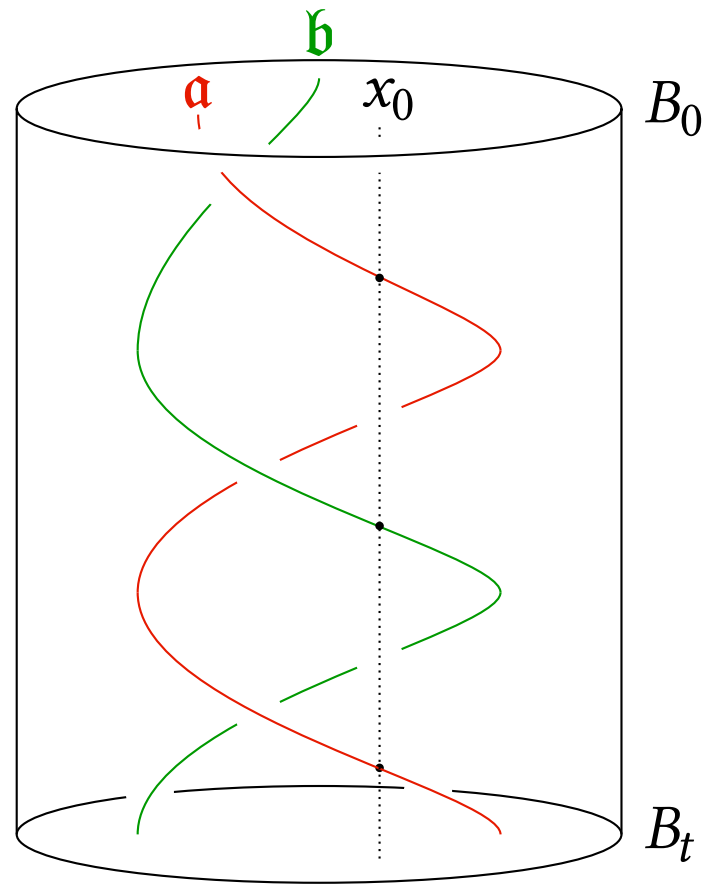

Figure 3: A topological picture of a nodal tangle used to produce Lagrangian knots. Here we slide two nodes $\mathfrak{a}, \mathfrak{b}$ alternatingly over $x_0$.

Using an invariant germ to distinguish Lagrangian tori has two advantages compared to using a holomorphic disc count as done in [Via16; Via17]: First, we can also treat non-monotone tori, that is, the point $x_0$ does not need to be the centre of a monotone moment polytope $\Delta$, since the displacement energy germ can be calculated also for non-monotone tori, whereas it is difficult to calculate the holomorphic disc count of non-monotone Lagrangian tori. Second, the combinatorics of constructing the Lagrangian knot are directly linked to the invariant germ via the integral affine structure of the base. This allows us to construct and distinguish Lagrangian knots over many more points $x$ in the base. See Example 4.7 and Sec-tion 5.4 for a discussion on how many families of knotted tori can be obtained using this method.

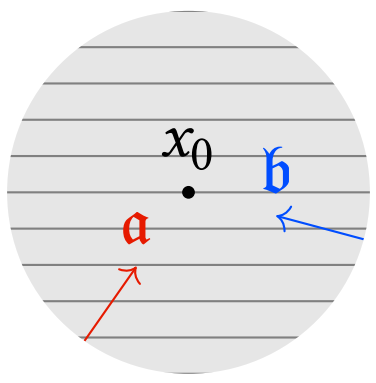

Figure 4: Two nodes $\mathfrak{a}, \mathfrak{b}$, marked by arrowheads, whose eigenlines pass through $x_0$. The eigenlines are parallel to the primitive integer vectors $v_{\mathfrak{a}}, v_{\mathfrak{b}}$ (the direction of the arrow tips). The level sets of the invariant $\mathcal{I}$ are in grey.

Theorem 4.6 describes a particular but still quite general situation in which we get infinitely many different Lagrangian tori over a point $x_0$, illustrated in Figure 4.

**Theorem D** (Theorem 4.6)**.** *Suppose that $\mathfrak{a}, \mathfrak{b}$ are two nodes of multiplicity $k_{\mathfrak{a}}, k_{\mathfrak{b}}$ with eigenlines running through $x_0$ in directions $v_{\mathfrak{a}}, v_{\mathfrak{b}}$, and that*

1. *$k_{\mathfrak{a}} k_{\mathfrak{b}} \det(v_{\mathfrak{a}}, v_{\mathfrak{b}})^2 \geq 4$,*
2. *The germ $[\mathcal{I}]_{x_0} : B_0 \dashrightarrow \mathbb{R}$ is affine and non-constant on an open dense set around $x_0$.*

*Then the invariant germs $[\mathcal{I}]_{x_0} \circ \tau_n^0$ with $n \in \mathbb{Z}_{\geq 0}$ are pairwise different, hence the Lagrangian tori $\{\pi_n^{-1}(\tau_0^n(x_0))\}_n$ are pairwise not related by symplectomorphism.*

*Remark* 1.6. The assumption that $[\mathcal{I}]_{x_0}$ is affine near $x_0$ is for example fulfilled for the displacement energy germ of all points in a Delzant polytope at which Theorem C holds.

Examples of infinite families of non-monotone tori have previously been constructed in [Bre25] in dimension six and in [BHS24] in dimension four, both also using the displacement energy germ as distinguishing invariant.

*Remark* 1.7. Another candidate for a good choice of $\mathcal{I}$ is the $\Psi$ invariant described in [STV24]. In this paper the authors calculate $\Psi$ for toric fibres of $\mu : X \to \Delta$, where $X$ is a Fano variety of any dimension, in which case they find $\Psi(\mu^{-1}(x)) = \mathcal{F}_\Delta(x)$, where $\mathcal{F}_\Delta$ is defined analogously for the higher dimensional polytope $\Delta$. Note however that $X$ being Fano means that $X$ is a monotone symplectic manifold; in particular if $X$ is four dimensional, then $X$ is a del Pezzo surface and the caustic of $\Delta$ is star shaped. Thus only the monotone fibre can have eigenlines of two nodes passing through its base point, and we are in the situation described in [Via17].

## 1.5 Lagrangian pinwheels

Another application of the above combinatorics used to construct Lagrangian knots is the construction of *Lagrangian pinwheels*. For coprime $p, q \in \mathbb{Z}$, the Lagrangian pinwheel $L_{p,q}$ is an embedded Lagrangian CW-complex which arises in algebraic geometry as the vanishing cycle of a $\mathbb{Q}$-Gorenstein smoothing of a Wahl singularity. [Eva23, Remark 7.9] Alternatively, they can be seen as *visible Lagrangians* lying over an eigenray of a node $\mathfrak{a}$ connecting it to $\partial B$. In the recent work [UZ25], the authors use methods from algebraic geometry to show which Wahl singularities can arise from degenerations of del Pezzo surfaces. In particular, [UZ25, Theorem 1.10] implies that any monotone symplectic $\mathbb{C}P^2\#n\overline{\mathbb{C}P^2}$ with $5 \leq n \leq 8$ admits an embedding of $L_{p,q}$ for any coprime $p, q \in \mathbb{Z}$: A Wahl singularity is a cyclic quotient T-singularity $\frac{1}{dp^2}(1, dpq - 1)$ with $d = 1$, which can be smoothed as in [Eva23, Section 7.4] to give a Lagrangian pinwheel. In Example 4.9 we give an alternative construction of these pinwheels using our methods.

### 1.6 Acknowledgements

I would like to thank Joé Brendel, Johannes Hauber and Felix Schlenk for listening to my incoherent rambles, letting me form them into something serviceable.

## 2 Nodal integral affine surfaces & almost toric fibrations

In Section 2.1 we define nodal integral affine surfaces independently from almost toric fibrations and setup some useful language to talk about them. In Section 2.2 we introduce nodal tangles and give a worked out examples how they can be manipulated. In Section 2.3 we relate nodal integral affine surfaces and nodal tangles to almost toric fibrations using [Sym03; Eva23], giving our first two translation theorems 2.27 and 2.30. In Section 2.4 we describe how invariant germs change under nodal tangles, giving our third translation theorem 2.35.

### 2.1 Nodal integral affine surfaces

Denote by $\mathrm{Aff}_n(\mathbb{Z}) = \mathrm{GL}_n(\mathbb{Z}) \rtimes \mathbb{R}^n$ the group of integral affine transformations of $\mathbb{R}^n$.

**Definition 2.1.** An **integral affine manifold** $B$ is a smooth manifold with corners with an *integral affine structure*, that is an equivalence class of atlases $\{\varphi_i : B \supset U_i \to \mathbb{R}^n_{\geq 0}\}_{i \in I}$ such that the transition functions $\varphi_i \circ \varphi_j^{-1}$ are elements of $\mathrm{Aff}_n(\mathbb{Z})$.

The **integral lattice bundle** of $B$ is the lattice bundle $\Lambda B \subset TB$ given by the pullback of $\mathbb{Z}^n \subset \mathbb{R}^n$ under integral affine charts. The **dual integral lattice bundle** $\Lambda^* B \subset T^* B$ is the lattice bundle dual to $\Lambda B$.

**Definition 2.2.** An **integral linear map** $T : \mathbb{R}^n \to \mathbb{R}^m$ is a linear map given as the linear extension of a linear map $\mathbb{Z}^n \to \mathbb{Z}^m$.

A **piecewise integral linear map** $R : \mathbb{R}^n \to \mathbb{R}^m$ is a continuous map such that there exists a complete fan $\Delta$ on $\mathbb{R}^n$ such that for any cone $c$ of $\Delta$ the restriction $R|_c$ is integral linear.

**Definition 2.3.** A **(piecewise) integral affine map** $\psi : A \to B$ between integral affine manifolds is a locally (piecewise) integral affine map. That is for any point $x \in A$ there exist integral affine charts $\varphi_A, \varphi_B$ on $A, B$ around $x, \psi(x)$ respectively, such that $\varphi_B \circ \psi \circ \varphi_A^{-1}$ is (piecewise) integral linear.

Note that in the literature a piecewise integral affine map as in Definition 2.3 would be usually called a *piecewise integral linear map*. This makes it impossible to distinguish piecewise integral linear maps and piecewise integral affine maps $\mathbb{R}^n \to \mathbb{R}^m$, so we use the non-standard but more sensible *piecewise affine* instead.

**Definition 2.4.** The **cut complex** of a piecewise integral affine map $\tau : A \to B$ is the set $\Delta$ of points where $\tau$ is not locally integral linear.

By Definition 2.2, the cut complex locally looks like the co-dimension one faces of a fan.

**Definition 2.5.** A **pre-nodal integral affine manifold** of dimension $n$ is a pair $(B, \mathfrak{N})$ where

- $B$ is a smooth $n$-dimensional manifold with corners,
- $\mathfrak{N} \subset B$ is a closed stratified subspace with smooth embedded strata and top dimension $n-2$,
- $B \setminus \mathfrak{N}$ is an integral affine manifold with corners.

The points of the set $\mathfrak{N}$ are called **nodes**.

We will be mostly interested in the case $n = 2$, where $\mathfrak{N}$ is a closed set of isolated points in $B$. The case $n = 3$ will be relevant for Definition 2.21, where we will be mostly interested in the case where $\mathfrak{N}$ is a closed 1-dimensional embedded submanifold.

Let $(B, \mathfrak{N})$ be a pre-nodal integral affine surface. An integral affine chart on $B \setminus \mathfrak{N}$ cannot describe the neighbourhood of a node. The following definition tries to address this drawback. The concept is similar to that of a base diagram used e.g. in [Sym03; Eva23].

**Definition 2.6.** Let $U \subset B$ be open. A **nodal chart** on $U$ is a homeomorphism onto its image $\varphi : U \to \mathbb{R}^n$ that is piecewise integral affine on $U \setminus \mathfrak{N}$.

An open $U \subset B$ might not admit a nodal chart.

**Definition 2.7.** For any nodal chart $\varphi : U \to \mathbb{R}^n$ and any point $x \in U \setminus \mathfrak{N}$, after a translation taking $\varphi(x)$ to 0, there is an invertible piecewise integral linear map $R_x : \mathbb{R}^n \to \mathbb{R}^n$ such that $R_x \circ \varphi$ is an integral affine chart on a neighbourhood of $x$. We call $R_x$ the **rectifying map of** $\varphi$ **at** $x$.

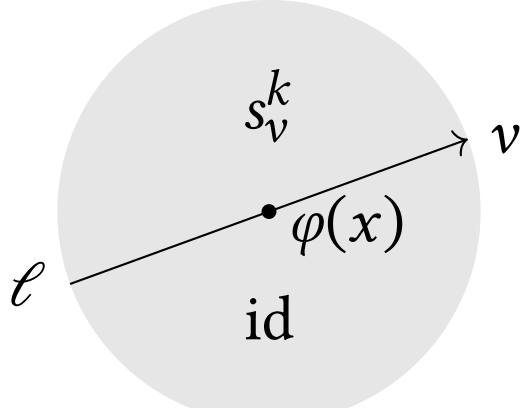

Figure 5: The rectifying map in the case of Remark 2.8.

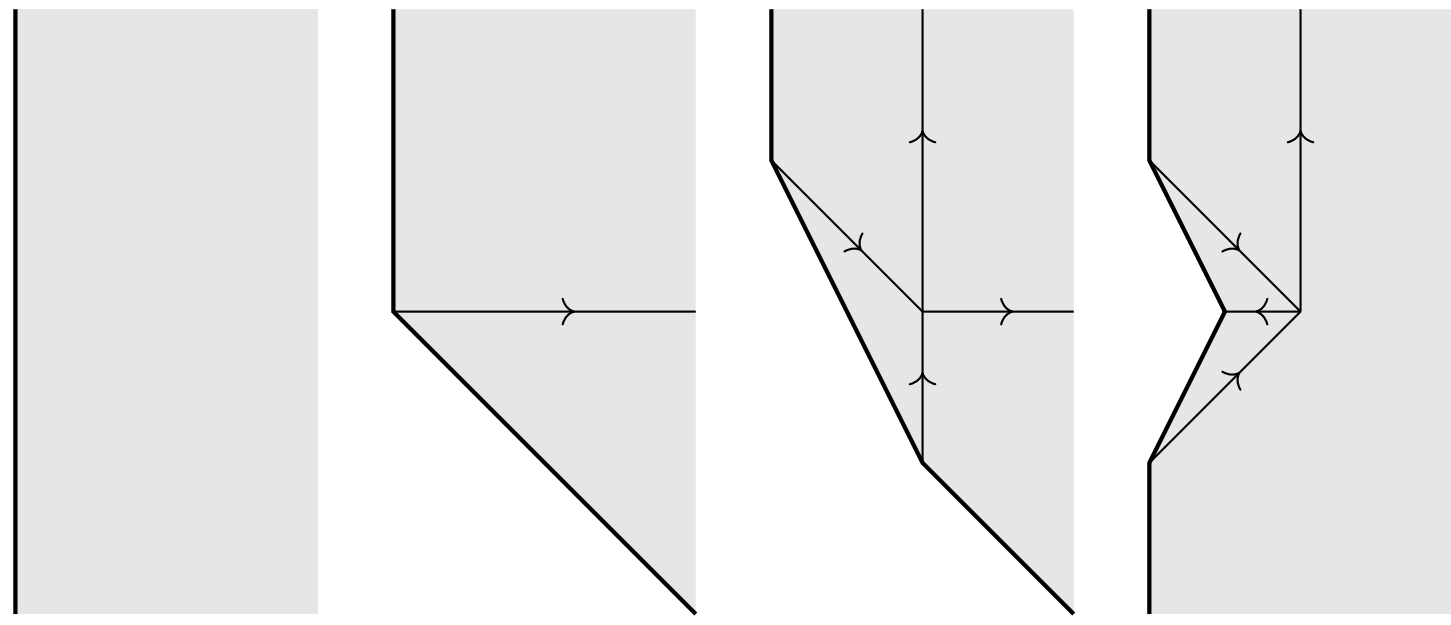

Figure 6: Four different nodal charts diagrams of nodal charts of $\mathbb{R}_{\geq 0} \times \mathbb{R}$ with the standard integral affine structure. On the left is the integral affine chart $\mathrm{id}_{\mathbb{R}_{\geq 0} \times \mathbb{R}}$. Each adjacent pair is related by a half-shear.

The rectifying map is unique up to composing with an integral linear map. For $x$ in the integral affine locus of $\varphi$, we may choose $R_x = \mathrm{id}$.

For the rest of this chapter we will look at the case $n = 2$.

Let $\varphi : U \to V$ be a nodal chart. The cut complex of $\varphi^{-1} : V \to U$ form an embedded graph $G$ in $V \subset \mathbb{R}^2$. The edges of $G$ have rational slope (see Remark 2.8 below). We call $G$ the **cut graph of** $\varphi^{-1}$ and the edges of $G$ **cuts of** $\varphi^{-1}$.[3]

*Remark* 2.8. Suppose $x \in B \setminus \mathfrak{N}$ is a point such that $\varphi(x)$ is in the interior of a cut $\ell$, as in Figure 5. Then its rectifying map $R_x$ is given by two integral affine maps that agree on $\ell$. W.l.o.g. we may take one of them to be the identity. Then the other one must fix $\ell$, and thus is the integral shear map

$$s_v^k(y) = y - kv \det(v, y) ,$$

where $v$ is a primitive vector along $\ell$, and $k$ a positive integer. Requiring $k$ to be positive induces an orientation on $\ell$, given by $v$. So $R_x$ must be equal to the **half-shear**

$$h_v^k(y) = \begin{cases} s_v^k(y) & \text{if } \det(v, y) \geq 0 \\ y & \text{if } \det(v, y) \leq 0 \end{cases} . \tag{1}$$

[3]The name comes from the similarity to *branch cuts* of almost toric base diagrams ([Eva23, Definition 8.3]). Since we're not cutting branches to get a fundamental domain of $B \setminus \mathfrak{N}$ in its universal cover, we just call them cuts, as we're still cutting the integral affine structure.

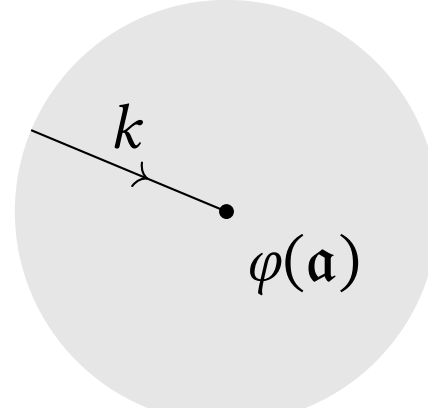

Figure 7: A nodal chart diagram of a neighbourhood of a node $\mathfrak{a}$ of multiplicity $k$.

This makes $G$ into a directed weighted graph, with direction on $\ell$ given by $v$ and weight given by $k$. The pair $(\operatorname{im}\varphi, G)$ lets one reconstruct the integral affine structure on $U \setminus \mathfrak{N}$. See Figure 6 for an example. We call the pair $(\operatorname{im}\varphi, G)$ a **nodal chart diagram**.

**Definition 2.9.** A nodal chart $\varphi : U \to \mathbb{R}^2$ is called **simple** if all vertices of the cut graph $G$ are 1-valent, i.e. $G$ is a disjoint union of line segments.

**Definition 2.10.** A **nodal integral affine surface** is a pre-nodal integral affine surface $(B, \mathfrak{N})$ such that every node $\mathfrak{a} \in \mathfrak{N}$ admits a simple nodal chart $\varphi : U_{\mathfrak{a}} \to \mathbb{R}^2$ such that $\varphi(\mathfrak{a})$ is a vertex of the cut graph of $\varphi^{-1}$ and the edge adjacent to $\varphi(\mathfrak{a})$ is an incoming edge, see Figure 7.

The weight of the incoming edge at $\varphi(\mathfrak{a})$ is called the **multiplicity** of $\mathfrak{a}$.

The requirement that a node $\mathfrak{a}$ admits a simple nodal chart as in Definition 2.10 is equivalent to a node admitting a model neighbourhood as in [Sym03, Proposition 4.14], meaning that Definition 2.10 is equivalent to [Sym03, Definition 5.1].

*Remark* 2.11. Let $\mathfrak{n} \in U$ be the only node in some pre-nodal integral affine $U$. It is not enough to fix the monodromy around $\mathfrak{n}$ to determine the integral affine structure of $U \setminus \{\mathfrak{n}\}$: Take $E = \mathbb{R}^2$ as a smooth manifold. We give $E$ the structure of a pre-nodal integral affine manifold. Let $\psi : \mathbb{R} \times \mathbb{R}_{>0} \to E \setminus \{0\}$ be polar coordinates on $E \setminus \{0\}$. Equip $E \setminus \{0\}$ with pullback integral affine structure along $\psi^{-1}$. This makes $E$ into a pre-nodal integral affine surface with a node at 0 and trivial monodromy around 0. This $E$ does not admit a nodal chart near 0, and clearly $E$ is not isomorphic to $\mathbb{R}^2$ with the standard integral affine structure as a pre-nodal integral affine manifold (for example $E$ has closed straight lines). The main point of Definition 2.10 is to restrict nodes to have the same model neighbourhood as in [Sym03, Definition 5.1].

*Remark* 2.12. One could be tempted to include "anti-nodes" in Definition 2.10, that is nodes with an outgoing edge of $G$ at $\varphi(\mathfrak{a})$. These anti-nodes do not seem to have a corresponding concept in symplectic geometry or algebraic geometry (e.g. we would have to rule them out in Theorems 2.27 and 2.30), so we already rule them out in the definition.

We may think of the nodes as discrete points of positive curvature in the otherwise flat $B$. We have the following "affine Gauss-Bonnet theorem":

**Theorem 2.13** ([KS06, Theorem 2]). *If $B$ is compact and without boundary,*

$$\frac{1}{12}\sum_{\mathfrak{n}\in\mathfrak{N}} k_{\mathfrak{n}} = \chi(B)$$

*where $k_{\mathfrak{n}}$ is the multiplicity of $\mathfrak{n}$, and $\chi$ denotes the Euler characteristic.*

*Remark* 2.14. This leads to a very limited topology on nodal integral affine surfaces: If $B$ is compact (perhaps with boundary/corners) then $B$ is either a disc, an annulus, a Möbius band, a Klein bottle with no nodes, a torus with no nodes, a projective plane with 12 nodes or a sphere with 24 nodes. See [Zun03, Proposition 3.3]. [LS10] determines the corresponding symplectic manifolds up to diffeomorphism.

Given a topological covering $p: E \to X$, for a chosen base point $x \in X$ we get the *monodromy action* $\pi_1(X, x) \to \operatorname{Aut}(p^{-1}(x))$ obtained by lifting loops based at $x$ to $E$. We are interested in the monodromy of the lattice bundle $\Lambda B \to (B \setminus \mathfrak{N})$ around a node.

Take a simple nodal chart $\varphi: U_{\mathfrak{a}} \to \mathbb{R}^2$ at $\mathfrak{a} \in \mathfrak{N}$ with the cut graph $G$ of $\varphi^{-1}$ and equip $U_{\mathfrak{a}}$ with the orientation induced by $\mathbb{R}^2$. Let $\ell$ be the incoming edge of $G$ at $\mathfrak{a}$, $v' \in \mathbb{Z}^2$ its primitive direction vector of $\ell$ and $k$ the weight of $\ell$. Then the monodromy of $\Lambda B$ induced by a positive loop around $\mathfrak{a}$ based at $x \in U_{\mathfrak{a}} \setminus \{\mathfrak{a}\}$ is given by the shear map $(D_x\varphi)^{-1} \circ s^k_{v'} \circ D_x\varphi$. Here $D\varphi_x$ is defined if $x$ does not lie on a cut of $\varphi$. As $v'$ is fixed by the monodromy, the pull-back $v = D\varphi^{-1}[v']$ gives a well-defined vector field on $U_{\mathfrak{a}}$, that is there is a unique continuous vector field $v$ such that $D\varphi[v] = v'$ whenever $D\varphi$ is defined. Similarly $\lambda = \varphi^*(\det(v', \cdot))$ gives a well-defined covector field on $U_{\mathfrak{a}}$. We can thus write the monodromy of a positive loop around $\mathfrak{a} \in B$ as

$$\operatorname{id} - kv\lambda\ .$$

The monodromy is independent of the choice of simple nodal chart, so the pair $(v, \lambda)$ is unique up to sign. We call the pair $(v, \lambda)$ the **monodromy pair of $\mathfrak{a}$**, $v$ the **monodromy vector field** and $\lambda$ the **monodromy covector field**.[4]

The integral lattice $\Lambda B$ also induces a connection on $TB \to (B \setminus \mathfrak{N})$: By lifting paths in $B \setminus \mathfrak{N}$ to $\Lambda B$, we get a parallel transport system on $\Lambda B \to B \setminus \mathfrak{N}$, which can be linearly extended to $TB \to B \setminus \mathfrak{N}$. A geodesic with respect to this connection is a path $\gamma$ such that $\dot{\gamma}$ is constant relative to $\Lambda B$, which motivates the following definition:

**Definition 2.15.** A **straight line** $\gamma$ in $B \setminus \mathfrak{N}$ is a geodesic with respect to the connection induced by the integral lattice $\Lambda B$.

- We say $\gamma$ is **maximal** if it cannot be extended in $B \setminus \mathfrak{N}$.
- It is **rational** if $r\dot{\gamma} \in \Lambda B$ for some real number $r > 0$.

[4]Note that here $\lambda$ has only an auxiliary function: If we choose an orientation around $\mathfrak{a}$, we can choose $\lambda = \det_{U_{\mathfrak{a}}}(v, \cdot)$, where $\det_{U_{\mathfrak{a}}}$ is given by pulling back det along a positive integral affine chart.

- It is **primitive** if $\dot{\gamma} \in \Lambda B$ is primitive.

*Remark* 2.16. A straight line might have self-intersections.

We think of primitive straight lines as having "unit speed":

**Definition 2.17.** Let $\gamma : (a, b) \to B \setminus \mathfrak{N}$ be a rational straight line with $r\dot{\gamma} \in \Lambda B$ primitive. Its **integral affine length** is given by $\frac{b-a}{r}$.

Note that this notion of length does not induce a well-defined metric.

**Definition 2.18.** Let $\mathfrak{a} \in \mathfrak{N}$. The **eigenline** $\gamma$ **of** $\mathfrak{a}$ is the one-dimensional subset of $B$ with $\mathfrak{a} \in \gamma$ such that $\gamma \setminus \{\mathfrak{a}\}$ is the union of the two maximal straight lines which are parallel to the monodromy vector field of $\mathfrak{a}$ near $\mathfrak{a}$ and whose closure includes $\mathfrak{a}$.

*Remark* 2.19. The eigenline of a node might be a circle, in which case the two maximal straight lines in the definition above coincide.

We will often want to know how images of straight lines in $B$ look like in a nodal chart diagram. The following remark is a useful tool for this:

*Remark* 2.20. Take a nodal chart $\varphi : B \supset U \to \mathbb{R}^2$ and let $\gamma : (-\varepsilon, \varepsilon) \to B \setminus \mathfrak{N}$ be a straight line segment with $\gamma(0) = x \in U$. Then $\varphi \circ \gamma$ is a piecewise affine path, and $R_x \circ \varphi \circ \gamma$ is an affine path. Denote by $(\varphi \circ \gamma)'(0^{\pm})$ the upper resp. lower derivative of $\varphi \circ \gamma$ at 0. Then we have

$$R_x((\varphi \circ \gamma)'(0^+)) = -R_x(-(\varphi \circ \gamma)'(0^-)) ,$$

since $R_x$ maps the vectors $\pm(\varphi \circ \gamma)'(0^{\pm})$ to $\pm(R_x \circ \varphi \circ \gamma)'(0)$.

Together with Remark 2.8, this translates to the following recipe: Let $\ell$ be a cut with primitive direction vector $v$. Approaching $\ell$ walking on a straight line in direction $u$, after going through $\ell$ we will exit on the other side of $\ell$ with direction $u + |\det(v, u)|v$.

## 2.2 Nodal tangles

A nodal slide as in [Sym03, §6.1] can be described as two nodal integral affine surfaces that are identical except for one node which has been displaced along its eigenline. The following definition captures a notion of isotopy between two nodal integral affine surfaces constructed by a sequence of nodal slides.

**Definition 2.21.** A **nodal tangle** is a pre-nodal integral affine manifold $(B \times I, \mathfrak{N})$, where $B$ is a surface with corners and $I$ is an interval with its integral affine structure inherited from $\mathbb{R}$, such that:

- $B \setminus \pi_B(\mathfrak{N})$ is an integral affine surface such that $(B \setminus \pi_B(\mathfrak{N})) \times I$ is an integral affine submanifold of $B \times I \setminus \mathfrak{N}$.

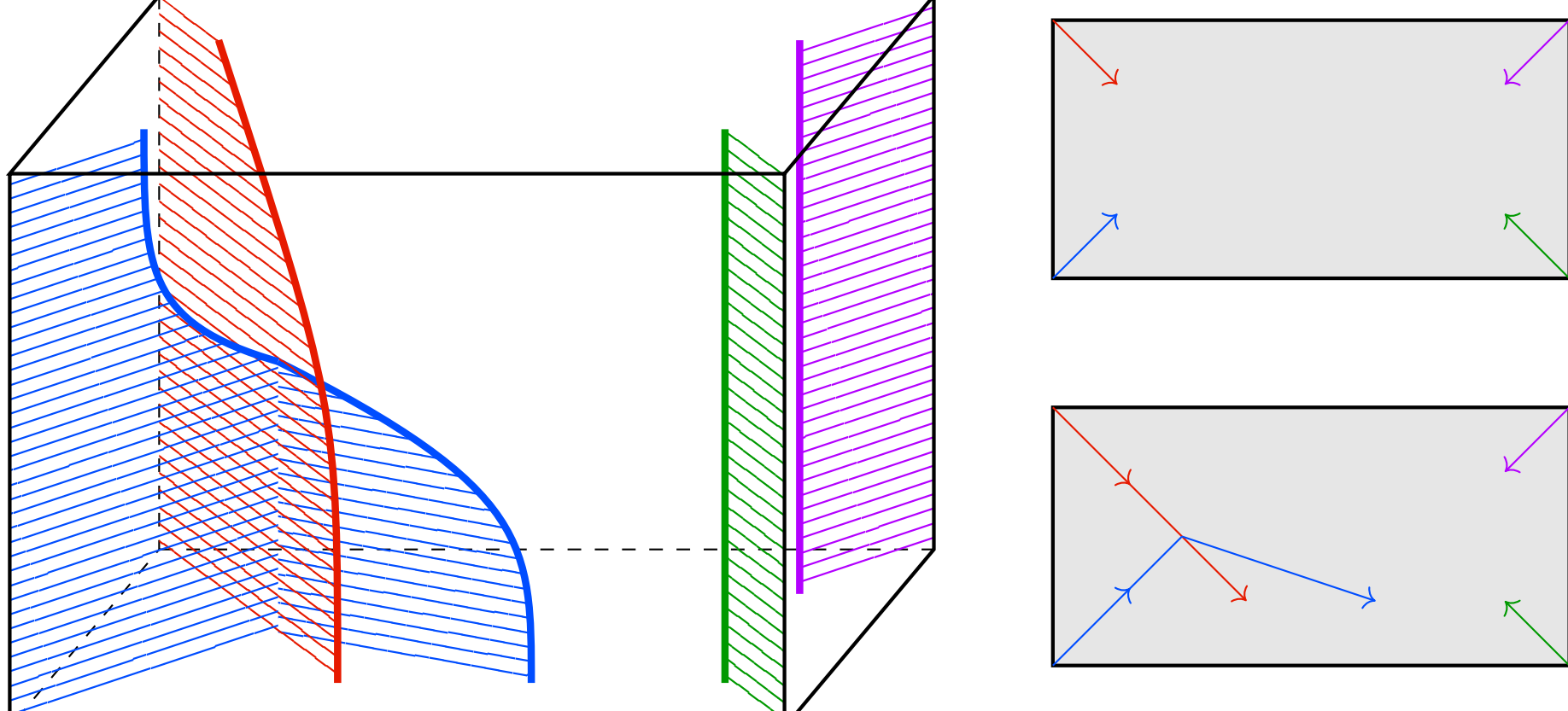

Figure 8: A nodal tangle. Horizontal slices are nodal chart diagrams. On the right are the top and bottom nodal chart diagrams. In the image on the left, the nodes $\mathfrak{N}$ are given by the thick coloured curves.

- For every $t \in I$, $(B_t, \mathfrak{N}_t) = (B \times \{t\}, \mathfrak{N} \cap B \times \{t\})$ is a nodal integral affine surface which intersects $\mathfrak{N}$ transversely.
- For $t \in I$ and $\mathfrak{a} \in \mathfrak{N}_t$, let $\lambda$ be the monodromy covector field of $\mathfrak{a} \in (B_t, \mathfrak{N}_t)$. Then for any 1-dimensional stratum $S$ of $\mathfrak{N}$ at $(\mathfrak{a}, t)$ we have $D\pi_B(T_{(\mathfrak{a},t)}S) \subset \ker \lambda$.

Here $\pi_B : B \times I \to B$ denotes the projection to $B$.

We can think of a nodal tangle as a movie in which nodes can slide along their eigenlines, and a node of multiplicity $k$ can split into multiple nodes with the same monodromy pair whose multiplicities add up to $k$. The first two conditions ensure that the integral affine structure on $B \setminus \pi_B(\mathfrak{N})$ remains fixed and that $\mathfrak{N}$ is not tangent to any $B_t$. The third condition ensures that nodes slide along their eigenline, as $\ker \lambda$ is spanned by the monodromy vector field $v$.

**Definition 2.22.** If $(B \times [a, b], \mathfrak{N})$ is a nodal tangle we say $(B_a, \mathfrak{N}_a)$ and $(B_b, \mathfrak{N}_b)$ are **connected by a nodal tangle**.

The following lemma states that this definition is compatible with [Sym03, §6.1]:

**Lemma 2.23.** *Two integral affine surfaces $(A, \mathfrak{N}_A), (C, \mathfrak{N}_C)$ are related by a finite sequence of nodal slides (see [Sym03, Definition 6.1/Exercise 6.7]) if and only if they are connected by a nodal tangle.*

*Proof.* Given $(A, \mathfrak{N}_A), (C, \mathfrak{N}_C)$ related by a single nodal slide we can easily construct a nodal tangle between $(A, \mathfrak{N}_A)$ and $(C, \mathfrak{N}_C)$. For a sequence of nodal slides we can stack the nodal tangles thus obtained on top of each other.

Given a nodal tangle $(B \times I, \mathfrak{N})$, for a subinterval $J \subset I$ let $(B \times J, \mathfrak{N}_J = \mathfrak{N} \cap B \times J)$ be the nodal tangle restricted to $J$. Restricting to a small enough subinterval $J = [a, b] \subset I$ the set $\pi_B(\mathfrak{N}_J)$ is a disjoint union of closed straight line segments in $B$, one for every connected component of $\mathfrak{N}_J$.

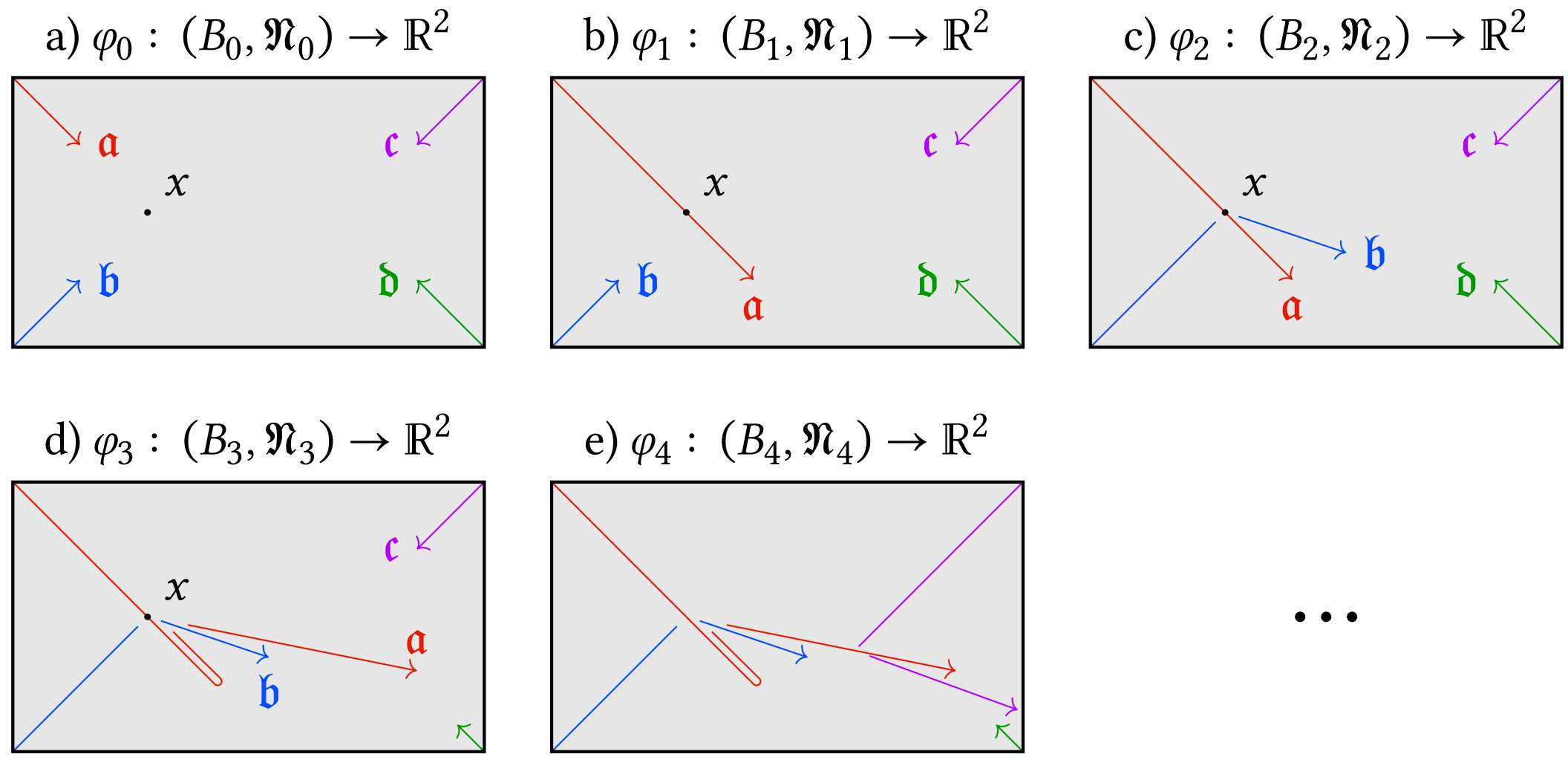

Figure 9: Building a tangle.

For connected components of $\mathfrak{N}_J$ which have no stratification (no nodes splitting), [Sym03, Definition 6.1] now gives a suitable nodal slide. For stratified components of $\mathfrak{N}_J$, we extend the definition of nodal slide [Sym03, Definition 6.1] to the case mentioned in [Sym03, Exercise 6.7].

Any order of these nodal slides gives a sequence of nodal slides transforming $B_a$ into $B_b$.

Since $I$ is compact by Definition 2.22, we only need finitely many subintervals $J$ to cover $I$. □

**Definition 2.24.** Let $(B \times I, \mathfrak{N})$ be a nodal tangle, and let $a, b \in I$. The **transition map** is the map $\tau_a^b : B_a \to B_b$ given by the composition $i_b \circ \pi_B|_{B_a}$, where $i_b$ is the natural homeomorphism $B \to B \times \{b\} = B_b$.

$\tau_a^b|_{B_a \setminus (\mathfrak{N}_a \cup \tau_b^a \mathfrak{N}_b)}$ is a piecewise integral affine map, with the restriction $\tau_a^b|_{B_a \setminus \pi_B \mathfrak{N}}$ being integral affine. It describes the changes in the integral affine structure induced by the nodal tangle. If $\varphi_a$ is an integral affine chart near $x_a \in B_a$, then $\varphi_b = \varphi_a \circ \tau_b^a$ is a nodal chart near $x_b = \tau_a^b(x_a)$, whose cut graph $G_b$ is contained in $i_a(\pi_B(\mathfrak{N}))$ (although $\varphi_b$ may still be an integral affine chart, in particular if $\pi_B(x_a) \notin \pi_B(\mathfrak{N})$).

**Example 2.25** (Building a nodal tangle)**.** Figure 9 a) shows a nodal chart diagram of a global simple nodal chart $\varphi_0 : (B_0, \mathfrak{N}_0) \to \mathbb{R}^2$, with four cuts ending in the nodes $\mathfrak{a}, \mathfrak{b}, \mathfrak{c}, \mathfrak{d}$ of multiplicity 1. The position of the nodes is marked by arrow tips. The direction of the arrow tips correspond to the orientation mentioned in Remark 2.8, and to the direction of the monodromy vector field of the node.

Let's modify $(B_0, \mathfrak{N}_0)$ by nodal slides to create a nodal tangle: First slide $\mathfrak{a}$ such that $x$ is contained in its cut, giving the new nodal chart diagram for $\varphi_1 := \varphi_0 \circ \tau_1^0 : (B_1, \mathfrak{N}_1) \to \mathbb{R}^2$ as seen in Figure 9 b). In doing so we modified the integral affine structure along $\mathfrak{a}$'s cut. The rectifying map at the marked point $x$ of $\varphi_1$ is now no longer trivial since $x$ sits on a cut. Concretely, let $v_1 = \left(\begin{smallmatrix} 1 \\ -1 \end{smallmatrix}\right)$ be the direction of $\mathfrak{a}$'s cut, then $R_x$ is the half-shear $h_{v_1}$, defined in (1).

Moving on, we can slide $\mathfrak{b}$ along its eigendirection, past the point $x$. As our nodal chart $\varphi_1$ has a cut through $x$, eigenline of $\mathfrak{b}$ appears broken in $\operatorname{im} \varphi_1$ at the point $x$. Using the recipe in Remark 2.20, we find that the incoming vector $v_2' = \left(\begin{smallmatrix} 1 \\ 1 \end{smallmatrix}\right)$ exits $x$ as $v_2 = \left(\begin{smallmatrix} 3 \\ -1 \end{smallmatrix}\right)$:

$$v_2 = v_2' + |\det(v_1, v_2)| v_1 = \begin{pmatrix} 1 \\ 1 \end{pmatrix} + 2 \begin{pmatrix} 1 \\ -1 \end{pmatrix} = \begin{pmatrix} 3 \\ -1 \end{pmatrix} .$$

We get the nodal chart diagram Figure 9 c) for $\varphi_2 : (B_2, \mathfrak{N}_2) \to \mathbb{R}^2$. The rectifying map at $x$ of $\varphi_2$ is now the composition $h_{v_2} \circ h_{v_1}$.

The tangle built up to this point is the one depicted in Figure 8, and we can obtain the diagram Figure 9 c) by projecting the left drawing in Figure 8 onto $B$. We indicate the order of the two nodal slides through $x$ by drawing the later slide "underneath" the first, as in knot diagrams.

Sliding $\mathfrak{a}$ back through $x$, it is once again deflected, giving Figure 9 d). Using Remark 2.20 again, we determine the outgoing vector $v_3$ pointing from $x$ to the node $\mathfrak{a}$ to be

$$v_3 = -v_1 + |\det(-v_1, v_2)| v_2 = \begin{pmatrix} -1 \\ 1 \end{pmatrix} + 2 \begin{pmatrix} 3 \\ -1 \end{pmatrix} = \begin{pmatrix} 5 \\ -1 \end{pmatrix} ,$$

giving the new rectifying map $R_x = h_{v_3} \circ h_{v_2} \circ h_{v_1}$. Here the "hook" at $x$ of $\mathfrak{a}$ going back along its path is technically not a cut of $\varphi_3$. However we still choose to draw it and think of the diagram as a projection of the nodal tangle $(B \times [0, 3], \mathfrak{N})$.

Proceeding in the same manner we can continue to create more and more complicated nodal tangles.

Another example of a more complicated nodal tangle can be seen in Figure 10.

Let $B \times [0, 1]$ be a nodal tangle, $U_0 \subset B_0$ open, and $\varphi_0 : U_0 \to \mathbb{R}^2$ a simple nodal chart. Perturb $\mathfrak{N}$ (disregarding the integral affine structure) such that $\pi_B|_{\mathfrak{N}}$ is an immersion. The image of $\varphi_0 \circ i_0 \circ \pi_B$ with over/under crossings of $\mathfrak{N}$ marked (viewed from $t = 0$ towards $t = 1$) as in a knot diagram is called a **nodal tangle diagram**.

See for example the diagrams in Figures 9 and 10. The perturbation of $\mathfrak{N}$ creates the hooks. Due to the imprecision introduced by perturbing $\mathfrak{N}$, nodal tangle diagrams are mainly useful for illustration purposes. However, nodal tangle diagrams are still enough to reconstruct the nodal tangle $B \times I$ up to an appropriate notion of isotopy of nodal tangles. We mainly use nodal tangle diagrams for illustrative purposes, as they can quickly become very hard to read in all but the simplest cases.

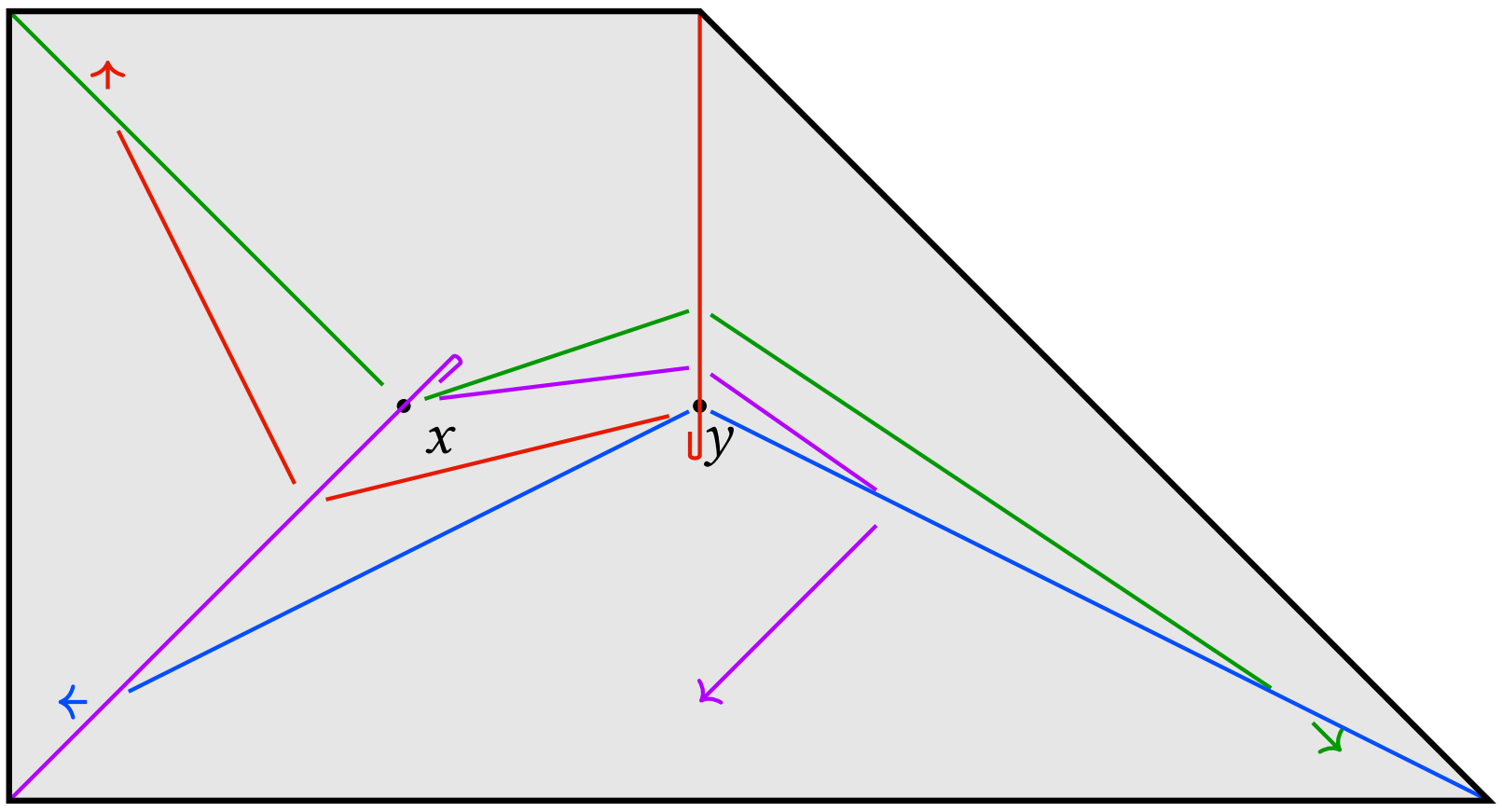

Figure 10: A more complicated tangle of nodes. The small hooks at $x$ and $y$ illustrate the node reversing its sliding direction, not an actual lateral movement.

## 2.3 Almost toric fibrations

**Definition 2.26.** Let $(X, \omega)$ be a 4-dimensional symplectic manifold, and $B$ a smooth surface. An **almost toric fibration** is an integrable system $\pi : X \to B$ with compact connected fibres, non-degenerate singularities and no hyperbolic singularities.

Here, by an integrable system we mean a fibration $\pi : X \to B$ such that for any two smooth functions $f, g$ on $B$, we have $\{\pi^* f, \pi^* g\} = 0$. In other words, for any smooth chart $\varphi : U \to \mathbb{R}^2$ on $B$, the map $\varphi \circ \pi|_{\pi^{-1}(U)}$ is a complete integrable system. See e.g. [Zun03].

An almost toric fibration induces on $B$ the structure of a nodal integral affine surface: Denote by $\mathfrak{N}$ the image of the focus-focus singularities of $\pi$. Then the Arnold-Liouville Theorem gives a set of charts on $B \setminus (\mathfrak{N} \cup \partial B)$, extending to $B \setminus \mathfrak{N}$, called *action-coordinates*, with transition maps in $\mathrm{Aff}_2(\mathbb{Z})$. The pullback $\pi^* f$ of a local smooth function $f : B \supset U \to \mathbb{R}$ induces a Hamiltonian action on $\pi^{-1}(U)$. Similarly, any 1-form $\lambda \in T^*B$ is locally exact, and thus induces a local Hamiltonian action. We have $\lambda \in \Lambda^* B$ if and only if $\lambda$ induces a local Hamiltonian $S^1$-action. The lattice $\Lambda^* B$ induces a local free Hamiltonian torus action near points of the **regular locus** $B \setminus (\partial B \cup \mathfrak{N})$.

The neighbourhood of a node $\mathfrak{n} \in \mathfrak{N}$ is described in [Zun97]; the existence of simple nodal charts is discussed in [Sym03, Section 4.4]. The fact that the monodromy covector field $\lambda$ of a node $\mathfrak{n}$ extends over $\mathfrak{n}$ corresponds to the fact that $\pi^* \lambda$ induces an $S^1$ action on a neighbourhood of $\pi^{-1}(\mathfrak{n})$.

The fibre over a point $x \in B$ is a Lagrangian torus if $x$ is a regular point, a circle if $x \in \partial B$ and not a corner, a point if $x$ is a corner of $B$, and if $x$ is a node of multiplicity $k$ the fibre is a torus pinched $k$ times along the monodromy covector field $\lambda$, that is we take $k$ disjoint orbits of $\lambda$ in a regular fibre nearby and collapse them to $k$ points.

In certain cases a nodal integral affine surface uniquely determines an almost toric fibration:

**Theorem 2.27.** *Let $(B, \mathfrak{N})$ be an nodal integral affine surface. Then an almost toric fibration $\pi : X \to (B, \mathfrak{N})$ exists.*

*If additionally $B$ is homotopy equivalent to a punctured surface, the almost toric fibration is unique in the following sense: Suppose two almost toric fibrations $\pi_1 : X_1 \to B$ and $\pi_2 : X_2 \to B$ induce the same nodal integral affine structure on $B$. Then for any neighbourhood $U$ of $\mathfrak{N}$ there exists a symplectomorphism $\psi : X_1 \to X_2$ such that*

$$\begin{array}{ccc} X_1 \setminus \pi_1^{-1}(U) & \xrightarrow{\psi} & X_2 \setminus \pi_2^{-1}(U) \\ & {\scriptstyle \pi_1}\searrow \quad \swarrow{\scriptstyle \pi_2} & \\ & B & \end{array}$$

*commutes, that is, except on a neighbourhood $U$ of the nodes $\mathfrak{N}$, $\psi$ is a* fibred symplectomorphism, .

The existence is shown in [Sym03, Theorem 5.2], and uniqueness is shown in [Eva23, Proof of Theorem 8.5].

*Remark* 2.28. The ambiguity for the uniqueness comes only from the regular part $\widetilde{B} = B \setminus (\mathfrak{N} \cup \partial B)$, and is discussed in [Dui80, Section 2]. We summarize the discussion: The quotient $T^*\widetilde{B}/\Lambda^*\widetilde{B}$ carries a natural symplectic structure. Let $\mathscr{L}(T^*\widetilde{B}/\Lambda^*\widetilde{B})$ denote the sheaf of Lagrangian sections $\widetilde{B} \to T^*\widetilde{B}/\Lambda^*\widetilde{B}$. Then for any element $\mu$ in the Čech cohomology $H^1(B, \mathscr{L}(T^*\widetilde{B}/\Lambda^*\widetilde{B}))$, the data of $(B, \mathfrak{N}, \mu)$ determines a unique almost toric fibration in the sense of Theorem 2.27. See also [Zun03] for putting this into the context of integrable systems with singularities.

*Remark* 2.29. [San03] describes the possible almost toric fibrations in the neighbourhood of a node up to fibred symplectomorphisms, which can be classified by a certain Taylor series. By allowing $\psi$ to be not fibred near the nodes, we can ignore this subtlety.

A nodal tangle can be lifted to a path of almost toric fibrations:

**Theorem 2.30.** *Let $(B \times I, \mathfrak{N})$ be a nodal tangle, $a \in I$ and $\pi_a : X \to (B_a, \mathfrak{N}_a)$ an almost toric fibration. Then for any neighbourhood $U \subset B$ of $\pi_B(\mathfrak{N})$, we find a fibration $\pi : X \times I \to B \times I$ such that:*

- $\pi|_{X \times \{a\}} = \pi_a$,
- *for all $t \in I$, $\pi_t = \pi|_{X \times \{t\}} : X \to (B_t, \mathfrak{N}_t)$ is an almost toric fibration,*
- *for all $t \in I$, $X \setminus \pi_t^{-1}(U) = X \setminus \pi_a^{-1}(U)$ and $\pi_t|_{X \setminus \pi_a^{-1}(U)} = \pi_a|_{X \setminus \pi_a^{-1}(U)}$.*

Using the correspondence Lemma 2.23, this theorem is a direct consequence of [Eva23, Proof of Theorem 8.10]. See also the earlier [Sym03, Theomrem 6.5].

**Definition 2.31.** Using the notation of Theorem 2.30, we say that we can **lift** a nodal tangle $(B \times I, \mathfrak{N})$ to a path of almost toric fibrations $\pi_t : X \to (B_t, \mathfrak{N}_t)$ **supported** on $U \subset B$.

The following useful lemma shows that although a lift of a nodal tangle modifies the fibres near $\pi_B(\mathfrak{N})$, for any given fibre outside of $\pi_B(\mathfrak{N})$ this can be "undone" by a symplectomorphism:

**Lemma 2.32.** *Let $\pi : X \times [0,1] \to (B \times [0,1], \mathfrak{N})$ be a lift of a nodal tangle. If $B$ is homotopy equivalent to a punctured surface, then for all $x \in B \setminus \pi_B(\mathfrak{N})$ there exists $\psi \in \mathrm{Symp}(X)$ such that $\pi_0 = \pi_1 \circ \psi$ near $\pi_0^{-1}(x)$.*

*Proof.* Extend the nodal tangle by concatenating it with its reverse: set $(B_t, \mathfrak{N}_t) = (B_{2-t}, \mathfrak{N}_{2-t})$ for $t \in [1,2]$, giving a nodal tangle $(B \times [0,2], \mathfrak{N})$.

Let $x \in B \setminus \pi_B(\mathfrak{N})$ and choose a neighbourhood $U$ of $\pi_B(\mathfrak{N})$ such that $x \notin U$. Then by Theorem 2.30 we may lift the nodal tangle $(B \times [1,2], \mathfrak{N})$ to $\pi_t : X \to (B_t, \mathfrak{N}_t), t \in [1,2]$ with the lift being supported on $U$. In particular we have $\pi_2 = \pi_1$ near $x$. Since $(B_0, \mathfrak{N}_0) = (B_2, \mathfrak{N}_2)$, applying Theorem 2.27 $\pi_0$ and $\pi_2$, we get a symplectomorphism $\psi \in \mathrm{Symp}(X)$ such that $\pi_0 = \pi_2 \circ \psi = \pi_1 \circ \psi$ near $x$. $\square$

*Remark* 2.33. Taking more care in the proof of Theorem 2.30, it should be possible to build lifts such that the Lagrangian isotopy $t \mapsto \pi_t^{-1}(x), t \in [0,1]$ is generated by a Hamiltonian isotopy, and to get $\psi \in \mathrm{Ham}(X)$ without the restriction on $B$ in Lemma 2.32.

## 2.4 Nodal tangles and invariant germs

Let $\mathcal{L}$ be the space of Lagrangians in $X$ equipped with the $\mathscr{C}^1$-topology, see [Ono08] for details on the $\mathscr{C}^1$-topology on the space of Lagrangians. Let $\mathcal{I} : \mathcal{L} \to A$ be a symplectic invariant of Lagrangian submanifolds with values in some set $A$, meaning that for all $\varphi \in \mathrm{Symp}(X)$ and $L \in \mathcal{L}$ we have $\mathcal{I}(\varphi(L)) = \mathcal{I}(L)$.

During a nodal tangle invariants don't change on fibres of $B \setminus \pi_B(\mathfrak{N})$:

**Lemma 2.34.** *Let $\pi : X \times I \to (B \times I, \mathfrak{N})$ be a lift of a nodal tangle, and $\mathcal{I} : \mathcal{L} \to A$ a symplectic invariant. If $B$ is homotopy equivalent to a punctured surface, then for $a, b \in I$ and $x_a \in B_a \setminus \pi_B(\mathfrak{N}), x_b = \tau_a^b x_a$, we have*

$$\mathcal{I}(\pi_a^{-1}(x_a)) = \mathcal{I}(\pi_b^{-1}(x_b)) .$$

*Proof.* Let $a, b \in I$ and $x \in B \setminus \pi_B(\mathfrak{N})$. Then by Lemma 2.32 there exists $\psi \in \mathrm{Symp}(X)$ such that $\pi_a = \pi_b \circ \psi$ near $\pi_a^{-1}(x)$. In particular $\psi(\pi_a^{-1}(x)) = \pi_b^{-1}(x)$, so

$$\mathcal{I}(\pi_a^{-1}(x)) = \mathcal{I}(\psi(\pi_a^{-1}(x))) = \mathcal{I}(\pi_b^{-1}(x)) . \qquad \square$$

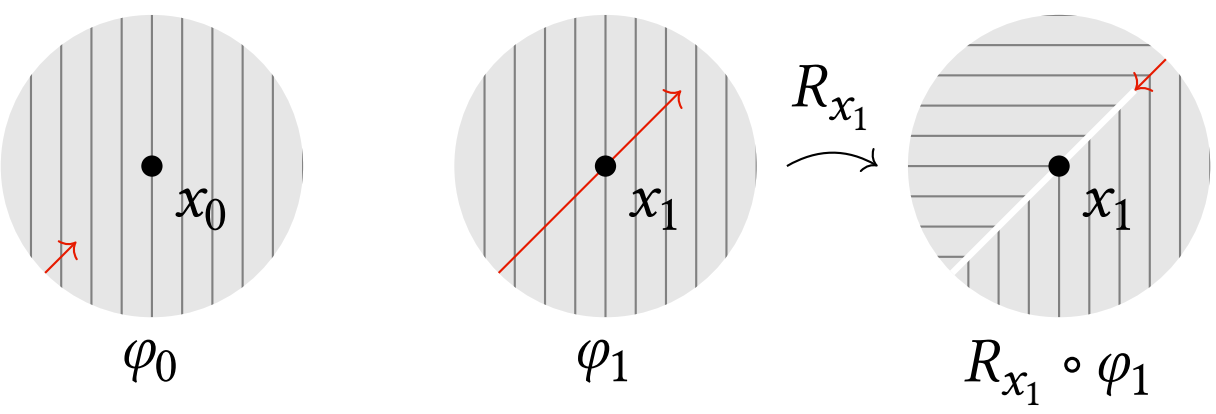

Figure 11: Modifying an invariant germ by a simple nodal tangle. We start with a nodal chart $\varphi_0$ which is affine around $x_0$ for $B_0$, then we modify $B_0$ by a nodal tangle to obtain $B_1$ and with $x_1 = \tau_0^1 x_0$, first giving us a nodal chart $\varphi_1 = \varphi_0 \circ \tau_1^0$. We then compose $\varphi_1$ with the rectifying map at $x_1$ to get a nodal chart $R_x \circ \varphi_1$ which is an integral affine chart near $x_1$. The grey lines are level sets of an invariant $\mathscr{I} : \mathscr{L} \to \mathbb{R}$. The difference of the germs $[\mathscr{I}]_{x_0}$ and $[\mathscr{I}]_{x_1}$ is captured by Corollary 2.35.

For an almost toric fibration $\pi : X \to B$ and a point $x \in B$, abusing notation, denote by $[\mathscr{I}]_x$ the germ of the function $(\mathscr{I} \circ \pi^{-1}) : B^{\mathrm{reg}} := B \setminus (\partial B \cup \mathfrak{N}) \to A$, where we think of $\pi^{-1}$ as a map $B^{\mathrm{reg}} \to \mathscr{L}$.

While we cannot directly describe how $\mathscr{I}$ changes for points $x \in \pi_B(\mathfrak{N})$ of a nodal tangle, we can partially describe how its germ $[\mathscr{I}]_x$ changes:

**Corollary 2.35.** *Let $\pi : X \times I \to (B \times I, \mathfrak{N})$ be a lift of a nodal tangle, and $\mathscr{I} : \mathscr{L} \to A$ a symplectic invariant. If $B$ is homotopy equivalent to a punctured surface, then for $a, b \in I$ and $x_a \in B_a$, $x_b = \tau_a^b x_a$, we have*

$$([\mathscr{I}]_{x_a})|_{B_a \setminus \pi_B(\mathfrak{N})} = ([\mathscr{I}]_{x_b} \circ \tau_a^b)|_{B_a \setminus \pi_B(\mathfrak{N})} .$$

*Remark* 2.36. More loosely we can say that $[\mathscr{I}]_{x_a} = [\mathscr{I}]_{x_b} \circ \tau_a^b$ on an open dense neighbourhood of $x_a$.

*Proof.* This follows directly from Lemma 2.34 by taking germs. □

See Figure 11 for an example.

We can show that the germ $[\mathscr{I}]_x$ only depends on the Lagrangian $\pi^{-1}(x) \subset X$. This requires a bit of setup:

Given a Lagrangian $L$, using a *versal deformation* $v : H^1(L; \mathbb{R}) \to \mathscr{L}$ we can identify a small neighbourhood of $L$ in $\mathscr{L}$ up to local Hamiltonian isotopy with a small neighbourhood of $0$ in $H^1(L; \mathbb{R})$. Versal deformations of Lagrangians were introduced in [Che96]. For details we refer the reader to the exposition in [BHS24, Section 3.3]. Denote by $[\mathscr{I}]_L : H^1(L; \mathbb{R}) \dashrightarrow A$ the germ of $\mathscr{I} \circ v$ at $0$, as defined in [BHS24, Definition 3.13], which is independent of the choice of versal deformation. This can be considered a refinement of $\mathscr{I}$ by considering values of $\mathscr{I}$ in a small neighbourhood of $L$.

Let $\pi : X \to (B, \mathfrak{N})$ be an almost toric fibration. For every regular value $x$ of $\pi$ we have a natural identification $\Lambda_x B \cong H^1(\pi^{-1}(x); \mathbb{Z})$: The map

$$\begin{aligned} \Lambda_x^* B \times H^1(\pi^{-1}(x); \mathbb{Z}) &\to \mathbb{Z} \\ (df, [\alpha]) &\mapsto \omega(X_{\pi^* f}, X_\alpha) \end{aligned}$$

is a perfect pairing, giving an isomorphism $H^1(\pi^{-1}(x); \mathbb{Z}) \cong \Lambda_x B$.[5] This extends to an integral linear isomorphism[6] $\alpha : H^1(\pi^{-1}(x); \mathbb{R}) \to T_x B$. Furthermore, the lattice bundle $\Lambda B$ induces a flat connection on $B$ such that the exponential map $\exp TB \to B$ locally gives integral affine maps $\exp_x : T_x B \dashrightarrow B$ near $0 \subset T_x B$. (We use $\dashrightarrow$ to denote a map only defined near 0.) Then the composition

$$H^1(\pi^{-1}(x); \mathbb{R}) \xrightarrow{\ \alpha\ } T_x B \overset{\exp_x}{\dashrightarrow} B \xrightarrow{\ \pi^{-1}\ } \mathscr{L}$$

is a versal deformation. See [BHS24, Section 3.4] for more details.

We summarize the discussion in the following lemma:

**Lemma 2.37.** *The germ $[\mathscr{J}]_x$ and the invariant germ $[\mathscr{J}]_{\pi^{-1}(x)}$ are related by the commuting diagram*

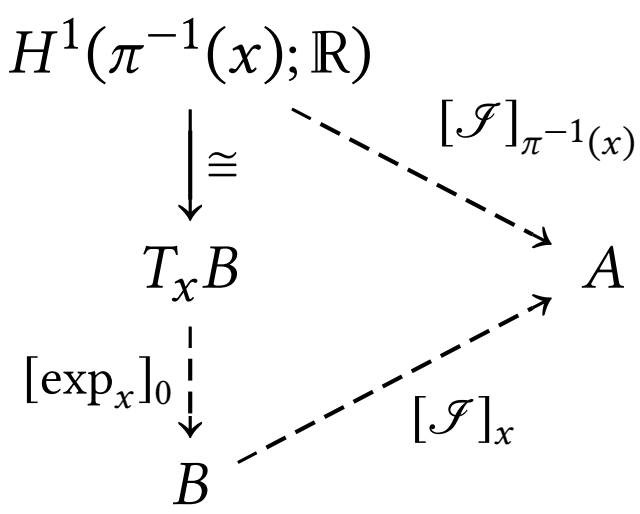

*where the vertical maps are integral affine. In particular $[\mathscr{J}]_x$ depends only on the Lagrangian $\pi^{-1}(x)$.*

# 3 Weakly Delzant polygonal domains

In this section, we combine the framework of Section 2 with [MS23] to develop a *canonical form* (Proposition 3.17) for *weakly Delzant polygonal domains*, of which Delzant polygons are important examples. In Section 3.2, we use this canonical form to prove [Sym03, Conjecture 6.8] in the case of closed toric symplectic four-manifolds. In Section 3.4, we use this canonical form and probes introduced in [McD11] to calculate the displacement energy of toric fibres.

[5] The dual isomorphism $\Lambda_x^* B \to H_1(\pi^{-1}(x); \mathbb{Z})$ has a nice interpretation: It maps an $S^1$-action in $\Lambda_x^* B$ to the class of its orbit.

[6] Meaning an isomorphism of pairs $\big(H^1(\pi^{-1}(x); \mathbb{R}), H^1(\pi^{-1}(x); \mathbb{Z})\big) \to (T_x B, \Lambda_x B)$.

## 3.1 Weakly Delzant polygonal domains and their caustics

For $c \in \mathbb{R}$ and a primitive vector $\lambda \in \mathbb{Z}^2$, define the half-plane

$$H_{\lambda,c} = \{x \in \mathbb{R}^2 \mid \langle \lambda, x \rangle + c \geq 0\} .$$

**Definition 3.1.** A set $\Delta \subset \mathbb{R}^2$ is a **polygonal domain** if it has non-empty interior and for every compact set $C \subset \mathbb{R}^2$, $\Delta \cap C$ is the intersection of finitely many half-planes and $C$.

We can write

$$\Delta = \bigcap_{H_{\lambda,c} \in \mathscr{H}} H_{\lambda,c}$$

for a minimal set of half-planes $\mathscr{H}$.

*Remark* 3.2. $\partial\Delta$ is made up of line segments. The vertices of $\Delta$ are isolated, and a small neighbourhood of a vertex may be represented by the intersection of two unique half-planes.

**Definition 3.3.** For $\lambda \in \mathbb{Z}^n$ write the **greatest common divisor** by

$$\gcd(\lambda) = \max\left\{ a \in \mathbb{N}_{>0} \;\middle|\; \frac{1}{a}\lambda \in \mathbb{Z}^n \right\} .$$

*Remark* 3.4. For all $A \in \mathrm{GL}_n(\mathbb{Z})$, $\gcd(A\lambda) = \gcd(\lambda)$.

**Definition 3.5.** A vertex $x$ of a polygonal domain $\Delta$ is of **type** $A_n$ if the two adjacent half-planes $H_{\lambda_1,c_1}, H_{\lambda_2,c_2}$ satisfy:

$$|\det(\lambda_1, \lambda_2)| = n + 1$$
$$\gcd(\lambda_1 - \lambda_2) = n + 1 ,$$

with $n \in \mathbb{Z}_{\geq 0}$. We call $x$ **weakly Delzant** if there exists $n \in \mathbb{Z}_{\geq 0}$ such that $x$ is of type $A_n$.

We say that $\Delta$ is **weakly Delzant** if all vertices of $\Delta$ are weakly Delzant. We say that $\Delta$ is **Delzant** if all vertices of $\Delta$ are of type $A_0$.

In [MS23], weakly Delzant polygonal domains are called *canonical tropical domains.*

*Remark* 3.6. A polygonal domain is an integral affine surface (with corners) if and only if it is Delzant.

*Remark* 3.7. If $\Delta$ is weakly Delzant, there is a unique symplectic toric orbifold with moment map $\mu : \tilde{X} \to \Delta$. An $A_n$-type vertex of $\Delta$ is the image of an $A_n$-singular point of $\mu$ in $\tilde{X}$, see [Eva23, Section 3.4] for details.

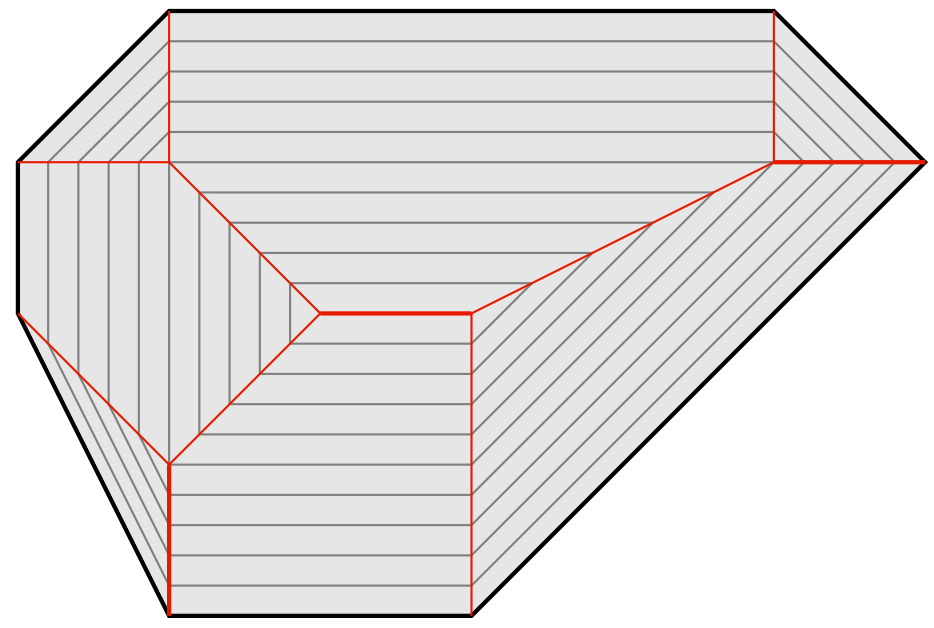

Figure 12: Level sets of $\mathscr{F}_\Delta$ in dark grey and the caustic $\mathscr{K}_\Delta$ in red.

Let $\Delta \subset \mathbb{R}^2$ be a weakly Delzant polygonal domain, and write $\Delta = \bigcap_{H_{\lambda,c} \in \mathscr{H}} H_{\lambda,c}$, where $\mathscr{H}$ is a minimal set of half-planes. Let

$$\begin{aligned} \mathscr{F}_\Delta : \mathbb{R}^2 &\to \mathbb{R} \cup \{-\infty\} \\ p &\mapsto \inf_{H_{\lambda,c} \in \mathscr{H}} \{\langle \lambda, p \rangle + c\}, \end{aligned}$$

be the **height function** of $\Delta$. We can write $\Delta = \mathscr{F}_\Delta^{-1}(\{s \geq 0\})$. For $t \in \mathbb{R}$, set $\Delta_t = \mathscr{F}_\Delta^{-1}(\{s \geq t\}) = \bigcap_{H_{\lambda,c} \in \mathscr{H}} H_{\lambda,c+t}$. See Figure 12.

The value $-\infty$ is never attained by $\mathscr{F}_\Delta$, as by Definition 3.1 only finitely many half-planes contribute to the infimum.

**Definition 3.8.** The **caustic $\mathscr{K}_\Delta$ of** $\Delta$ is the set of points of $\Delta$ where $\mathscr{F}_\Delta$ is not affine. See Figure 12.

*Remark* 3.9. The caustic $\mathscr{K}_\Delta$ is a tropical curve: By Definition 3.1, the restriction of $\mathscr{F}_\Delta$ to any compact set is a tropical polynomial. See e.g. [BS14, Section 2.1] for the definition of a tropical curve.

**Theorem 3.10** ([MS23, Theorem 26])**.** *Let $\Delta$ be weakly Delzant. Then $\Delta_t$ is a weakly Delzant polygonal domain for all $0 \leq t < \sup \mathscr{F}_\Delta$.*

**Definition 3.11.** We say that an edge $\ell$ of $\mathscr{K}_\Delta$ has **weight** $n$ if for any $x$ in the interior of $\ell$, $x$ is an $A_{n-1}$-type vertex of the weakly Delzant polygonal domain $\Delta_{\mathscr{F}_\Delta(x)}$.

*Remark* 3.12. For edges where the weight is defined (namely those not contained in the max-set of $\mathscr{F}_\Delta$), this is equal to the weight of an edge of a tropical curve given by a tropical polynomial, see for example [BS14, Definition 2.1].

The following theorem describes the local behaviour of $\mathscr{F}_\Delta$ near edges and vertices of $\mathscr{K}_\Delta$. For edges, this is a reformulation of Theorem 3.10, while the statement for vertices is a reformulation of [MS23, Theorem 46].

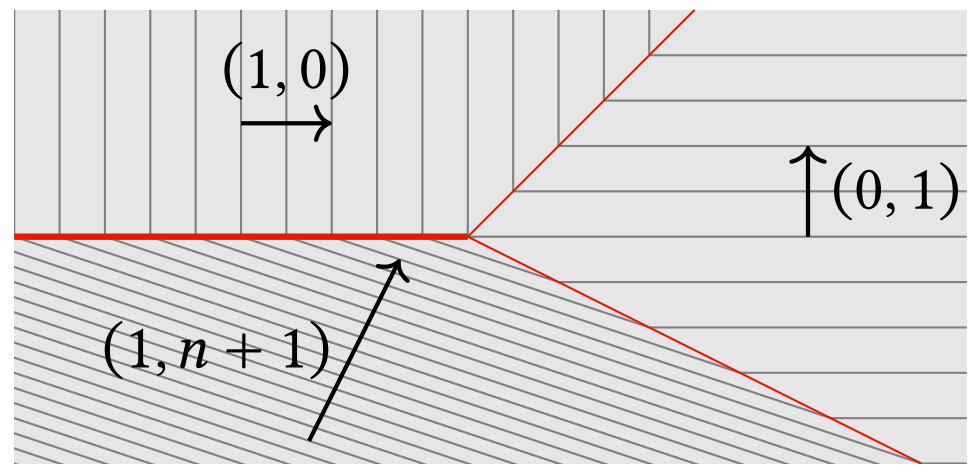

Figure 13: Level sets of $\mathscr{F}_\Delta$ near a vertex of $\mathscr{K}_\Delta$ outside of $\Delta_M$.

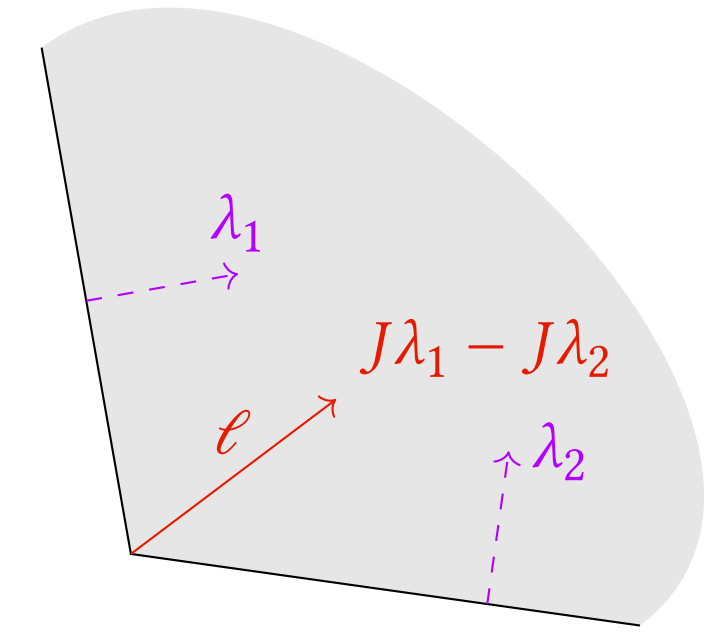

Figure 14: Smoothing an $A_n$ corner.

**Theorem 3.13** ([MS23, Theorems 26, 46])**.** *Let $\Delta$ be weakly Delzant and $M = \sup \mathscr{F}_\Delta$. If $\ell$ is an edge of $\mathscr{K}_\Delta$ in $\Delta \setminus \Delta_M$ of weight $n+1$, $x$ is in the interior of $\ell$ and $a \in \mathbb{R}^2$ small enough, then, after an integral affine transformation, we can write*

$$\mathscr{F}_\Delta(x+a) = \mathscr{F}_\Delta(x) + \min\{\langle(1,0),a\rangle, \langle(1,n+1),a\rangle\} \, .$$

*For a vertex $x$ of $\mathscr{K}_\Delta$ in $\Delta \setminus \Delta_M$ and $a \in \mathbb{R}^2$ small enough, after an integral affine transformation, we can write*

$$\mathscr{F}_\Delta(x+a) = \mathscr{F}_\Delta(x) + \min\left\{\langle(1,0),a\rangle \, , \langle(1,n+1),a\rangle \, , \langle(0,1),a\rangle\right\}$$

*with $n \in \mathbb{Z}_{\geq 0}$.*

The case of a vertex is illustrated in Figure 13.

**Corollary 3.14.** *If $\Delta$ is a Delzant polygonal domain, then so is $\Delta_t$ for all $0 \leq t < \sup \mathscr{F}_\Delta$ and in Theorem 3.13 we have $n = 0$ for all vertices and edges of $\mathscr{K}_\Delta$.*

**Definition 3.15.** We can **smooth** an $A_n$ type vertex $x$ of a weakly Delzant polygonal domain by decorating it as follows: Let $H_{\lambda_1,c_1}, H_{\lambda_2,c_2}$ be the two incident half planes at $x$ with $\det(\lambda_1,\lambda_2) = n+1$, and let $J$ be the rotation by $\pi/2$. We add a small cut $\ell$ from $x$ in direction $J(\lambda_1 - \lambda_2)$ ending in a node $\mathfrak{n}$ of multiplicity $n+1$.

By interpreting this as a nodal chart diagram, this makes $\Delta$ into a nodal integral affine surface near $x$.

Let $\Delta$ be a weakly Delzant polygonal domain with exactly one vertex $x$. Smooth $x$ to obtain a nodal integral affine surface $(B, \{\mathfrak{n}\})$ with a simple nodal chart $\varphi : B \to \Delta$ with nodal chart diagram given by Definition 3.15. Since $B$ is contractible, by Theorem 2.27, it determines a unique symplectic manifold $X$ with an almost toric fibration $\pi : X \to B$. We also can obtain $X$ by smoothing the $A_n$ singularity corresponding to the vertex $x$ of the toric orbifold $\mu : \tilde{X} \to \Delta$ determined by Remark 3.7 as in [Eva23, Section 7.3]. For a neighbourhood $U$ of the cut $\ell$, we can also require that the fibrations $\varphi \circ \pi : X \setminus (\varphi \circ \pi)^{-1}(U) \to \Delta$ and $\mu : \tilde{X} \setminus \mu^{-1}(U) \to \Delta$ are identical. If $n = 0$, the $A_0$ singularity is smooth already and the smoothing corresponds to a *nodal trade* (see [Sym03, Lemma 6.3] or [Eva23, Section 8.2]), and we have $X = \tilde{X}$.

*Remark* 3.16. The cut $\ell$ obtained by smoothing a vertex $x$ lies in the caustic $\mathcal{K}_\Delta$: Near $x$ the caustic is given by the equation

$$\langle \lambda_1, \cdot \rangle + c_1 = \langle \lambda_2, \cdot \rangle + c_2 \Leftrightarrow \langle \lambda_1 - \lambda_2, \cdot \rangle + c_1 - c_2 = 0$$

and by construction $\ell$ is contained in this set. In particular, by modifying the nodal integral affine surface $B$ by small nodal slides of $\mathfrak{n}$, the node stays on the caustic.

**Proposition 3.17.** *Let $(B_0, \mathfrak{N}_0)$ be the nodal integral affine surface obtained by smoothing all vertices of a weakly Delzant polygonal domain $\Delta$, and let $\varphi_0 : B_0 \to \Delta$ be the associated nodal chart.*

*For any $0 < M_\varepsilon < M := \sup \mathcal{F}_\Delta$ there exists a nodal tangle $(B \times [0,1], \mathfrak{N})$ starting at $(B_0, \mathfrak{N}_0)$ such that $\mathcal{F}_1 := \mathcal{F}_\Delta \circ \varphi_1$ is an integral affine function on $B_1 \setminus (\varphi_1^{-1}(\Delta_{M_\varepsilon}) \cup \mathfrak{N}_1)$, where $\varphi_1 := \varphi_0 \circ \tau_1^0$.*

In applications we take $M_\varepsilon$ to be arbitrarily close to $M$.

*Proof.* To construct the nodal tangle $(B \times [0,1], \mathfrak{N})$ we proceed as illustrated in Figure 15 b): By Remark 3.16 all nodes introduced by the smoothing slide on the caustic. Sort the nodes by decreasing multiplicity. Start with the node $\mathfrak{n}_1$ and slide it up (with respect to the height function $\mathcal{F}_\Delta$) until it slides off the caustic at a vertex of $\mathcal{K}_\Delta$, or it is contained in $\Delta_{M_\varepsilon}$. If it slides off $\mathcal{K}_\Delta$ at a vertex we stop it shortly after the vertex and call it **parked**. Now continue with the next node $\mathfrak{n}_2$. If it reaches a vertex $x$ of $\mathcal{K}_\Delta$ where a node previously slid off $\mathcal{K}_\Delta$, Theorem 3.13 assures that the multiplicity of $\mathfrak{n}_2$ is 1 and, using Remark 2.20, $\mathfrak{n}_2$ stays on $\mathcal{K}_\Delta$ sliding through $x$. See Figure 16. Continuing this procedure, we see inductively that every node either becomes parked or stays on the caustic and eventually ends up in $\Delta_{M_\varepsilon}$. This procedure provides a nodal tangle $(B \times [0,1], \mathfrak{N})$.

Now we show that $\mathcal{F}_1 : (B_1, \mathfrak{N}_1) \to \mathbb{R}$ is integral affine on $B_1 \setminus (\varphi_1^{-1}(\Delta_{M_\varepsilon}) \cup \mathfrak{N}_1)$.

Let $G$ be the cut graph of $\varphi_1^{-1}$ in $\Delta \setminus \Delta_{M_\varepsilon}$. By construction of $(B_1, \mathfrak{N}_1)$ we have $\mathcal{K}_\Delta \cap (\Delta \setminus \Delta_{M_\varepsilon}) \subset G$ and that $\mathcal{F}_1$ is integral affine on $\varphi_1^{-1}(\Delta \setminus (\Delta_{M_\varepsilon} \cup G))$.

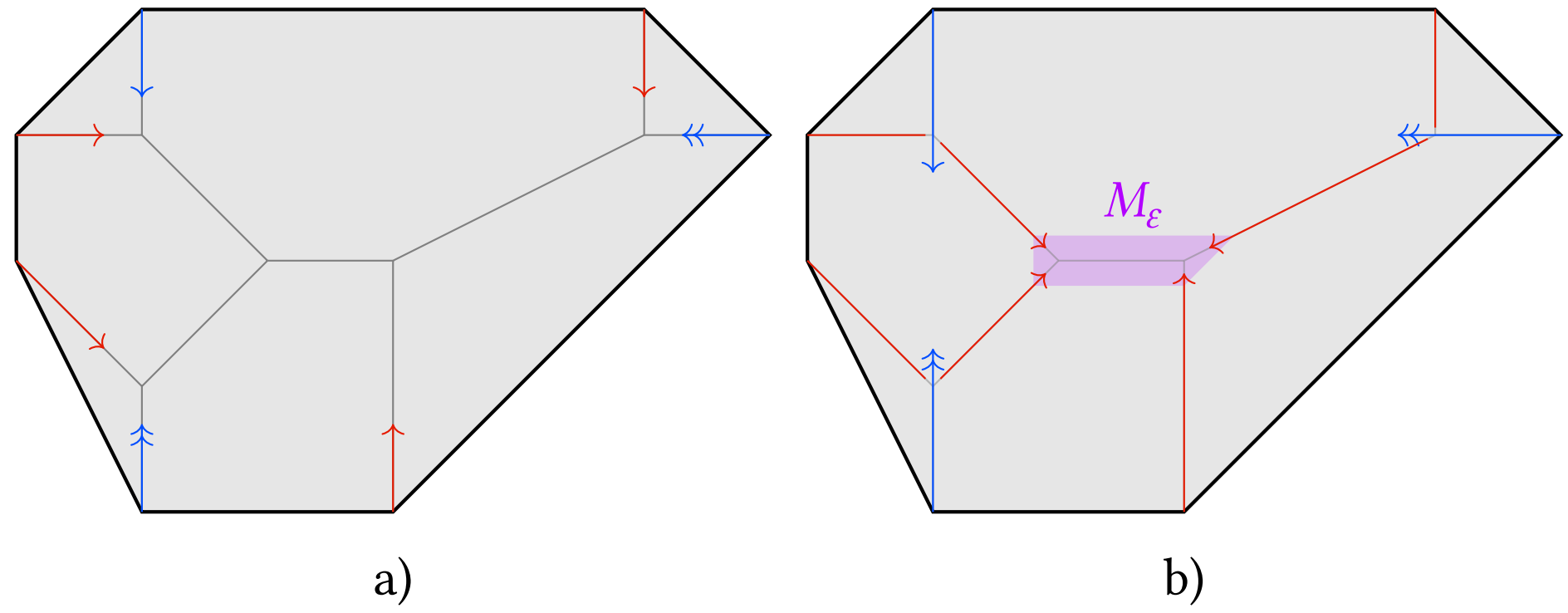

Figure 15: a): smoothing a weakly Delzant polygon $\Delta$.
b): filling the caustic $\mathscr{K}_\Delta$ with cuts.
Nodes of multiplicity 2 have a double arrow head. Nodes parked in b) and their cuts are drawn in blue.

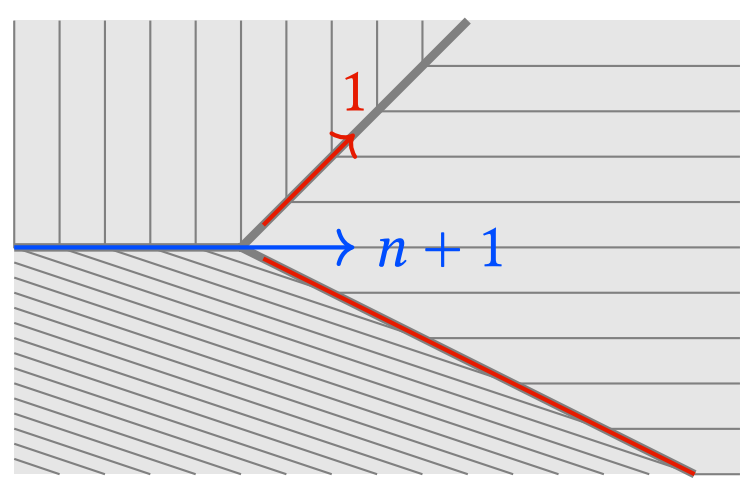

Figure 16: Filling $\mathscr{K}_\Delta$ with cuts near a vertex of $\mathscr{K}_\Delta$.

Let $\varphi_1(x) \in G \setminus \varphi_1(\mathfrak{N}_1)$. Then $\varphi_1(x)$ is either in the interior of an edge $\ell$ of $\mathscr{K}_\Delta$, a vertex of $\mathscr{K}_\Delta$ or in the interior of a cut of $G$ connecting a vertex of $\mathscr{K}_\Delta$ to a parked node.

Suppose $\varphi_1(x) \in \mathscr{K}_\Delta$ is a vertex of $\mathscr{K}_\Delta$ and, after a translation of $\Delta$, $\varphi_1(x) = 0$. See Figure 16. After an integral affine transformation, using Theorem 3.13 and the construction of the nodal tangle $(B \times [0,1], \mathfrak{N})$,

$$\mathscr{F}_\Delta(\underbrace{\varphi_1(x)}_{=0} + a) = \mathscr{F}_\Delta(a) = \mathscr{F}_1(x) + \min\{\langle(1,0),a\rangle, \langle(1,n+1),a\rangle, \langle(0,1),a\rangle\} \quad (2)$$

$$R_0 = h_{(1,n+1)} \circ h_{(1,0)}^{n+1} ,$$

where $R_0$ is the rectifying map at 0 of $\varphi_1$, and $h_v$ denotes a half-shear as in (1) Taking the integral affine chart $R_0 \circ \varphi_1$, using (2) and applying the half-shears of $R_0$ to Figure 16, it is quickly verified that for small $a \in \mathbb{R}^2$,

$$(\underbrace{\mathscr{F}_\Delta \circ \varphi_1}_{=\mathscr{F}_1} \circ (R_0 \circ \varphi_1)^{-1})(a) = \mathscr{F}_\Delta(R_0^{-1}a) = \mathscr{F}_1(x) + \langle(1,0),a\rangle ,$$

showing that $\mathscr{F}_1 = \mathscr{F}_\Delta \circ \varphi_1$ is integral affine near $x$.

We proceed similarly in the remaining two cases. □

*Remark* 3.18. We may choose a different nodal chart for $(B_1, \mathfrak{N}_1)$ of Proposition 3.17 in which $\mathscr{F}_1$ is visibly integral affine:

Take the annulus $A = \varphi_1^{-1}(\Delta \setminus \Delta_M) \subset (B_1, \mathfrak{N}_1)$, and remove a straight line $\ell$ such that $U = A \setminus \ell$ is simply connected. On $U$ we may choose a nodal chart $\varphi'$ as in Figure 17 as follows: The level sets $\mathscr{F}_1^{-1}(h)$ with $h \leq M_\varepsilon$ are all straight lines. For $x \in B_1$ and $h = \mathscr{F}_1(x)$, let $f(x)$ be the integral affine length of $\mathscr{F}_\Delta^{-1}(h)$ between $\mathscr{F}_\Delta^{-1}(h) \cap \ell$ and $x$ in the counter-clockwise direction. Then $\varphi' = (f, \mathscr{F}_1)$ is a nodal chart of $B_1 \setminus (\ell \cup \varphi_0^{-1}(\Delta_M))$. See Figure 17.

In this nodal chart the parked nodes have horizontal cuts. Hence $\mathscr{F}_1$ is indeed integral affine on $A \setminus \mathfrak{N}_1$.

## 3.2 Toric actions on symplectic four manifolds

Theorem 3.13 describes vertices $x$ of $\mathscr{K}_\Delta$ with $\mathscr{F}_\Delta(x) < \sup \mathscr{F}_\Delta$. [MS23, Propositions 54, 55] describe vertices $x$ of $\mathscr{K}_\Delta$ with $\mathscr{F}_\Delta(x) = \sup \mathscr{F}_\Delta$. Here we restate the case of Delzant polygonal domains:

**Proposition 3.19.** *Let $\Delta$ be a Delzant polygonal domain. Then near a vertex $x$ of $\mathscr{K}_\Delta$ with $\mathscr{F}_\Delta(x) = \sup \mathscr{F}_\Delta$, up to an integral affine transformation, $\mathscr{F}_\Delta$ is given by one of the seven possibilities in Figure 18.*

Let $\Delta$ be a compact Delzant polygonal domain. Since then $M = \sup \mathscr{F}_\Delta < \infty$, we may choose $M_\varepsilon = M - \varepsilon$ such that $\Delta_{M_\varepsilon} \setminus \Delta_M$ contains no vertices of $\mathscr{K}_\Delta$.

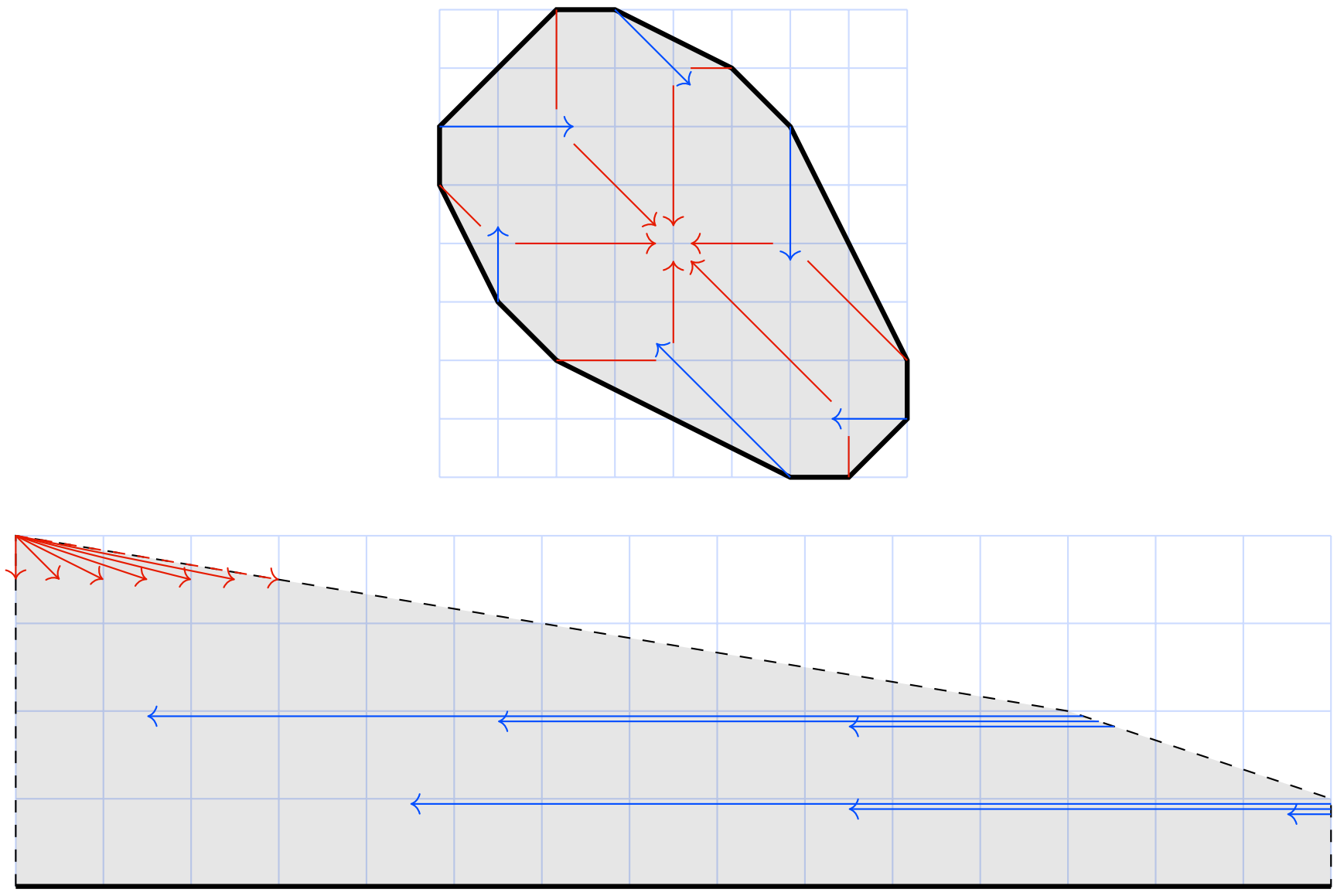

Figure 17: Constructing an “open” nodal chart. Parked nodes are drawn in blue.

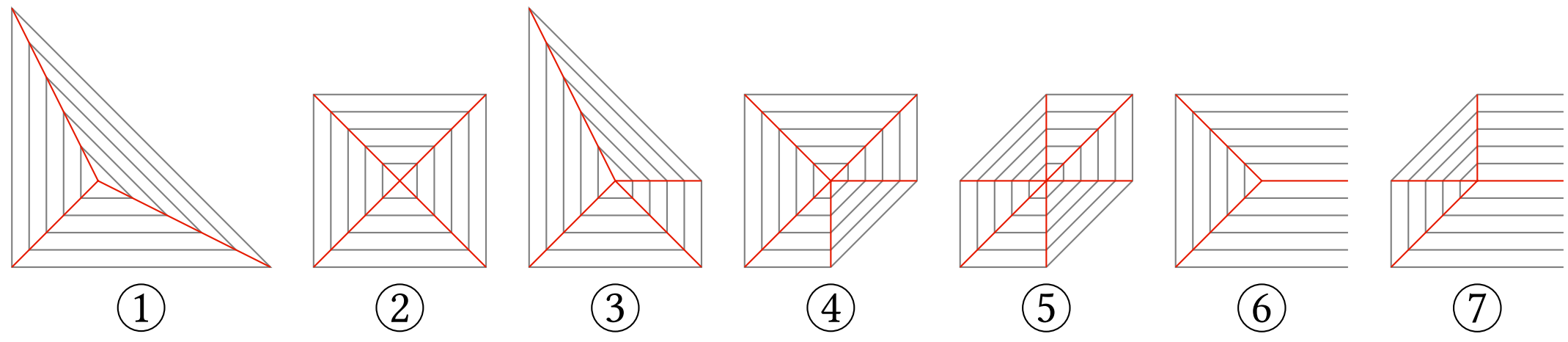

Figure 18: Possible caustics near $\sup \mathscr{F}_\Delta$ for Delzant polygonal domains. $\mathscr{K}_\Delta$ in red, level sets of $\mathscr{F}_\Delta$ in grey.

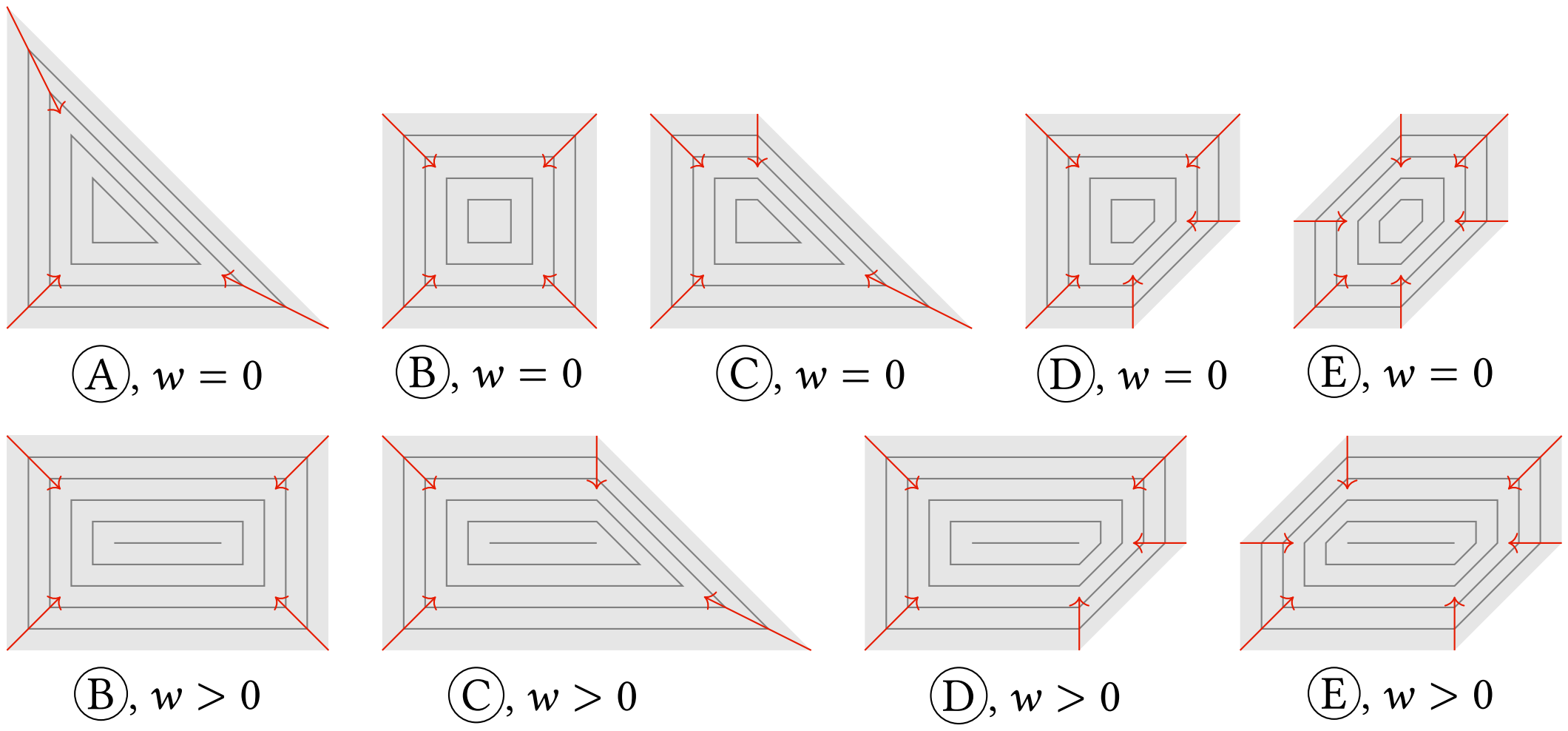

Figure 19: All possible $\varepsilon$-hats up to nodal tangle.

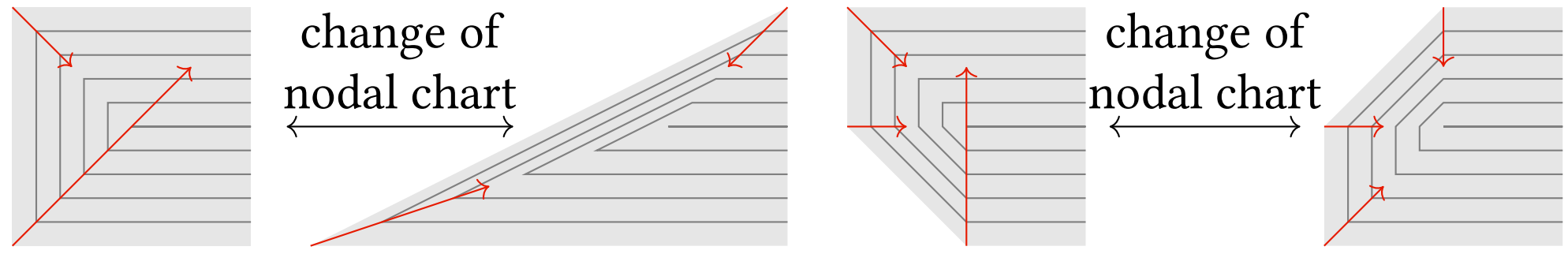

Figure 20: Nodal tangle transforming part of an $\varepsilon$-hat of positive width.

**Definition 3.20.** In this case, we say $B_1$ as constructed in Proposition 3.17 is in **canonical form** for $\Delta$. The $\varepsilon$-**hat** of $B_1$ is given by the nodal integral affine surface

$$H_\varepsilon = \{x \in B_1 \mid \mathscr{F}_1(x) > M - \varepsilon\} .$$

The **hat width** is given by the integral affine length of $\Delta_M$.

**Definition 3.21.** An $\varepsilon_1$-hat $H_1$ is equivalent to an $\varepsilon_2$-hat $H_2$ if, after a nodal tangle, $H_1 \subset H_2$ or $H_2 \subset H_1$. We call the equivalence class of an $\varepsilon$-hat its **hat class**.

**Lemma 3.22.** *Let $w \geq 0$.*

- *If $w = 0$, there are five hat classes of width $w$,*
- *If $w > 0$, there are four hat classes of width $w$.*

*They are listed in Figure 19.*

*Proof.* If $w = 0$, this follows immediately from Proposition 3.19, with the five hats Ⓐ–Ⓔ in the top row of Figure 19 corresponding to the caustics ①–⑤ in Figure 18.

If $w > 0$, then the $\varepsilon$-hat must be assembled from ⑥ and ⑦ in Figure 18, and by using the moves illustrated in Figure 20 it is easy to check that there exists a nodal tangle relating it to one of the $\varepsilon$-hats in the bottom row of Figure 19. □

**Definition 3.23.** The **canonical type** of a compact Delzant polygonal domain $\Delta$ is given by $(H, M, \alpha_1, \ldots, \alpha_n)$, where

- $H$ is the hat class of an $\varepsilon$-hat of $\Delta$,
- $M = \sup \mathcal{F}_\Delta$,
- $\alpha_i = \mathcal{F}_1(\mathfrak{a}_i)$ where $\mathfrak{a}_1, \ldots, \mathfrak{a}_n$ are the parked nodes of $(B_1, \mathfrak{N}_1)$ such that $\alpha_1 \geq \alpha_2 \geq \ldots \geq \alpha_n$.

We can essentially reconstruct a canonical form of $\Delta$ from its canonical type up to nodal tangle. In Figure 17 the $\alpha_i$'s correspond to the heights of the parked nodes.

**Definition 3.24** ([KK17, Definition 1.2])**.** Let $k \in \mathbb{Z}_{\geq 0}$. A vector $v = (\lambda, \delta_1, \ldots, \delta_k) \in \mathbb{R}^{k+1}$ is **reduced** if

$$\delta_1 \geq \delta_2 \geq \ldots \geq \delta_k, \quad \delta_1 + \delta_2 + \delta_3 \leq \lambda \ .$$

If additionally

$$\delta_1 \geq \delta_2 \geq \ldots \geq \delta_k > 0 \quad \text{and} \quad \lambda^2 - \sum_{i=1}^{k} \delta_i^2 > 0 \ ,$$

we call it **symplectic**.

Let $X_k = \mathbb{C}P^2 \# k \overline{\mathbb{C}P^2}$. In [KK17] it is shown that symplectic reduced vectors classify symplectic forms $\omega$ on $X_k$ up to diffeomorphism. In a symplectic reduced vector, $\lambda$ encodes the symplectic size of $\mathbb{C}P^2$, and the $\delta_i$ encode the symplectic sizes of the blow-ups used to obtain $X_k$ from $\mathbb{C}P^2$. For a reduced vector $v$ denote by $\omega_v$ the corresponding symplectic form.

We call a symplectic reduced vector $v$ **toric** if the corresponding $(X, \omega_v)$ admits a toric moment map.

Let $\mathcal{R}$ be the set of all toric reduced vectors, and let $\mathcal{C}$ be the set of all canonical types obtained from compact Delzant polygonal domains.

**Theorem 3.25.** *There is a bijection $\alpha : \mathcal{C} \to \mathcal{R}$ such that if $\mu : (X, \omega_v) \to \Delta$ is a toric moment map, the canonical type of $\Delta$ is sent to $v$.*

*In particular, every four-manifold admitting a toric moment map has a unique canonical type.*

*Proof.* Take $C = (H, M, \alpha_1, \ldots, \alpha_n) \in \mathcal{C}$ and let $w$ be the width of $H$. Depending on the hat class $H$ we define $\alpha(C)$ as follows:

| $H$ | $\lambda,$ | $\delta_1, \ldots, \delta_k$ |
|---|---|---|
| Ⓐ | $3M,$ | $\alpha_1, \ldots, \alpha_n$ |
| Ⓑ | $4M + w - \alpha_1,$ | $(2M + w - \alpha_1), (2M - \alpha_1), \alpha_2, \ldots, \alpha_n$ |
| Ⓒ | $3M + w,$ | $M + w, \alpha_1, \ldots, \alpha_n$ |
| Ⓓ | $3M + w,$ | $M + w, M, \alpha_1, \ldots, \alpha_n$ |
| Ⓔ | $3M + w,$ | $M + w, M, M, \alpha_1, \ldots, \alpha_n$ |

(We used the same letters as in Figure 19 to denote the type of hat.) The rest of the proof will motivate this definition. It is easy to check that in each case the vector $\alpha(C)$ is reduced and symplectic: The inequalities $\delta_1 \geq \ldots \geq \delta_k > 0$ and $\delta_1+\delta_2+\delta_3 \leq \lambda$ are quickly verified using $M > \alpha_1 \geq \alpha_2 \geq \ldots \geq \alpha_n$ and $w \geq 0$. Moreover, the quantity $\lambda^2 - \sum_{i=1}^{k} \delta_i^2$ is equal to twice the area of $B$, which is positive. This can be most easily verified by examining a nodal chart diagram as in Remark 3.18, which we do for the hat type Ⓔ, as pictured in Figure 17 in the case $w = 0$. For the case Ⓔ we have

$$\lambda^2 - \sum_{i=1}^{k} \delta_i^2 = (3M + w)^2 - (M + w)^2 - 2M^2 - \sum_{i=1}^{n} \alpha_i^2 = 6M^2 + 4Mw - \sum_{i=1}^{n} \alpha_i^2 . \quad (3)$$

If there were no parked nodes, Figure 17 would be the union of a triangle of height $M$ and width $6M$ and a rectangle of height $M$ and width $2w$. Every time we introduce a parked node we remove $\alpha_i^2/2$ area, so twice the area of $B$ is given by the right hand side of (3).

Let $(B, \mathfrak{N})$ be the canonical form given by $C$ and $\pi : (Y, \omega) \to (B, \mathfrak{N})$ an almost toric fibration. Our goal is to show that $(Y, \omega)$ is symplectomorphic to $(X, \omega_v)$ for some reduced vector $v$ uniquely determined by $C$. We proceed as follows: To start, suppose for simplicity that $H$ is not of type Ⓑ. First, we blow down $(Y, \omega)$ symplectically with sizes $\alpha_n, \ldots, \alpha_1$. This can be done either as *almost toric blow downs* ([Eva23, Section 9.1]) at the parked nodes $\mathfrak{a}_n, \ldots, \mathfrak{a}_1$, or as blow downs along the exceptional *visible symplectic spheres* obtained by connecting $\mathfrak{a}_n, \ldots, \mathfrak{a}_1$ to $\partial B$ by smooth paths transversal to the height function $\mathscr{F}_1$. In this manner we arrive at an almost toric base diagram without parked nodes. Pushing the remaining nodes to the boundary, we get a toric base diagram. Depending on $H$, the resulting symplectic space is

| Ⓐ | Ⓒ | Ⓓ | Ⓔ |
|---|---|---|---|
| $\mathbb{C}P^2$ | $\mathbb{C}P^2\#\overline{\mathbb{C}P^2}$ | $\mathbb{C}P^2\#2\overline{\mathbb{C}P^2}$ | $\mathbb{C}P^2\#3\overline{\mathbb{C}P^2}$ |

Blowing down up to three more times with appropriate sizes (see the table above), we get $\mathbb{C}P^2$. The symplectic exceptional spheres along which we blow down can also be seen as visible symplectic curves in $(B, \mathfrak{N})$, see Figure 21 for the hat type Ⓔ. In total, we see that all symplectic spheres along which we blow down can be chosen to be disjoint — indeed, they are all visible in $B$ as disjoint visible (tropical) symplectic curves — and the sequence of blow down sizes produces a symplectic reduced vector $v = \alpha(C)$ as in the table above. Thus $(V, \omega) = (X, \omega_v)$.

For the case of a hat $H$ of type Ⓑ, we only blow down $\alpha_n, \ldots, \alpha_2$, arriving at $\mathbb{C}P^2\#2\overline{\mathbb{C}P^2}$. Blowing down two more times at visible symplectic spheres as in Figure 21, we get $\mathbb{C}P^2$. (Here we need to avoid blowing down $\alpha_1$, as otherwise we get stuck at $S^2 \times S^2$.)

The size $\lambda$ of $\mathbb{C}P^2$ thusly obtained can be calculated from $H$ and $M$ (and $\alpha_1$ in the case of a hat of type Ⓑ).

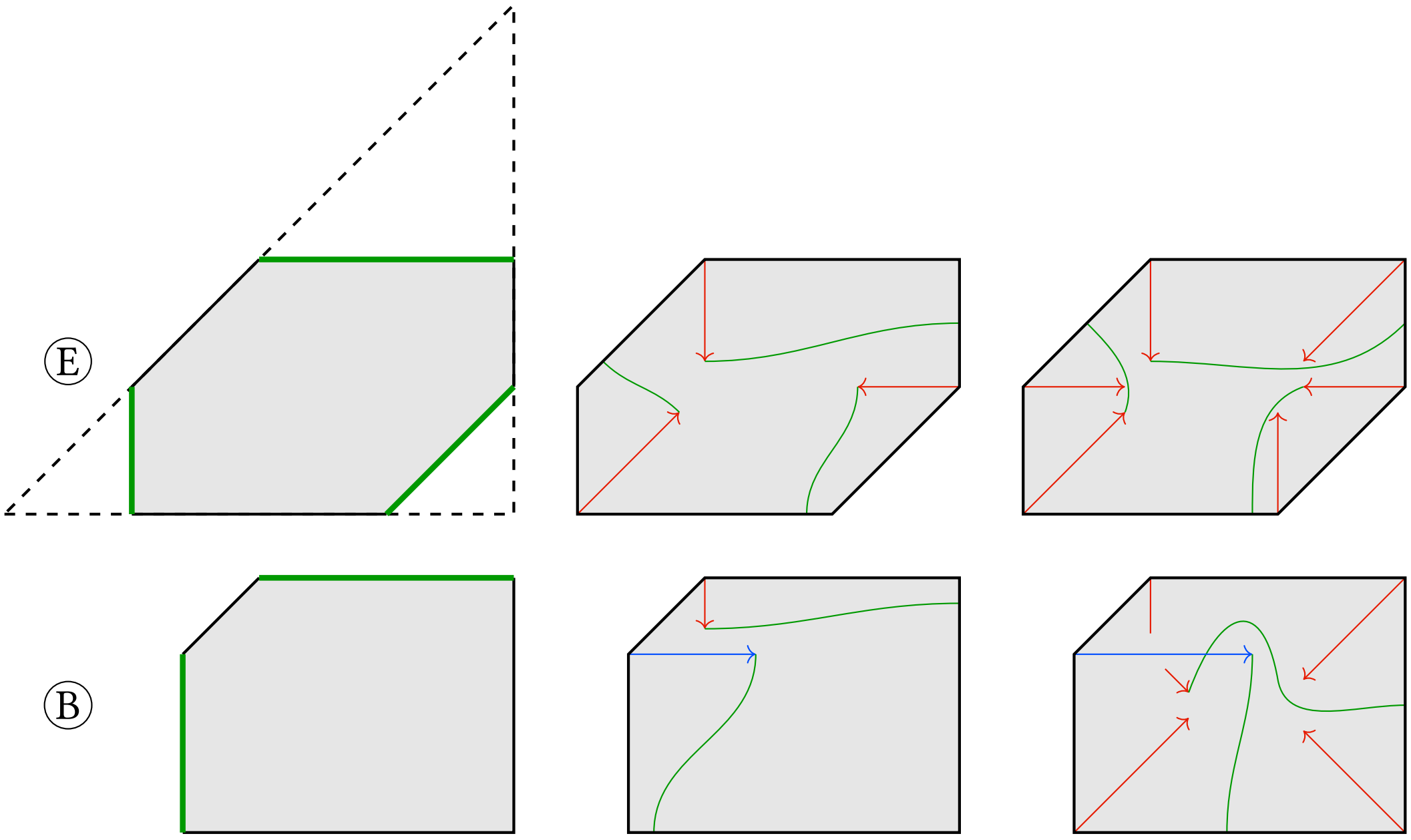

Figure 21: Constructing visible tropical spheres for hats of type Ⓔ (top) and Ⓑ (bottom). From left to right: in a toric model we first mark which spheres we would like to blow down in green. Then we perform nodal trades to obtain the middle diagram where the marked spheres are visible symplectic curves. Performing the remaining nodal trades and sliding the nodes to get the canonical form, we can read of the desired visible symplectic spheres.

To show that $\alpha$ gives a one-to-one correspondence, we construct an inverse $\beta : \mathscr{R} \to \mathscr{C}$. Given $\lambda, \delta_1, \ldots, \delta_k$, set

$$M = \min\left\{\frac{\lambda}{3}, \frac{\lambda - \delta_1}{2}\right\} = \begin{cases} \frac{\lambda}{3} & \text{if } \frac{\lambda}{3} \geq \delta_1 \\ \frac{\lambda-\delta_1}{2} & \text{if } \delta_1 \geq \frac{\lambda}{3} \end{cases}.$$

Depending on how many $\delta_i$ are bigger or equal than $M$, we define the canonical type $\beta(\lambda, \delta_1, \ldots, \delta_k) = (H, M, \alpha_1, \ldots, \alpha_n)$ as follows:

| | $H$ | $\alpha_1, \ldots, \alpha_n$ |
|---|---|---|
| $M > \delta_1$ | Ⓐ, $w = 0$ | $\delta_1, \ldots, \delta_k$ |
| $\delta_1 \geq \delta_2 > M > \delta_3$ | Ⓑ, $w = \delta_1 - \delta_2$ | $(\lambda - \delta_1 - \delta_2), \delta_3, \ldots, \delta_k$ |
| $\delta_1 \geq M > \delta_2$ | Ⓒ, $w = \frac{3\delta_1 - \lambda}{2}$ | $\delta_2, \ldots, \delta_k$ |
| $\delta_1 \geq \delta_2 = M > \delta_3$ | Ⓓ, $w = \frac{3\delta_1 - \lambda}{2}$ | $\delta_3, \ldots, \delta_k$ |
| $\delta_1 \geq \delta_2 = \delta_3 = M > \delta_4$ | Ⓔ, $w = \frac{3\delta_1 - \lambda}{2}$ | $\delta_4, \ldots, \delta_k$ |

(This table can be constructed by following the above process in reverse, starting with $\mathbb{C}P^2$ and performing toric blow-ups for all $\delta_i \geq M$, then reading off the hat type.) It is easily checked that $\beta \circ \alpha = \mathrm{id}_{\mathscr{C}}$ and $\alpha \circ \beta = \mathrm{id}_{\mathscr{R}}$. □

**Corollary 3.26.** *Let $\mu_i : X \to \Delta_i$ with $i \in \{1, 2\}$ be two toric moment maps on a four dimensional closed symplectic manifold. Then smoothings (Definition 3.15) of $\Delta_1$ and $\Delta_2$ are connected by a nodal tangle.*

*Proof.* Since the canonical type is uniquely determined by $X$ (Theorem 3.25), both $\Delta_1$ and $\Delta_2$ have the same canonical form, and are thus their smoothings are connected by a nodal tangle. □

See Section 5.2 for possible generalizations.

## 3.3 Lagrangian Poincaré non-recurrence

The goal of this section is to prove Theorem B. That is, for a closed non-monotone toric symplectic four manifold $(X, \omega)$, we want to construct a Lagrangian $L \subset (X, \omega)$ and Hamiltonian diffeomorphism $\psi$ such that $\psi^n(L) \cap L = \emptyset$ for all $n \in \mathbb{Z}_{>0}$.

*Proof of Theorem B.* We give the proof in three parts:

**The base space.** Let $(B_0, \mathfrak{N}_0)$ be a nodal integral affine surface for $(X, \omega)$ as in Proposition 3.17 with canonical type $(H, M, \alpha_1, \ldots, \alpha_n)$ obtained by modifying a Delzant polygon by a nodal tangle. Recall that $(B_0, \mathfrak{N}_0)$ admits a natural nodal chart $\varphi_0 : B_0 \to \mathbb{R}^2$ as in Remark 3.18, illustrated in Figure 22 a). Let $\mathscr{F} : B_0 \to [0, M)$ be the height function on $B_0$ inherited from the Delzant polygon and let $\varepsilon > 0$. We

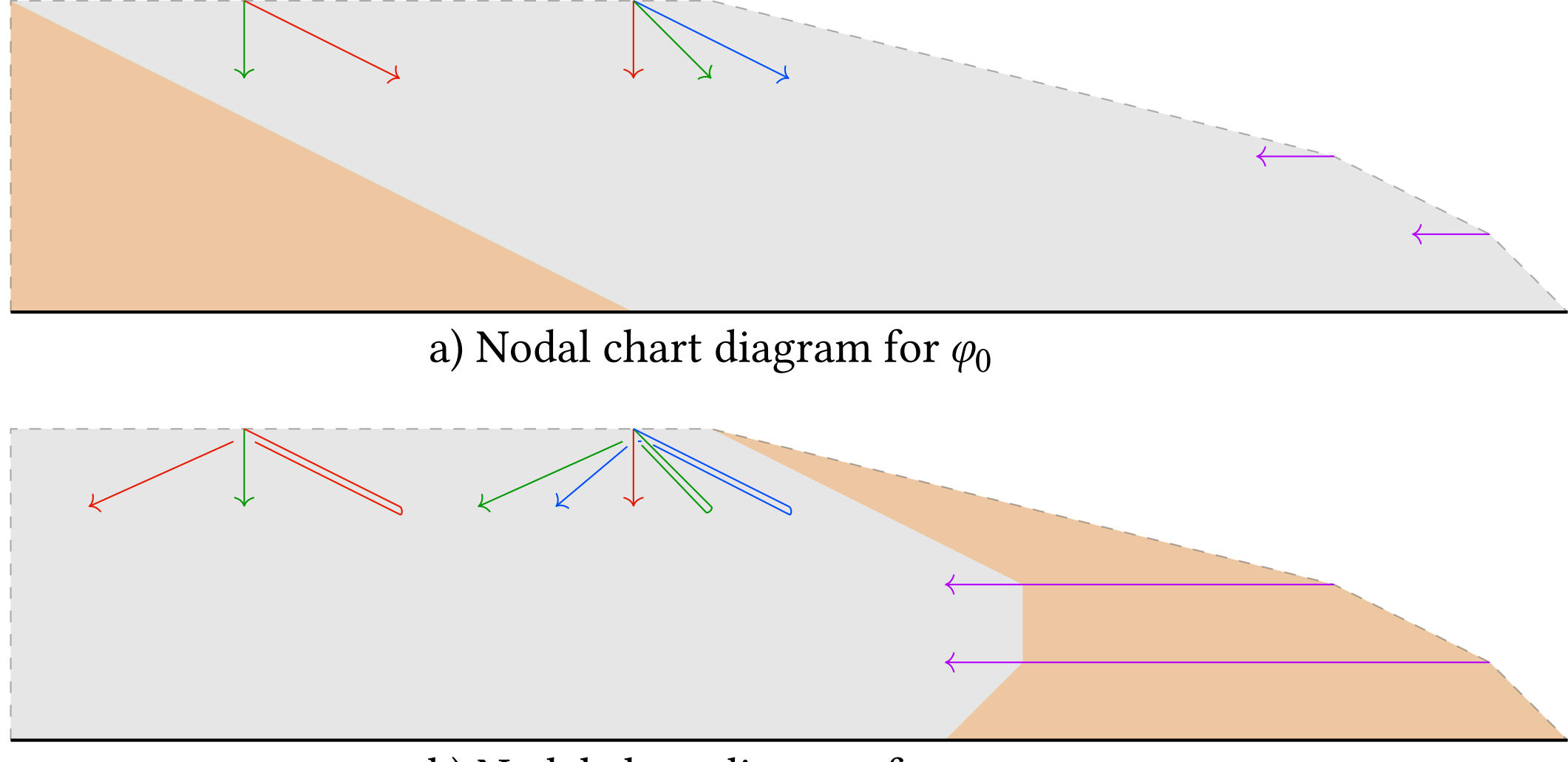

a) Nodal chart diagram for $\varphi_0$

b) Nodal chart diagram for $\varphi_1$

Figure 22: Two nodal charts for the same canonical type. (Slopes are not to scale.)

may choose the nodes in the hat $H$ to have at least height $M - \varepsilon$. The level sets $\mathscr{F}^{-1}(h)$ with $h \in [0, M - \varepsilon)$ are closed straight lines in $B_0$. We have $M = \max \mathscr{F}$, and we call $\Delta_M := \mathscr{F}^{-1}(M)$ the **ridge** of $B_0$. The integral affine length of $\Delta_M$ is the width $w$ of the hat $H$ (Definition 3.20). The integral affine length of $\mathscr{F}^{-1}(h)$ is given by the function

$$g(h) = 2w + k(M - h) + \sum_{i=1}^{n} \min\{h - \alpha_i, 0\},$$

where $k$ is a constant depending on the hat class of $H$ as follows:

| Ⓐ | Ⓑ | Ⓒ | Ⓓ | Ⓔ |
|---|---|---|---|---|
| $k = 9$ | $k = 8$ | $k = 8$ | $k = 7$ | $k = 6$ |

where we used the same notation as in Figure 19.

**The nodal tangle.** Since $X$ is non-monotone, $B_0$ has either $w > 0$ or at least one parked node. For the second case, let $\mathfrak{a}_n$ be a parked node with smallest height $\alpha_n = \mathscr{F}(\mathfrak{a}_n)$. Sliding $\mathfrak{a}_n$ once around the curve $\mathscr{F}^{-1}(\alpha_n)$ gives us the construction alluded to in [Sch24, Remark 3.2] leading to a counterexample to the Lagrangian Poincaré recurrence conjecture.

Here we give a different approach that provides slightly different examples in the case $w > 0$. Then $H$ contains the line segment $\Delta_M$. As can be seen from Figure 19, there are two types of ends of $\Delta_M$: either there are three nodes pointing to the end or two. Depending on the types of ends of $\Delta_M$, modify $H$ by a nodal tangle as in Figure 23. Let $\mathfrak{a}_1, \ldots, \mathfrak{a}_n$ be the parked nodes of $B_0$ with height $\mathscr{F}(\mathfrak{a}_i) = \alpha_i$.

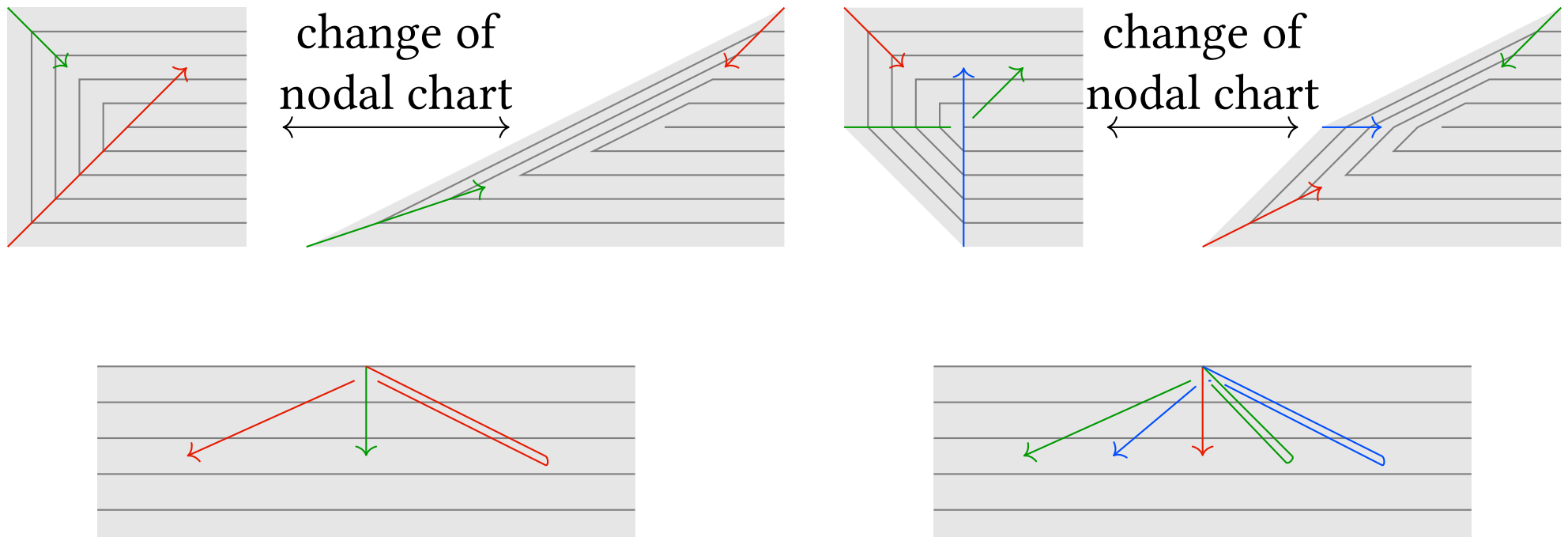

Figure 23: Nodal tangle transforming part of an $\varepsilon$-hat of positive width. On the left the nodal tangle for an end of $\Delta_M$ with two incoming nodes, on the right the one for an end with three incoming nodes. The top row shows the nodal tangle in a nodal chart as in Figure 19, the bottom row shows the same nodal tangle in a nodal chart as in Figure 22.

Further modify $B_0$ by a nodal tangle by sliding each parked node $\mathfrak{a}_i$ to the left by $2(M - \alpha_i)$. Denote by $(B_1, \mathfrak{N}_1)$ the resulting nodal integral affine surface, and by $(B \times [0,1], \mathfrak{N})$ the nodal tangle connecting $(B_0, \mathfrak{N}_0)$ to $(B_1, \mathfrak{N}_1)$. We also get the natural nodal chart $\varphi_1 : (B_1, \mathfrak{N}_1) \to \mathbb{R}^2$ as in Figure 22 b). There is a integral affine isomorphism $\tau : (B_0, \mathfrak{N}_0) \to (B_1, \mathfrak{N}_1)$: The nodal chart diagrams in Figure 22 a) and b) of $(B_0, \mathfrak{N}_0)$ and $(B_1, \mathfrak{N}_1)$ are related by applying the shear matrix $\left(\begin{smallmatrix} 1 & 2 \\ 0 & 1 \end{smallmatrix}\right)$ to Figure 22 a) and cutting the orange triangle off Figure 22 a) and regluing it as indicated in Figure 22 b). We quickly verify that this indeed gives an integral affine isomorphism: For points in the interior of the nodal chart diagrams in Figure 22 this is obvious. For points in a neighbourhood of $\Delta_M$ this can be verified in a nodal chart as in Figure 23.

**The Hamiltonian diffeomorphism.** Let $U$ be an $\varepsilon$-neighbourhood in $[0, M]$ of the points $\{M, \alpha_1, \ldots, \alpha_n\}$. Then we have $\mathcal{F}(\pi_B(\mathfrak{N})) \subset U$. Let $\pi : X \times [0,1] \to (B \times [0,1], \mathfrak{N})$ be a lift of the nodal tangle $(B \times [0,1], \mathfrak{N})$ supported on $\mathcal{F}^{-1}(U)$ (Theorem 2.30). The symplectic manifold $X$ is determined by $(H, M, \alpha_1, \ldots, \alpha_n)$, see Theorem 3.25. Since $\tau$ is an integral affine isomorphism identifying $B_0$ and $B_1$, Theorem 2.27 gives a symplectomorphism $\psi \in \mathrm{Symp}(X)$ such that

$$\begin{array}{ccc} X & \xrightarrow{\;\psi\;} & X \\ \downarrow{\scriptstyle \pi_0} & & \downarrow{\scriptstyle \pi_1} \\ (B_0, \mathfrak{N}_0) & \xrightarrow{\;\tau\;} & (B_1, \mathfrak{N}_1) \end{array}$$

commutes outside of the set $\mathcal{F}^{-1}(U)$. Let $x \in \mathcal{F}^{-1}([0,M] \setminus U)$ in $B_0$, and $h = \mathcal{F}(x)$ its height. Then $\psi(\pi_0^{-1}(x)) = \pi_1^{-1}(\tau(x)) = \pi_0(2(M - h) \cdot x)$, where $t \cdot$ denotes the $\mathbb{R}$-action on $B_0 \setminus \Delta_M$ translating a point $x$ by integral affine distance $t$ in the clockwise direction along the level set $\mathscr{F}^{-1}(h)$. The level set $\mathscr{F}^{-1}(h)$ has integral affine length $g(h)$, so for $x \in \mathscr{F}^{-1}([0, M] \setminus U)$ with $\frac{g(h)}{2(M-h)}$ irrational we have $\psi^n(\pi_0^{-1}(x)) \cap \pi_0^{-1}(x) = \emptyset$ for all $n \geq 0$. [LLW22, Corollary 4.6] says that the symplectic mapping class group of $X$ is finite, thus $\psi^m \in \mathrm{Symp}_0(X) = \mathrm{Ham}(X)$ for some $m \geq 0$. □

## 3.4 Probes and displacement energy

We want to use probes to give an upper bound on the displacement energy $e$ of a fibre of an almost toric fibration $\pi : X \to B$. Probes where introduces by McDuff in [McD11] in the context of toric symplectic manifolds. Here we consider the following slight generalisation:

**Definition 3.27.** Let $(B, \mathfrak{N})$ be a nodal integral affine surface. A **probe** $P : [0, d) \to B \setminus \mathfrak{N}$ is a rational straight line segment such that:

1. $P$ intersects $\partial B$ integrally transversely, i.e. $T_{P(0)}\partial B \cap \Lambda_{P(0)}$ and $P'(0)$ span the lattice $\Lambda_{P(0)}B$. In particular $P$ does not intersect $\partial B$ in a corner of $B$.
2. $P$ has no self-intersections.

Note that it follows from point 1 that $P'$ is primitive. With this definition we can choose action-coordinates on a contractible neighbourhood of $P$ where one can perform the symplectic reduction needed for the results in [McD11]. In particular we get the following version of [Bre23, Proposition 3.4]:

**Lemma 3.28.** *If $P : [0, d) \to B \setminus \mathfrak{N}$ is a probe and $\pi : X \to (B, \mathfrak{N})$ is an almost toric fibration, then for $t < \frac{d}{2}$,*

$$e(\pi^{-1}(P(t))) \leq t \, .$$

Let $\Delta$ be a weakly Delzant polygonal domain. Let $(B_0, \mathfrak{N}_0)$ be the nodal integral affine surface obtained by smoothing all vertices of $\Delta$, and $\varphi_0 : (B_0, \mathfrak{N}_0) \to \Delta$ the associated nodal chart, see Figure 15 a).

**Theorem 3.29.** *Set $M = \sup(\mathscr{F}_\Delta)$. Let $\varphi_0 : (B_0, \mathfrak{N}_0) \to \Delta$ as above, and $x \in \Delta \setminus \mathscr{K}_\Delta$ with $\mathscr{F}_\Delta(x) < \frac{1}{2}M$. Then there exists a nodal tangle $(B \times [0, 1], \mathfrak{N})$ connecting $(B_0, \mathfrak{N}_0)$ to $(B_1, \mathfrak{N}_1)$ with $x \notin \varphi_0(\pi_B(\mathfrak{N}))$ and a probe $P : [0, d) \to (B_1, \mathfrak{N}_1)$ with $P(\mathscr{F}_\Delta(x)) = \varphi_1^{-1}(x)$ and $2d > \mathscr{F}_\Delta(x)$.*

Using Lemma 3.28 we immediately get the symplectic version of Theorem 3.29:

**Corollary 3.30.** *Let $\varphi_0 : (B_0, \mathfrak{N}_0) \to \Delta$ be as above, $x \in \Delta \setminus \mathscr{K}_\Delta$ with $\mathscr{F}_\Delta(x) < \frac{1}{2}M$, and let $\pi_0 : X \to (B_0, \mathfrak{N}_0)$ be the unique almost toric fibration over $(B_0, \mathfrak{N}_0)$.*

*Then $e((\varphi_0 \circ \pi_0)^{-1}(x)) < \mathscr{F}_\Delta(x)$.*

*Proof of Theorem 3.29.* Take $M_\varepsilon < M$ so close to $M$ that $2\mathscr{F}_\Delta(x) < M_\varepsilon$. Let $(B \times [0,1], \mathfrak{N})$ be the nodal tangle constructed in Proposition 3.17. Since $x \notin \mathscr{K}_\Delta$, we may choose an open neighbourhood $V$ of $\mathscr{K}_\Delta$ such that $x \in \Delta \setminus V$. By the construction of $(B \times [0,1], \mathfrak{N})$ we may assume that $\pi_B \mathfrak{N} \subset V$. Let $\pi : X \times [0,1] \to (B \times [0,1], \mathfrak{N})$ be a lift of $(B \times [0,1], \mathfrak{N})$ supported on $V$. In particular we have that $(\varphi_0 \circ \pi_0)^{-1}(x) = (\varphi_1 \circ \pi_1)^{-1}(x)$, where $\varphi_1 = \varphi_0 \circ \tau_1^0$.

We proceed to construct a probe $P$ displacing $(\varphi_0 \circ \pi_0)^{-1}(x)$.

Since $x \notin \mathscr{K}_\Delta$, for small $a \in \mathbb{R}^2$ we can write $\mathscr{F}_\Delta(x + a) = \mathscr{F}_\Delta(x) + \langle \lambda, a \rangle$ for some primitive $\lambda \in \mathbb{Z}^2$. Let $u$ be a primitive integral vector in $\Lambda_{\varphi_1^{-1}(x)} B_1$ such that

$$\partial_{D_{\varphi_1^{-1}(x)}\varphi_1(u)} \mathscr{F}_\Delta = \langle \lambda, D_{\varphi_1^{-1}(x)}\varphi_1(u) \rangle = 1 \ . \tag{4}$$

Here $D\varphi_1$ is well-defined at $\varphi_1^{-1}(x)$ since $\varphi_1$ is integral affine near $\varphi_1(x)$. Set $t_x = \mathscr{F}_\Delta(x)$. Take $P$ to be the straight line with initial condition $P(t_x) = x$ and $P'(t_x) = u$. We take $P$ to have maximal domain such that $\operatorname{im} P \subset \varphi_1^{-1}(\Delta \setminus \Delta_{M_\varepsilon})$.

Since, by Proposition 3.17, $\mathscr{F}_1$ is integral affine on $\varphi_1^{-1}(\Delta \setminus (\Delta_{M_\varepsilon} \cup \mathfrak{N}_1)) \subset (B_1, \mathfrak{N}_1)$, so is $\mathscr{F}_1 \circ P$. In particular $(\mathscr{F}_1 \circ P)'$ is constant equal to $(\mathscr{F}_1 \circ P)'(t_x) = 1$ by (4).

The only nodes in $\Delta \setminus \Delta_{M_\varepsilon}$ are the ones we parked near vertices of $\mathscr{K}_\Delta$ in Proposition 3.17. By performing small nodal slides at the parked nodes, we may assume that $P$ does not hit a node. (The parked nodes slide in level sets of $\mathscr{F}_1$, in particular not parallel to the direction of $P$.)

Then the domain of $P$ is $[0, M_\varepsilon)$ and $\mathscr{F}_1 \circ P = \mathrm{id}_{[0,M_\varepsilon)}$. Furthermore $(\mathscr{F}_1 \circ P)'(0) = 1$ ensures that the intersection of $P$ and $\partial B = \mathscr{F}_1^{-1}(0)$ is integrally transverse. Since $\mathscr{F}_1 \circ P$ is injective, $P$ has no self intersections. So $P$ is a probe passing through $P(t_x) = \varphi_1^{-1}(x)$ with total length $M_\varepsilon > 2t_x$. $\square$

Corollary 3.30 gives an upper bound on the displacement energy $e$. In case $\Delta$ is Delzant, we get a matching lower bound:

**Lemma 3.31** ([Bre23, Proposition 3.4]). *Let $\mu : X \to \Delta$ be a toric symplectic manifold with $\Delta$ a Delzant polygonal domain. Then for all $x \in \Delta$ we have $\mathscr{F}_\Delta(x) \leq e(\mu^{-1}(x))$.*

Combining Corollary 3.30 and Lemma 3.31 we get:

**Corollary 3.32.** *Let $\mu : X \to \Delta$ be a toric symplectic manifold over a Delzant polygonal domain $\Delta$. Then for all $x \in \Delta \setminus \mathscr{K}_\Delta$ with $\mathscr{F}_\Delta(x) < \frac{1}{2} \sup \mathscr{F}_\Delta$ we have $\mathscr{F}_\Delta(x) = e(\mu^{-1}(x))$.*

*Remark* 3.33. It may be possible to get rid of the restriction $\mathscr{F}_\Delta(x) < \frac{1}{2} \sup \mathscr{F}_\Delta$ in Theorem 3.29 and Corollary 3.32. Indeed, by continuing the straight line $P$ in the proof of Theorem 3.29 through $\varphi_1^{-1}(\Delta_{M_\varepsilon})$, we can often get a *symmetric probe*, see [ABM14, §1.2.1], allowing us to displace all fibres of $\Delta \setminus \mathscr{K}_\Delta$ with the desired displacement energy. The challenge here is to avoid self-intersections of $P$. In all toric examples I examined this approach works, however I do not know a general method that always works.

# 4 Lagrangian knots

We summarize the procedure used to construct new Lagrangian torus knots using almost toric fibrations. This method was first used by Vianna in [Via16; Via17] to obtain monotone torus knots. The method goes as follows:

1. Create a nodal tangle $(B \times I, \mathfrak{N})$ and a lift $\pi : X \times I \to (B \times I, \mathfrak{N})$.
2. Maybe $\pi_0^{-1}(x) \not\cong \pi_t^{-1}(x)$ for some $x \in B$? Use your favourite invariant to try to distinguish them.

To distinguish the monotone tori in [Via16; Via17], Vianna uses holomorphic disc counts. We will use the displacement energy germ, mostly relying on Corollary 3.32 to compute it. The displacement energy germ was first used for this purpose in [CS10].

Lemma 2.32 tells us that this approach will only work for $x \in \pi_B(\mathfrak{N})$, as other fibres are mapped to symplectomorphic fibres under nodal tangles. For fibres over $x \in \pi_B(\mathfrak{N})$, Corollary 2.35 tells us that the invariant germs of $\pi_0^{-1}(x)$ and $\pi_1^{-1}(x)$ are related by the transition map $\tau_0^1$ created by the nodal tangle, so our goal is to create a tangle that provides a "complicated" $\tau_0^1$ in order to change the given invariant germ.

In general, we are given $x \in B_0$ and some nodes $\mathfrak{n}_1, \ldots, \mathfrak{n}_n$ whose eigenlines pass through $x$. Then we slide some of the nodes $\mathfrak{n}_1, \ldots, \mathfrak{n}_n$ back and forth over $x$. This situation can be described by using (rank 2) cluster algebras, see for example [STW16]. To get the setup of [STW16], take straight line segments connecting $x$ with the nodes $\mathfrak{n}_1, \ldots, \mathfrak{n}_n$. The *visible Lagrangians* (see e.g. [Eva23, Section 5.1]) over these straight line segments are Lagrangian discs $D_1, \ldots, D_n$ with boundary on $\pi^{-1}(x)$. The collection $(\pi^{-1}(x), D_1, \ldots, D_n)$ gives the Lagrangian skeleton of [STW16]. A *mutation* at $D_i$ is given by a *la-disc surgery* introduced in [Yau13] and corresponds to sliding the node $\mathfrak{n}_i$ over $x$. See also [ABR23] for a cluster theoretic setup that can be adopted almost directly to our case. Here we take a more ad hoc approach for the cases of one or two nodes sliding over $x$.

## 4.1 Rectifying map for one node

Let $\pi_0 : X \to (B_0, \mathfrak{N}_0)$ be an almost toric fibration and let $\mathfrak{n} \in \mathfrak{N}_0$ be a node with eigenray $\ell$ passing through a point $x_0 \in B_0 \setminus \mathfrak{N}_0$.

Let $\varphi_0 : U \to \mathbb{R}^2$ be a nodal chart with $\{\mathfrak{n}, x_0\} \subset U$ which is integral affine near $x_0$, and take $(B \times I, \mathfrak{N})$ to be a nodal tangle starting at $(B_0, \mathfrak{N}_0)$ such that $\mathfrak{n}$ slides over $x_0$ once. Using $\varphi_0$ we identify $U$ with $\varphi_0(U)$ for the rest of this section. Let $v \in \mathbb{Z}^2$ be the monodromy vector of $\mathfrak{n}$ and $k$ the multiplicity of $\mathfrak{n}$. Then, using Remark 2.8, the rectifying map at $x_1 = \tau_0^1 x_0$ of $\varphi_1 = \varphi_0 \circ \tau_1^0$ is given by the half-shear $h_v^k$. Let $\pi : X \times I \to (B \times I, \mathfrak{N})$ be a lift of the nodal tangle. Now we may use Corollary 2.35 and Lemma 2.37 to distinguish $\pi_0^{-1}(x)$ and $\pi_1^{-1}(x)$.

Suppose we have an invariant $\mathcal{I} : \mathcal{L} \to \mathbb{R}$ which is affine near $x$. For example this is the case if $\mathcal{I}$ is the displacement energy and we are in the setting of

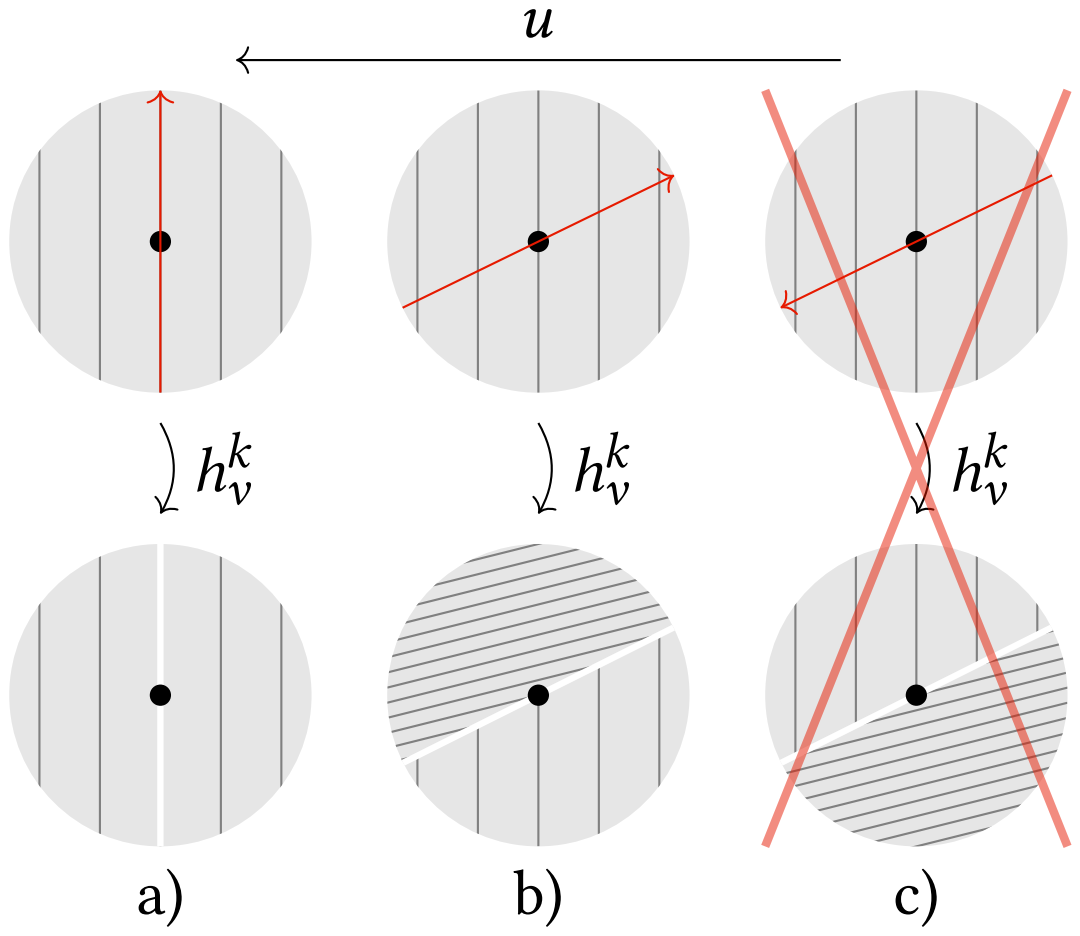

Figure 24: Modifying invariant germs by sliding over one node.

Corollary 3.32 with $x$ not being in the caustic.

Using Lemma 2.37 we may identify the invariant germ $[\mathcal{I}]_{\pi_0^{-1}(x_0)}$ with $[\mathcal{I}]_{x_0}$, and similarly for $x_1$. Using Corollary 2.35 we can write $[\mathcal{I}]_{x_1}|_{U\setminus\ell} = [\mathcal{I}]_{x_0}|_{U\setminus\ell} \circ h_v^k$, where we have used that the rectifying map of $\varphi_1$ at $x_1$ is given by $h_v^k$.

Since $\mathcal{I}$ is linear near $x_0$, there is $u \in \mathbb{R}^2$ such that $\mathcal{I}(a) = \langle u, a\rangle$ for $a \in U$ near $x_0$. We distinguish the following cases, illustrated in Figure 24:

a) If $\langle u, v\rangle = 0$ then $[\mathcal{I}]_{x_1}|_{U\setminus\ell} = [\mathcal{I}]_{x_0}|_{U\setminus\ell} \circ h_v^k = [\mathcal{I}]_{x_0}|_{U\setminus\ell}$, so $\pi_0^{-1}(x)$ cannot be distinguished from $\pi_1^{-1}(x)$ using Corollary 2.35.

b) If $\langle u, v\rangle < 0$, then the germs $[\mathcal{I}]_{x_1}|_{U\setminus\ell} = [\mathcal{I}]_{x_0}|_{U\setminus\ell} \circ h_v^k$ and $[\mathcal{I}]_{x_0}$ are not related by an integral linear transformation, so $\pi_0^{-1}(x)$ and $\pi_1^{-1}(x)$ are not related by symplectomorphism. Let $p = -\langle u, v\rangle, q = \det(u, v)$. Then the invariant germ of $\pi_1^{-1}(x)$ is the same (away from $\ell$) as the displacement energy germ of the torus $T_{pq}^k(\mathcal{I}(x))$ constructed in [BHS24].

c) We think that $\langle u, v\rangle > 0$ never happens (at least when $\mathcal{I}$ is the displacement energy), as it would lead to an invariant germ which does not look like a tropical polynomial, see [Bre25, Question 5.4 and Remark 2.18].

## 4.2 Entangling nodes: Rectifying maps for two nodes

Let $\pi_0 : X \to (B_0, \mathfrak{N}_0)$ be an almost toric fibration and $\mathfrak{a}, \mathfrak{b} \in \mathfrak{N}_0$ two nodes. Assume that the eigenlines $\gamma_0, \gamma_1$ of $\mathfrak{a}$ and $\mathfrak{b}$ respectively have some intersection point $x \in B$. Modifying $B_0$ by a nodal tangle if necessary, we may assume that $\mathfrak{a}$ and $\mathfrak{b}$ are contained in a small contractible neighbourhood $U_0$ of $x$ containing no other nodes.

We want to study the nodal tangle obtained by alternatingly sliding $\mathfrak{a}$ and $\mathfrak{b}$ through $x$. We define it recursively: For $n \in 2\mathbb{Z}_{\geq 0}$, let $(B_{n+1}, \mathfrak{N}_{n+1})$ be the nodal integral affine surface obtained from $(B_n, \mathfrak{N}_n)$ by sliding $\mathfrak{a}$ through $x$ once, and let $(B_{n+2}, \mathfrak{N}_{n+2})$ be the nodal integral affine surface obtained from $(B_{n+1}, \mathfrak{N}_{n+1})$ by sliding $\mathfrak{b}$ through $x$ once. We obtain a nodal tangle $(B \times \mathbb{R}_{\geq 0}, \mathfrak{N})$. See Figure 25.

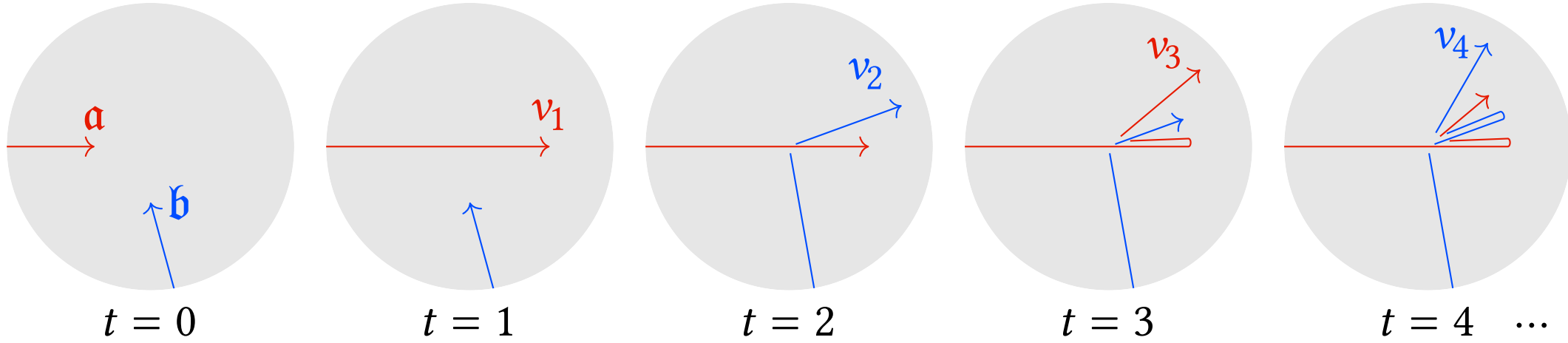

Figure 25: Entangling two nodes

Let $\varphi_0 : U \to \mathbb{R}^2$ be an integral affine chart of $(B_0, \mathfrak{N}_0)$ at $x$. Then $\varphi_t = \varphi_0 \circ \tau_0^t$ is a nodal chart of $(B_t, \mathfrak{N}_t)$ for all $t \geq 0$. For the rest of the discussion, identify $U$ with its image $\varphi_t(U) \subset \mathbb{R}^2$.

For $n \in \mathbb{Z}_{\geq 0}$ odd let $v_n$ be the primitive vector such that at time $t = n$ the node $\mathfrak{a}$ lies on the ray in direction $v_n$ emanating from $x$. For $n \in \mathbb{Z}_{\geq 0}$ even let $v_n$ be the primitive vector such that at time $t = n$ the node $\mathfrak{b}$ lies on the ray in direction $v_n$ emanating from $x$. In both cases $v_n$ coincides with the monodromy vector field of $\mathfrak{a}$ resp. $\mathfrak{b}$ near $\mathfrak{a}$ resp. $\mathfrak{b}$ at time $t = n$.

For times $t$ where $x \notin \mathfrak{N}_t$, denote by $R_t$ a choice of rectifying map of $\varphi_t$ at $x$, so that $R_t \circ \varphi_t$ is an integral affine chart of $(B_t, \mathfrak{N}_t)$ at $x$. From Figure 25 and Remark 2.8 we see that $R_1$ is given by the half-shear $h_{\mathfrak{a}} := h_{v_1}^{k_{\mathfrak{a}}}$, and $R_2$ is given by the composition of half-shears $h_{\mathfrak{b}} \circ h_{\mathfrak{a}}$, where $h_{\mathfrak{b}} := h_{-v_0}^{k_{\mathfrak{b}}}$. Set

$$R = -\operatorname{id} \circ R_2 = -h_{\mathfrak{b}} \circ h_{\mathfrak{a}},$$

which is another choice of rectifying map of $\varphi_2$ at $x$. Note that the nodes $\mathfrak{a}, \mathfrak{b}$ at time $t = 0$ are on the rays given by the pair $(-v_1, v_0)$, and for time $t = 2$ they are on the rays $v_1$ and $v_2 = -v_0 + k_{\mathfrak{a}} \det(v_0, v_1) v_1$, where we used Remark 2.20 to calculate $v_2$. We quickly check that $Rv_1 = -v_1$ and $Rv_2 = v_0$, so the images of $\varphi_0$ and $R\varphi_2$ are identical. Iterating, we can write the other rectifying maps as $R_{2n} = (-R)^n = (h_{\mathfrak{b}} \circ h_{\mathfrak{a}})^n$ and $R_{2n+1} = h_{\mathfrak{a}} \circ R_{2n} = h_{\mathfrak{a}} \circ (h_{\mathfrak{b}} \circ h_{\mathfrak{a}})^n$. We also get an explicit formula for $v_n$:

$$\begin{aligned} v_{2n-1} &= R^{-n}(-v_1) \\ v_{2n} &= R^{-n}(v_0) \end{aligned} \tag{5}$$

Set $a = \det(v_0, v_1) k_{\mathfrak{a}}, b = \det(v_0, v_1) k_{\mathfrak{b}}$ and take $(v_1, -v_0)$ as basis of $\mathbb{R}^2$. In this basis the piecewise linear map $R$ can be written as in Figure 26, where

$$A := s_{v_1}^{k_{\mathfrak{a}}} = \begin{pmatrix} 1 & -a \\ 0 & 1 \end{pmatrix}, \quad B := s_{v_0}^{k_{\mathfrak{b}}} = \begin{pmatrix} 1 & 0 \\ b & 1 \end{pmatrix}.$$

The following proposition says when we can expect "interesting" rectifying maps. Recall that $a, b \in \mathbb{Z}_{\geq 1}$.

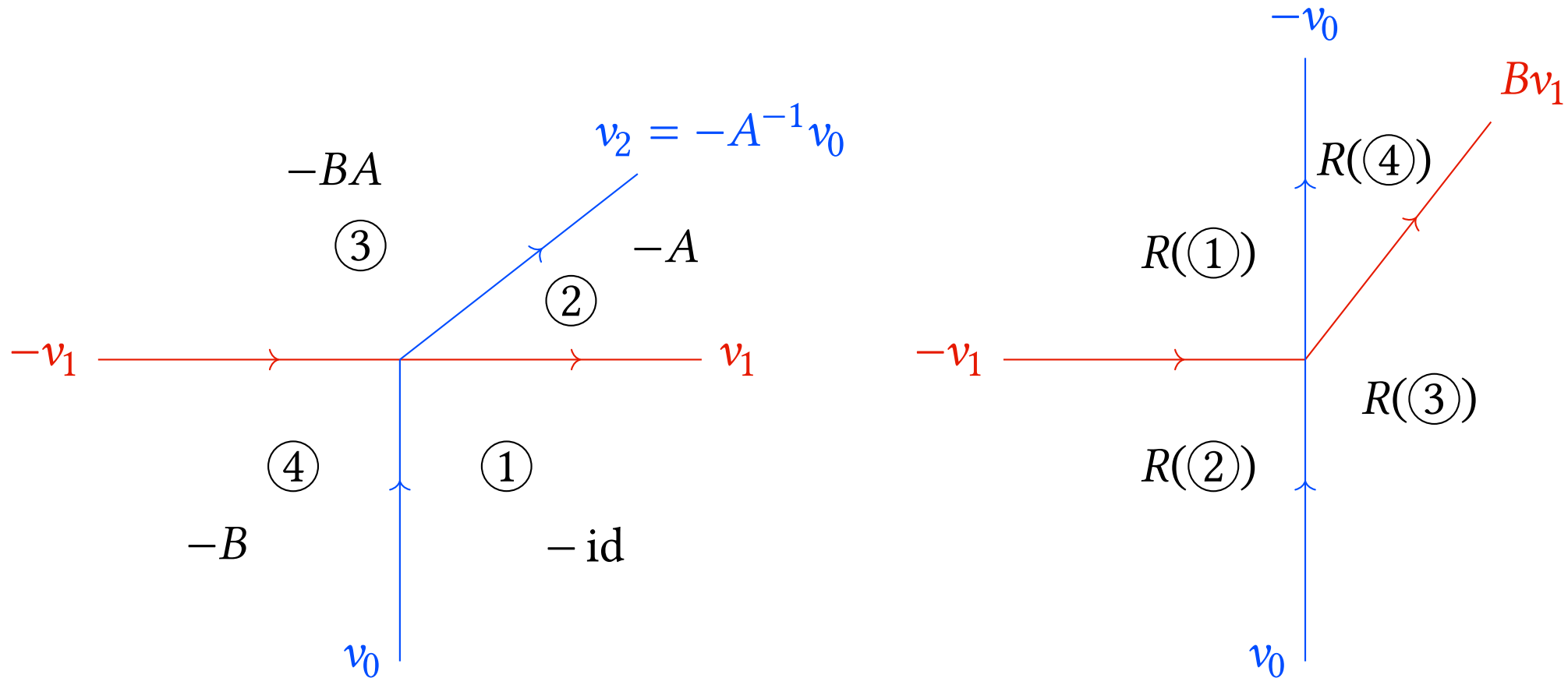

Figure 26: The rectifying map $R$. On the left the domain of $R$ and on the right its image. The domain is divided into four closed cones ①, ②, ③, ④ on which $R$ acts linearly by the matrices $-\mathrm{id}, -A, -BA, -B$ respectively. The rays dividing these regions form the cut graph of $R$ (and of $R^{-1}$ on the right), and are labelled by the corresponding orientation.

**Proposition 4.1.** *We have*

- *If* $ab = 1$ *then* $R_n$ *is linear iff* $n \in 5\mathbb{Z}_{\geq 0}$.
- *If* $ab = 2$ *then* $R_n$ *is linear iff* $n \in 6\mathbb{Z}_{\geq 0}$.
- *If* $ab = 3$ *then* $R_n$ *is linear iff* $n \in 8\mathbb{Z}_{\geq 0}$.
- *If* $ab \geq 4$ *then* $R_n$ *is not linear for any* $n > 0$. *The rectifying map* $R_{2n} = (-1)^n R^n$ *can be written as in Figure 28, and we may write*

$$\begin{aligned} v_{2n-1} &= bv_{2n-2} - v_{2n-3} = (-BA)^{1-n} v_1 \\ v_{2n} &= av_{2n-1} - v_{2n-2} = (-BA)^{1-n} v_2 \ . \end{aligned} \tag{6}$$

The proof will require some preparations.

*Remark* 4.2. The cases $ab = 1, 2, 3$ appear in [CV21, Section 3.2], where they are related to $A_2, B_2, G_2$ cluster relations.

This is a general phenomenon: Using the setup of [ABR23], we see that a cluster relation corresponds to a local nodal braid, that is a nodal tangle near $x$ such that $(U_0, \mathfrak{N}_0) = (U_1, \mathfrak{N}_1)$. The nodal braid for $ab = 1$ is illustrated in Figure 27.

*Remark* 4.3. In the last case of Proposition 4.1, if $a = b \geq 2$ then the recurrence relation in (6) is a linear recurrence relation with constant coefficients:

$$v_n = av_{n-1} - v_{n-2} \ . \tag{7}$$

By the using e.g. generating functions we can thus find a closed form for $v_n$: Let $\mu^{\pm 1} = \frac{a \pm \sqrt{a^2-4}}{2}$ be the solutions to the characteristic polynomial of (7). Then we

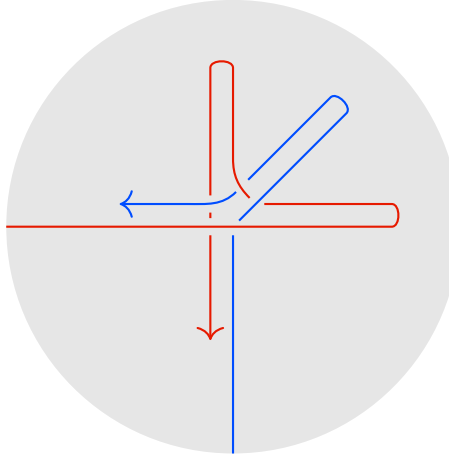

Figure 27: Nodal tangle diagram for the case $ab = 1$ at $t = 5$. Here the rectifying map at $x$ is integral linear.

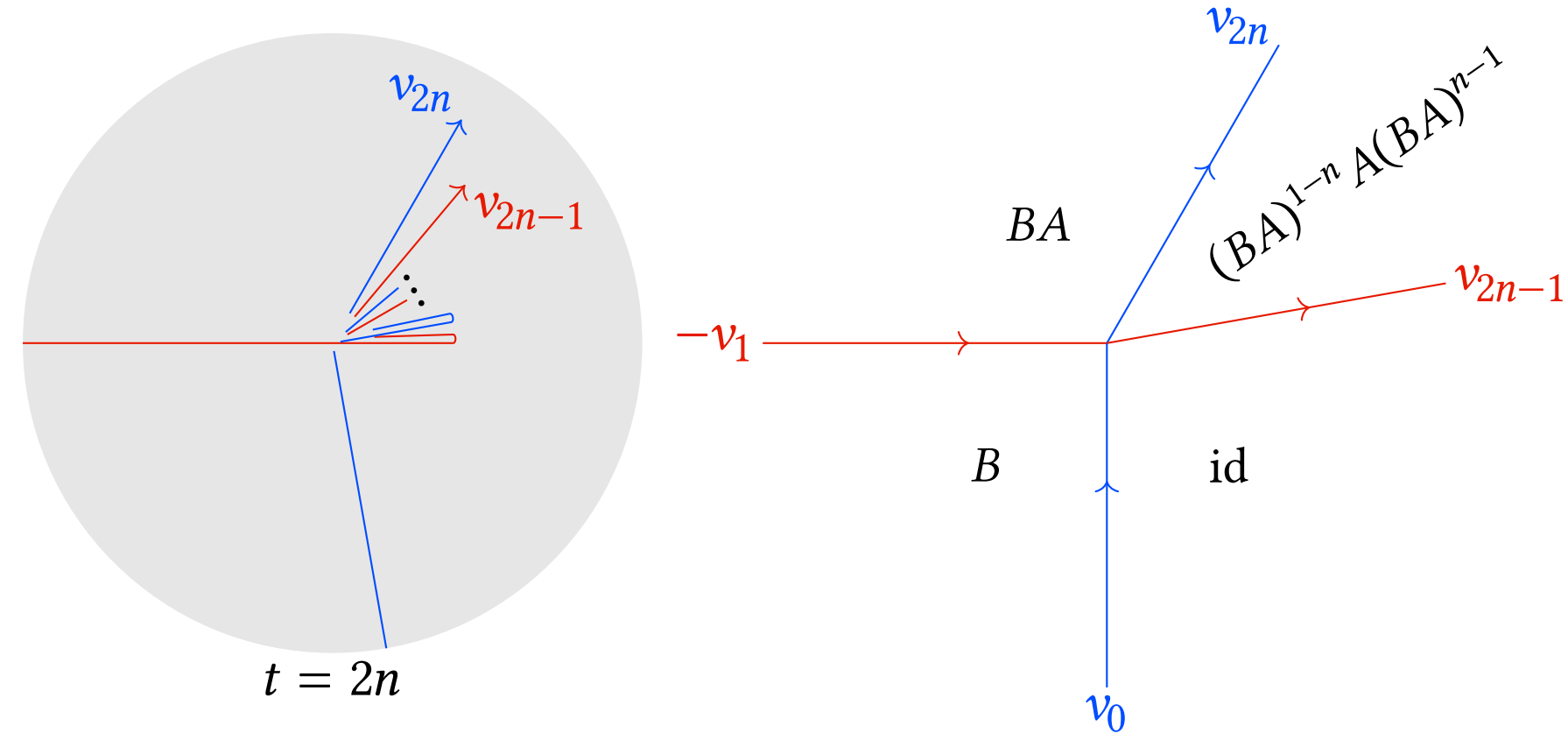

Figure 28: The rectifying map $R_{2n}$

have

$$v_n = \begin{cases} v_0 + n(v_1 - v_0) & \text{if } a = 2 \\ \frac{1}{\sqrt{a^2-4}}\big((\mu^n - \mu^{-n})v_1 - (\mu^{n-1} - \mu^{1-n})v_0\big) & \text{if } a > 2 \end{cases}.$$

To prove Proposition 4.1 we use the following definition and lemma.

**Definition 4.4.** A vector $u \neq 0$ is an **eigenvector** of a piecewise linear map $T$ if

$$Tu = \lambda u$$

for some **eigenvalue** $\lambda \in \mathbb{R}_{>0}$.

**Lemma 4.5.** *Let $C$ be the cone* ③ $\cap R($③$)$.

- *If $ab < 4$ then $R$ has no eigenvector,*
- *if $ab = 4$ then $R$ has exactly one eigenvector, which is contained in $C$,*
- *if $ab > 4$ then $R$ has exactly two eigenvectors, which are contained in $C$.*

*Proof.* By definition, we have $C \subset$ ③. We see from Figure 26 that

$R($①$) \subset$ ③ $\qquad\qquad R($②$) =$ ④

$R($③$) \subset$ ① $\cup$ ② $\cup$ ③ $\qquad\qquad R($④$) \subset$ ③ ,

so ③ is the only region which may intersect its image under $R$, and there cannot be any eigenvectors in the regions ①, ②, ④. Since $R|_{③} = -BA$, if $R$ has eigenvectors in ③ then they are also eigenvectors of the integral linear map $-BA$.

In the basis $(v_1, -v_0)$, we can write $-BA$ as the matrix

$$\begin{pmatrix} -1 & a \\ -b & ab-1 \end{pmatrix} .$$

The discriminant of the characteristic polynomial of $-BA$ is $ab(ab-4)$, showing that $-BA$ has the desired number of eigenvectors. We claim they are contained in $C$ if they exist.

Let $ab \geq 4$. Calculating the eigenvalues and eigenvectors, we get

$$\lambda^{\pm 1} = \frac{ab - 2 \pm \sqrt{ab(ab-4)}}{2} , \quad u_{\lambda^{\pm 1}} = \begin{pmatrix} ab \mp \sqrt{ab(ab-4)} \\ 2b \end{pmatrix}$$

and using $\det(-BA) = 1$ we get $\lambda > 1 > \lambda^{-1} > 0$.

The cone $C$ is spanned by $(v_2, Bv_1) = ((a,1),(1,b))$. We can write

$$(ab-1)u_{\lambda^{\pm 1}} = 2b\lambda^{\mp 1} v_2 + 2(\lambda^{\pm 1} + 1)Bv_1$$

where the coefficients $2b\lambda^{\mp 1}$ and $2(\lambda^{\pm 1} + 1)$ are both positive, so $u_{\lambda^{\pm 1}}$ is indeed contained in $C$. □

*Proof of Proposition 4.1.* The cases $ab = 1, 2, 3$ are straightforward to check by hand using the fact that $ab = 1, 2, 3$ implies that $\det(v_0, v_1) = 1$ and $\{k_{\mathfrak{a}}, k_{\mathfrak{b}}\} = \{a, b\} = \{1, ab\}$. The case $ab = 1$ is illustrated in Figure 27.

So let $ab \geq 4$. Then by Lemma 4.5, $R^n$ has at least one eigenvector, and all of its eigenvectors are contained in the cone $C$. This means $R^n$ is not linear: If $R^n$ was linear, it would also have at least one eigenvector in the cone $-C$. Thus $R_{2n} = (-R)^n$ is also not linear. The cuts of $R_{2n+1}$ are given by the rays spanned by the vectors $-v_1, v_0, v_{2n}, v_{2n+1}$. Below we will see that these are all distinct, and thus $R_{2n+1}$ is also not linear.

By (5),

$$\begin{aligned} v_{2n-1} &= R^{-n}(-v_1) = R^{1-n}(v_1) \\ v_{2n} &= R^{-n}(v_0) \quad = R^{1-n}(v_2) . \end{aligned}$$

The pair $v_1, v_2$ is contained in the positive cone $\langle v_0, u_{\lambda^{-1}} \rangle_+$ spanned by $v_0, u_{\lambda^{-1}}$, with $u_{\lambda^{-1}}$ as in the proof of Lemma 4.5. (This follows from the fact that $u_{\lambda^{-1}}$ is contained in $C$ which is spanned by $(v_2, Bv_1)$.) On $\langle v_0, u_{\lambda^{-1}} \rangle_+ \subset R(③)$ the map $R^{-1}$ is given by $(-BA)^{-1}$. Under repeated iterations of $R^{-1}$, $u_{\lambda^{-1}}$ is the attractive eigenvector of $(-BA)^{-1}$, meaning that the sequence $\frac{v_n}{|v_n|}$ in $S^1$ converges monotonically to $\frac{u_{\lambda^{-1}}}{|u_{\lambda^{-1}}|}$, so that $v_n$ is contained in $\langle v_0, u_{\lambda^{-1}}\rangle_+$ for all $n \geq 1$. This lets us give an explicit formula for the $v_n$:

$$\begin{aligned} v_{2n-1} &= (-BA)^{1-n} v_1 \\ v_{2n} &= (-BA)^{1-n} v_2 \ , \end{aligned}$$

from which we easily deduce the recurrence

$$\begin{aligned} v_{2n-1} &= -v_{2n-3} + b v_{2n-2} \\ v_{2n} &= -v_{2n-2} + a v_{2n-1} \ . \end{aligned}$$

The cones where $R_{2n}$ is linear are given by the vectors $-v_1, v_0, v_{2n-1}, v_{2n}$ as indicated in Figure 28. Using Figure 26, we see that in the lower right cone $\langle v_0, v_{2n-1}\rangle_+$ the map $R_{2n} = (-R)^n$ is given by the identity.

Note that shear maps transform under linear bijections like

$$s_{Tv} = T s_v T^{-1} \ .$$

Now we can use this identity and Remark 2.8 to determine the map $R_{2n}$ on the remaining cones: The ray given by $v_{2n-1} = (-BA)^{1-n} v_1$ is a cut of weight $k_{\mathfrak{a}}$, and the associated shear map is given by

$$\begin{aligned} s^{k_{\mathfrak{a}}}_{(-BA)^{1-n} v_1} &= (-BA)^{1-n} s^{k_{\mathfrak{a}}}_{v_1} (-BA)^{n-1} \\ &= (-BA)^{1-n} A (-BA)^{n-1} = (BA)^{1-n} A (BA)^{n-1} \ . \end{aligned}$$

The linear maps on the two other cones are determined similarly. Similarly we can also determine the rectifying map for $R_{2n+1}$. ☐

The following theorem describes a situation where we get an infinite family of Lagrangian tori which are not Hamiltonian isotopic, see Figure 4.

**Theorem 4.6.** *Let $\pi : X \to (B_0, \mathfrak{N}_0)$ be an almost toric fibration, and $x_0 \in B_0$. Suppose that*

1. *There are two nodes $\mathfrak{a}, \mathfrak{b}$ incident at $x_0$ with*

$$k_{\mathfrak{a}} k_{\mathfrak{b}} \det(v_{\mathfrak{a}}, v_{\mathfrak{b}})^2 \geq 4$$

   *where $k_{\mathfrak{a}}, k_{\mathfrak{b}}$ are the multiplicities of $\mathfrak{a}$ and $\mathfrak{b}$ and $v_{\mathfrak{a}}, v_{\mathfrak{b}}$ are the monodromy vector fields of $\mathfrak{a}$ and $\mathfrak{b}$.*
2. *$\mathscr{I} : \mathscr{L} \to \mathbb{R}$ is an invariant of Lagrangians under symplectomorphisms such that the invariant germ $[\mathscr{I}]_{x_0} : B_0 \to \mathbb{R}$ is affine and non-constant on an open dense set near $x_0$.*

*Let $(B \times \mathbb{R}_{\geq 0})$ be the nodal tangle obtained as above by sliding $\mathfrak{a}$ and $\mathfrak{b}$ alternatingly over $x$.*

*Then the invariant germs $[\mathscr{I}]_{x_0} \circ \tau_n^0$ are pairwise different, hence the Lagrangian tori $\{\pi_n^{-1}(x)\}$ are pairwise not related by symplectomorphism.*

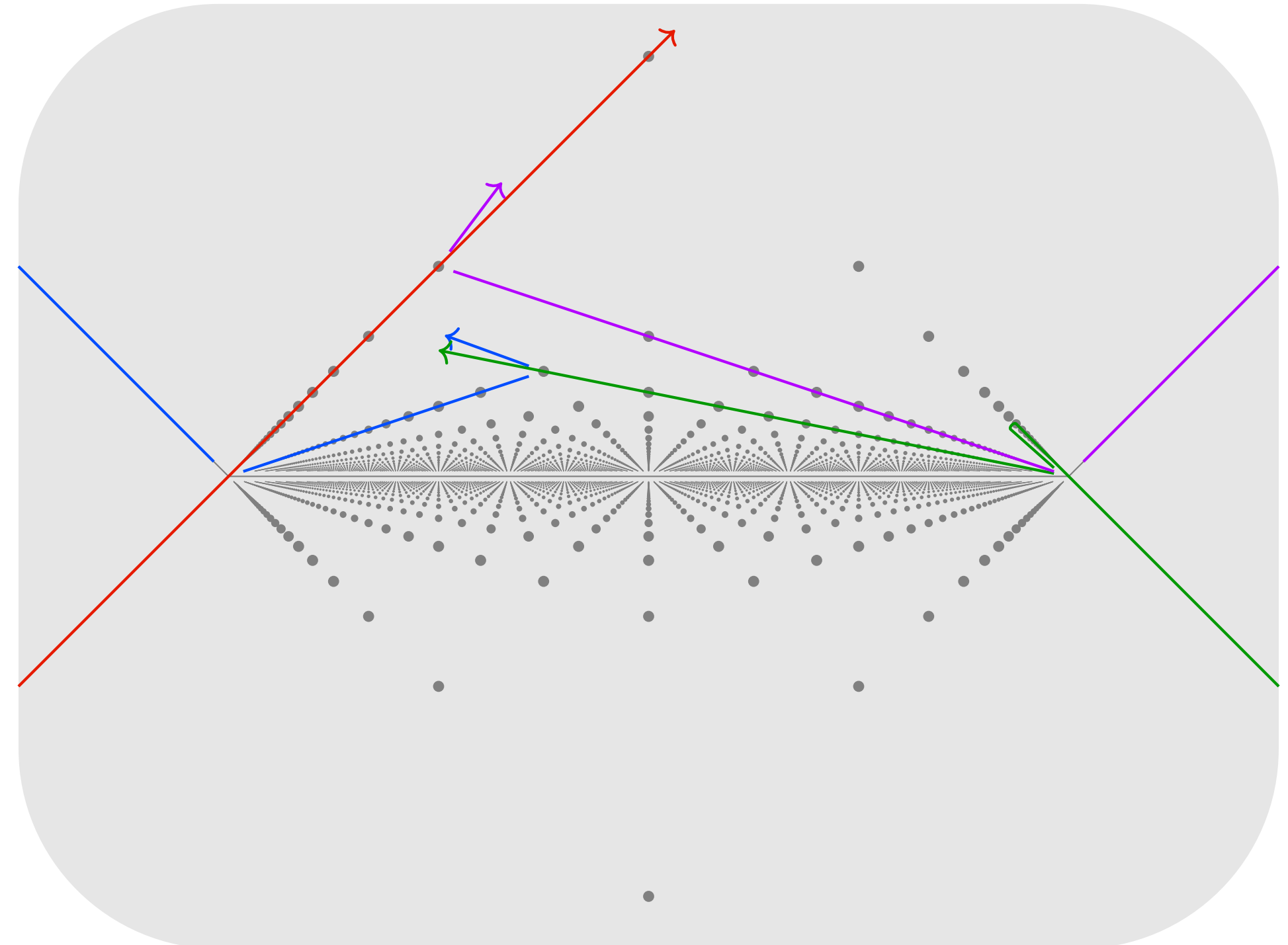

Figure 29: Some Lagrangian knots in $S^2 \times S^2$; we only picture the interior of the moment polytope. Grey dots correspond to points in $\mathscr{D}'$.

*Proof.* Condition 1. puts us in the case $ab \geq 4$ in Proposition 4.1.

Let $\varphi_0$ be an integral affine chart at $x_0$, and let $\varphi_n := \varphi_0 \circ \tau_n^0$ and $x_t = \tau_0^t x_0$. An integral affine chart at $x_n$ is given by $R_n \circ \varphi_n$, where $R_n$ is the rectifying map of $\varphi_n$ at $x_n$. Identifying a neighbourhood of $x_n$ with its image under $R_n \circ \varphi_n$, the invariant germ $[\mathscr{I}]_{x_0} \circ \tau_n^0$ is given by $[\mathscr{I}]_{x_0} \circ R_n$.

Since $[\mathscr{I}]_{x_0}$ is affine, the cut graph of $[\mathscr{I}]_{x_0} \circ R_n$ is the same as the one of $R_n$, except for possibly one edge missing. The cut graph of $R_n$ is given by the rays spanned by the vectors $-v_1, v_0, v_{n-1}, v_n$, see Figure 28 for the case where $n$ is even. The cut graphs of $R_n$ are all not related by an integral affine transformation by Equation (6), and removing one edge does not change that. By Lemma 2.37 this suffices to show that the Lagrangian tori $\pi_n^{-1}(x_n)$ are not related by symplectomorphism. $\square$

## 4.3 Examples

**Example 4.7** ($S^2 \times S^2$). Take $S^2 \times S^2$ with the symplectic form $\omega_\alpha = 2(1+\alpha)\omega_{S^2} \oplus 2\alpha\omega_{S^2}$ where $\alpha > 0$, and let $\square_\alpha = [-(1+\alpha), 1+\alpha] \times [-\alpha, \alpha] \subset \mathbb{R}^2$ be the moment image under the usual Hamiltonian torus action $\mu\colon (S^2 \times S^2, \omega_\alpha) \to \square_\alpha$. This normalization corresponds to the one chosen in [BK25]. It means that the caustic $\mathscr{K}_{\square_\alpha}$ is independent of $\alpha$.

After nodal trades at the corners of $\square_\alpha$, the points $(-1,0), (1,0)$ both lie on the intersection of eigenlines of nodes as in Section 4.2, see Figure 29. With the notation of Proposition 4.1, at $x = (-1, 0)$ we have two incoming nodes $\mathfrak{a}, \mathfrak{b}$ with $a = b = 2$ and $v_0 = (-1, 1), v_1 = (1, 1)$. Entangling these nodes at $x$, we see from Remark 4.3 that we can get nodes exiting $x$ in the directions $(2n + 1, 1)$ for $n \geq 0$. Swapping the roles of $\mathfrak{a}, \mathfrak{b}$ (i.e. sliding $\mathfrak{b}$ over $x$ first, resulting in setting $v_0 = (-1, -1), v_1 = (1, -1)$) we also get nodes exiting $x$ in the directions $(2n + 1, -1)$. Similarly, we can get nodes exiting $(1, 0)$ in the directions $(-(2n + 1), \pm 1)$ for $n \geq 0$.

Eigenrays of these nodes may again intersect in one of the points of

$$\mathscr{D}' = \left\{ \left( \frac{k - 2l - 1}{k}, \pm \frac{1}{k} \right) \,\middle|\, k \in \mathbb{Z}_{\geq 1}, l \in \{0, ..., k - 1\} \right\} .$$

Over these points, entangling two nodes gives infinitely many Lagrangian torus knots. The displacement energy germs of two families over two points in $\mathscr{D}'$ have the same displacement energy germs if and only if their $y$-coordinate is the same. (This means that in $\mathscr{D}'$ they share the same value of $k$, and we have that $a = 2k$ in Remark 4.3, resulting in the same displacement energy germs.)

Of course there are many other points in $\square_\alpha$ where we can entangle two nodes, which are however more difficult to determine.

The set $\mathscr{D}'$ coincides with the set of toric fibres which are not *ball non-monotone*, see [BK25, Theorem 1.23].

Note that at $(\pm 1, 0)$ entangling the nodes has no effect on the displacement energy germ. Indeed, using probes it can be shown that all tori constructed by entangling at $(\pm 1, 0)$ are Hamiltonian isotopic.

**Example 4.8** (Vianna Tori). A **Markov triple** $(p_1, p_2, p_3) \in \mathbb{N}^3$ is a triple satisfying the Markov equation

$$p_1^2 + p_2^2 + p_3^2 = 3 p_1 p_2 p_3 .$$

In [Via16], Vianna constructed a monotone Lagrangian knot in $\mathbb{C}P^2$ for every Markov triple by following the outline given at the start of Section 4. Denote by $T_{p_1, p_2, p_3}$ Vianna's torus corresponding to the Markov triple $p_1, p_2, p_3$. See [Eva23, Appendix I] for a nice exposition on the connection between the base diagrams and the Markov equation. We do not give a full description of the construction, but just describe the rectifying map obtained by following Vianna's sequence of nodal slides starting from the standard toric base diagram for $\mathbb{C}P^2$ and interpreting it as a nodal tangle.

In Figure 30 a), the fibre over the central point $x$ is the torus $T_{5,2,1}$. Realize a sequence of mutations to arrive at the Markov triple $(p_1, p_2, p_3)$ as a nodal tangle as done by Vianna and denote by $\varphi_1$ the corresponding nodal chart. We arrive at a situation where the rectifying map at $x$ of $\varphi_1$ is given as in Figure 30 b) where $(p_i, q_i)$ form an "extended Markov triple", including the $q_i$ described in [ES18]. Concretely, $p = (p_1, p_2, -p_3)$ and $q = (q_1, q_2, -q_3)$ satisfy the equation

$$p \times q = 3(p_1, p_2, \hat{p}_3) , \tag{8}$$

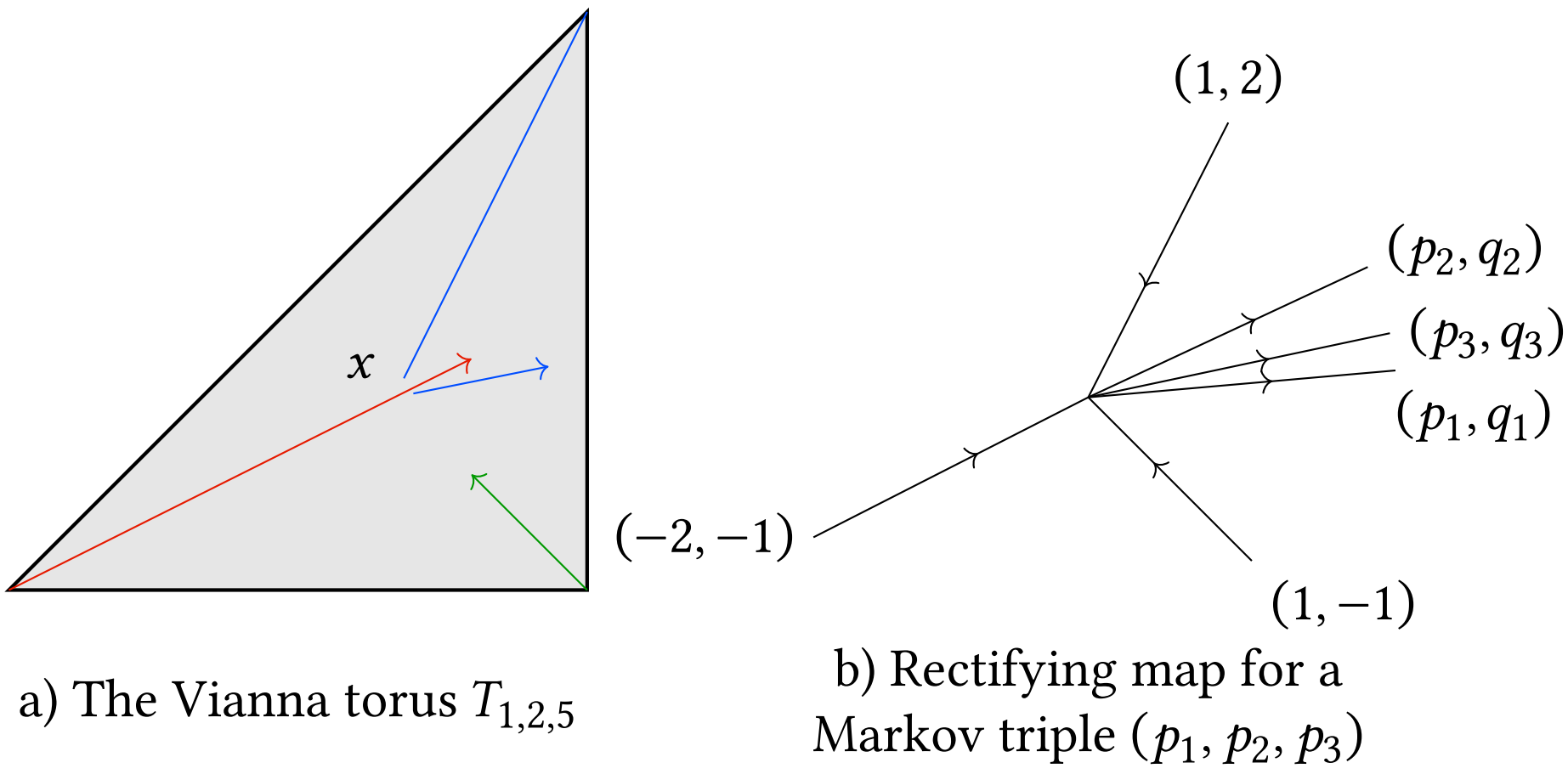

a) The Vianna torus $T_{1,2,5}$

b) Rectifying map for a Markov triple $(p_1, p_2, p_3)$

Figure 30: Constructing Vianna tori via nodal tangles

where $\hat{p}_3 = 3p_1p_2 - p_3$ is the Markov number given by mutation at $p_3$. This equation follows from the fact that $\frac{1}{p_1}(p_1, q_1), \frac{1}{p_2}(p_2, q_2), \frac{1}{\hat{p}_3}(-p_3, -q_3)$ form the corners of a *Vianna triangle*, see [Eva23, Appendix I] for details, especially [Eva23, Corollary I.13]. The Markov equation for $(p_1, p_2, p_3)$ follows also from (8).

In $\mathbb{C}P^2$, the displacement energy of toric fibres is given by $\mathscr{F}_\Delta$, except for the central fibre at $x$, which is non-displaceable. Applying the rectifying map $R_x$ of $\varphi_1$, as pictured in Figure 30 b), the level sets of $\mathscr{F}_\Delta \circ R_x^{-1}$ give smaller copies of the Vianna triangle with corners $\frac{1}{p_1}(p_1, q_1), \frac{1}{p_2}(p_2, q_2), \frac{1}{\hat{p}_3}(-p_3, -q_3)$, which allows us to distinguish the Vianna tori using the displacement energy germ.

**Example 4.9** (Farey tree & Lagrangian Pinwheels)**.** Let $\mathfrak{a} \in \mathfrak{N}$ be a node of a nodal integral affine surface, and let $\gamma$ be an eigenray of $\mathfrak{a}$ intersecting $\partial B$ in $x$. Near $x$, let $u$ be a primitive vector parallel to $\partial B$ and extend the monodromy pair $(v, \lambda)$ of $\mathfrak{a}$ along $\gamma$ to $x$. The *visible Lagrangian* over $\gamma$ (see [Eva23, Chapter 5]) is the *Lagrangian pinwheel* $L_{p,q}$ (see [ES18, Definition 2.1]), where $p = |\det(u, v)|$ and $q = \langle u, v \rangle \mod p$. Here the scalar product is defined by some choice of basis of $\Lambda_x B$, thus $q$ is only well-defined up to sign (since we don't fix an orientation of $B$).

In [ES18] it is shown that $\mathbb{C}P^2$ does only allow embeddings for certain pinwheels, namely if $p$ is a Markov number.

Alternatively, Lagrangian pinwheels may be obtained as the vanishing cycles of a $\mathbb{Q}$-Gorenstein smoothing of a Wahl singularity. For a del Pezzo surface $X$, in the recent [UZ25] it is shown which Wahl singularities arise from degenerations of $X$, meaning which Lagrangian pinwheels can be realized "algebraically". In particular [UZ25, Theorem 1.10] states that for $X = \mathbb{C}P^2 \# n\overline{\mathbb{C}P^2}$ with $n \geq 5$ all Lagrangian pinwheels $L_{p,q}$ admit embeddings into $X$.

The Farey tree is a well-known procedure to generate all primitive vectors in the positive quadrant of $\mathbb{Z}^2$:

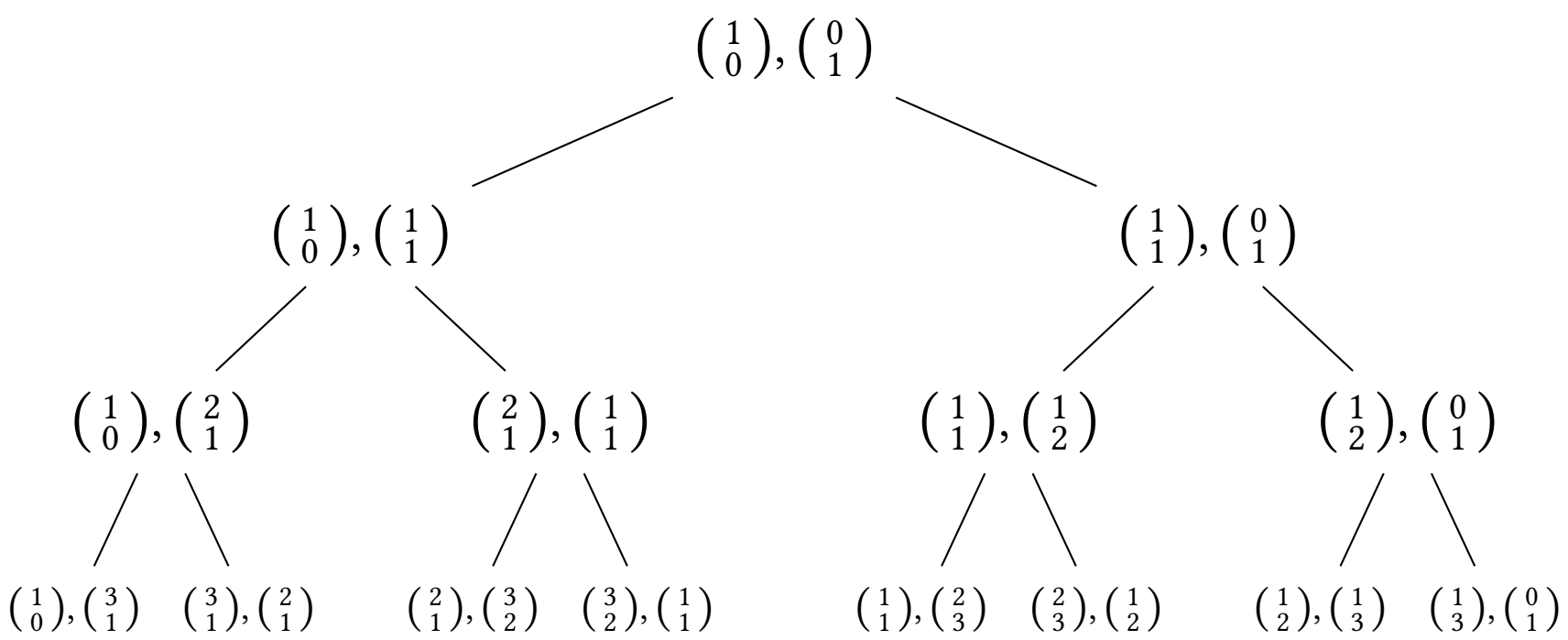

Figure 31: The first four rows of the Farey tree.

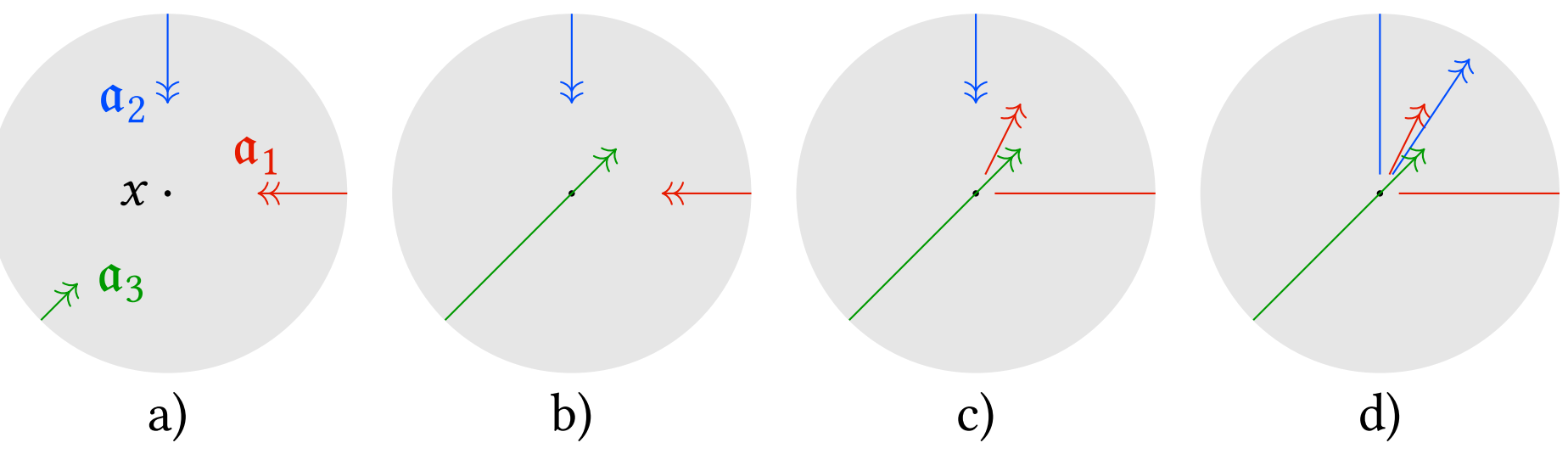

Figure 32: Stepping through the Farey tree with nodal slides. The double arrow heads mark that the nodes are of multiplicity two.

**Definition 4.10.** A **Farey pair** $(v_0, v_1)$ is a positive basis of $\mathbb{Z}^2$, i.e. $\det(v_0, v_1) = 1$. Its children are the Farey pairs $(v_0, v_0 + v_1)$ and $(v_0 + v_1, v_1)$.

The **Farey tree** is the binary tree of Farey pairs with root $\left(\binom{1}{0}, \binom{0}{1}\right)$. See Figure 31.

Reading the $n$-th row of the tree from left to right ignoring repeated entries, the vectors encode the numerator and denominator of the fractions in the *Farey sequence* $F_n$.

Suppose that $B$ is a nodal integral affine surface, $x \in B$ with three incoming nodes $\mathfrak{a}_1, \mathfrak{a}_2, \mathfrak{a}_3$ of multiplicity 2 on the rays originating from $x$ spanned by primitive vectors $v_1, v_2, v_3$ such that $\det(v_1, v_2) = \det(v_2, v_3) = \det(v_3, v_1) = 1$. This means that $x$ admits a nodal chart diagram as in Figure 32 a). Then $v_1 + v_2 + v_3 = 0$.

Pushing node $\mathfrak{a}_3$ through $x$, we get nodes on the rays $(v_1, v_1 + v_2, v_2)$, see Figure 32 b). Note that the outer vectors $(v_1, v_2)$ form a Farey pair. Pushing one of the "outer nodes" on $v_1$ or $v_2$ through $x$ by using the recipe in Remark 2.20, we get nodes either on the rays $(v_1 + v_2, (v_1 + v_2) + v_2, v_2)$ or $(v_1, v_1 + (v_1 + v_2), v_1 + v_2)$, with the outer vectors corresponding to the children of the Farey pair $(v_1, v_2)$. Iterating, we see that we can reconstruct the Farey tree through nodal tangles at $x$.

In particular, for any primitive vector $v \in \Lambda_x B$, there exists a nodal tangle at $x$ such that a node of multiplicity 2 slides along the ray spanned by $v$.

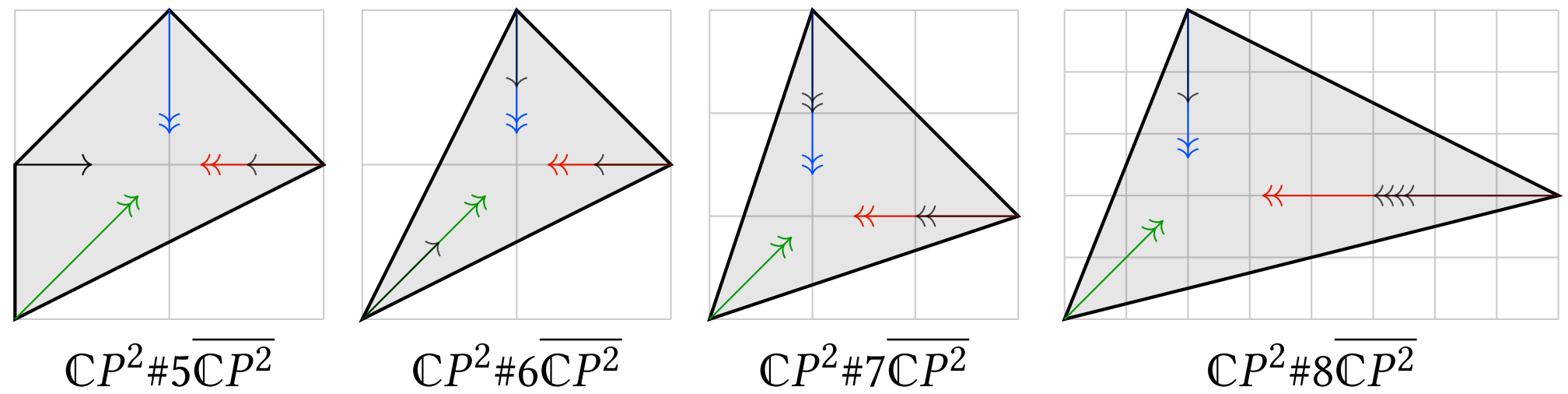

Figure 33: Nodal chart diagrams for monotone $\mathbb{C}P^2\#n\overline{\mathbb{C}P^2}$ with $n \geq 5$ showing the desired nodal configuration.

This configuration of nodes can be found in the monotone $\mathbb{C}P^2\#n\overline{\mathbb{C}P^2}$ with $n \geq 5$, as illustrated in Figure 33, see [Via17, Section 3] for a construction of these diagrams. In each case we can quickly check that for each primitive vector $(p, q)$, we can find a $L_{p,q}$ pinwheel: Take a nodal chart $\varphi$ as in Remark 3.18. The nodal tangles constructed above allow us to slide a node of multiplicity two along any primitive vector $v$ pointing from $x$ to the boundary. In the nodal chart $\varphi$ we can directly read off the type of pinwheel from the coordinates of $v = (v_1, v_2)$ in $\varphi$:

$$p = |\det(v, (1, 0))| = |v_2|, \quad q = \langle v, (1, 0) \rangle \mod p = v_1 \mod p.$$

Since any primitive vector $v$ can be realized by a nodal tangle as described above, we can find Lagrangian pinwheel $L_{p,q}$ for any coprime $p, q$.

# 5 Some open questions

## 5.1 The piecewise integral linear group

In this paper all piecewise integral linear transformations we encountered were generated by half-shears.

*Question* 5.1. Is the group of piecewise (integral) linear automorphisms of $\mathbb{R}^2$ generated by half-shears?

*Remark* 5.2. The answer for this question in dimensions higher than two is false. Take for example the piecewise integral linear automorphism

$$T : \mathbb{R}^3 \to \mathbb{R}^3$$
$$(x, y, z) \mapsto \begin{cases} (x, y, z) & \text{if } x \leq 0, y \leq 0 \\ (x, y, z + x) & \text{if } x \geq 0, x \geq y \\ (x, y, z + y) & \text{if } y \geq 0, y \geq x \end{cases}.$$

Indeed, if $T$ was a composition of half-shears, it would have an even number of linear domains.

## 5.2 Generalizing Symington's conjecture

Consider the set $\mathscr{B}_X$ of nodal integral affine surfaces that arise as almost toric bases of a symplectic four manifold $X$. Extending the definition of nodal tangle to include nodal trades, Theorem A says that if $X$ is toric, any two Delzant polygons $\Delta_0, \Delta_1 \in \mathscr{B}_X$ are connected by a nodal tangle.

*Question* 5.3. If $X$ is toric, are any two nodal integral affine surfaces in $\mathscr{B}_X$ related by a nodal tangle?

Note that there might be an almost toric base of $X$ which is not connected to a toric base by a nodal tangle.

We ask a stronger question:

*Question* 5.4. If $X$ is a symplectic rational surface, are any two nodal integral affine surfaces in $\mathscr{B}_X$ related by a nodal tangle?

For both these questions it might be possible to generalize the proof strategy of Theorem A, namely trying to deform elements of $\mathscr{B}_X$ into a canonical form which only depends on $X$.

Being a symplectic rational surface exactly means that all bases in $\mathscr{B}_X$ are topologically discs with marked points. This follows from [LS10, Table 1], which classifies the closed symplectic 4-manifolds $X$ that admit an almost toric fibration and their bases up to diffeomorphism. We can also ask the same question for the other types of $X$ appearing in the list given in [LS10]:

*Question* 5.5. If $X$ is a closed symplectic four manifold admitting an almost toric fibration, are any two nodal integral affine surfaces in $\mathscr{B}_X$ related by a nodal tangle?

## 5.3 Displacement energy of (almost) toric fibres

*Question* 5.6. If $\Delta$ is a Delzant polygonal domain, is the displacement energy of toric fibres given by $\mathscr{F}_\Delta$ for all points outside the caustic $\mathscr{K}_\Delta$?

See also Remark 3.33.

## 5.4 Tangling points

To apply our recipe for constructing Lagrangian knots, we need to find a point $x \in B$ where, perhaps after modifying $B$ by a nodal slide, at least two eigenlines of nodes intersect.

*Question* 5.7. Given a nodal integral affine surface $(B, \mathfrak{N})$, what are its "entangling points" $\mathscr{D}$ where we can entangle two (or more) nodes as in Section 4.2? What are the accumulation points of $\mathscr{D} \subset B$? Can $\mathscr{D}$ be dense in a non-empty open subset of $B$?

These seem to be very hard combinatorial questions. An accumulation point $x$ of $\mathscr{D}$ gives a fibre where there are infinitely many families of almost Hamiltonian isotopic Lagrangian knots nearby (in terms of Lagrangian flux). The simplest case is for monotone polygonal domains $\Delta$ such as in Examples 4.8 and 4.9, where $\mathscr{D}$ contains only the central point. In the simplest non-monotone example, Example 4.7, we already have an infinite subset $\mathscr{D}' \subset \mathscr{D}$, which has $\Delta_M$ as accumulation points.